\documentclass[twoside]{amsart}
\usepackage{float}
\usepackage[bookmarksnumbered,plainpages,hypertex]{hyperref}
\usepackage{graphicx}
\usepackage{amssymb}
\usepackage[usenames]{color}

\numberwithin{equation}{section}

\input{xy}
\xyoption{all}
\newtheorem{theorem}{\sc Theorem}[section]
\newtheorem{proposition}[theorem]{\sc Proposition}
\newtheorem{lemma}[theorem]{\sc Lemma}

\theoremstyle{definition}
\newtheorem{definition}[theorem]{\sc Definition}

\theoremstyle{remark}
\newtheorem{remark}[theorem]{\sc Remark}
\newtheorem{remarks}[theorem]{\sc Remarks}

\def\B{{\bf B}}
\def\A{{\bf A}}
\def\R{{\bf R}}
\def\L{{\bf L}}
\def\bH{{\bf H}}
\def\cM{\mathcal M}
\def\cB{\mathcal B}
\def\hH{\mathcal H}
\def\ot{\otimes} 
\def\Hom{\mathrm{Hom}}
\def\tw{{\rm \textsf{tw}}} 

\def\stac#1{\raise-.2cm\hbox{$\stackrel{\,\displaystyle\otimes\,}{\scriptscriptstyle{#1}}$}}
\def\sstac#1{\otimes_{#1}}

\begin{document}
\title{Hopf Algebroids}
\author{Gabriella B\"{o}hm}
\address{Research Institute for Particle and Nuclear Physics, Budapest, H-1525
  Budapest 114, P.O.B.49\\ 
\indent Hungary } 
\email{G.Bohm@rmki.kfki.hu}
\urladdr{http://www.rmki.kfki.hu/$\sim$bgabr}
\date{May 2008}
\maketitle

\tableofcontents

\pagestyle{headings}

\section{Introduction}

If it should be formulated in one sentence what a Hopf algebroid is, it should
be described as a generalisation of a Hopf algebra to a non-commutative base
algebra. More precisely, the best known examples of Hopf algebroids are Hopf
algebras. Clearly, the notion of a Hopf algebroid turns out to be a
successful generalisation only if a convincing
amount of results about Hopf algebras extend to Hopf algebroids. But one
expects more: working with Hopf algebroids should be considered to be useful
if in this way one could solve problems that could not be solved in terms of
Hopf algebras. Hopf algebroids provide us with results of both types, ones
which extend known results about Hopf algebras and also ones which are
conceptually new. 
 
Hopf algebras have been intensively studied, and successfully applied in
various 
fields of mathematics and even physics, for more than fifty years. Without
aiming at a complete list, let us mention a few applications. Hopf algebras
were used to construct invariants in topology and knot theory. 
In connection with solutions of the quantum Yang-Baxter equation, quantum
groups i.e. certain Hopf algebras play a central role. In (low
dimensional) quantum field theory Hopf algebras are capable of describing
internal symmetry of some models. In non-commutative differential geometry
(faithfully flat) Galois extensions by a Hopf algebra are interpreted as
non-commutative principal bundles.
Although the theory of Hopf algebras was (is!) extremely successful, in the
1990's there arose more and more motivations for a generalisation. 

Originally the term `Hopf algebroid' was used for cogroupoid objects in the 
category of commutative algebras. These are examples of Hopf algebroids in 
our current note, with commutative underlying algebra structure.
They found an application e.g. in algebraic topology \cite{Mor:alg.top}.
As a tool of a study of the geometry of principal fiber bundles with groupoid
symmetry, recently more general, non-commutative Hopf algebroids have been
used, but still over commutative base algebras \cite{Mrcun:hgd}. 
For some applications it is still not the necessary level of generality. 
In Poisson geometry, solutions of the dynamical Yang-Baxter equation
correspond to dynamical quantum groups,
which are not Hopf algebras \cite{EtiNik:Dynqgr}, \cite{NikVai:FinQgr},
\cite{Lu:hgd}, \cite{Xu:QGr}, \cite{DonMud:Qgr}, \cite{Kar:review}. In
topology, invariants obtained in \cite{NiTuVa:Inv} do not fit the Hopf
algebraic framework. 
In transverse geometry extensions of Hopf algebras by non-commutative base
algebras occurred \cite{ConMos:trans.geo}.
In low dimensional quantum field theories non-integral
values of the statistical (also called quantum-) dimensions in some models
exclude a Hopf algebra symmetry \cite{BohmSzl:WHAII}. Another field where
important questions could 
not be answered in the framework of Hopf algebras is non-commutative geometry,
i.e. Hopf Galois theory. Thinking about classical Galois extensions of fields
by a finite group, such an extension can be characterised without explicitly
mentioning the Galois group. A (unique, upto isomorphism) Galois group is
determined by a Galois field extension. In the case of Hopf Galois extensions
no such intrinsic characterisation, without explicit use of a Hopf algebra, is
known. Also, although the Hopf algebra describing the symmetry of a given Hopf
Galois extension is known to be non-unique, the relation between the possible
choices is not known. These questions have been handled by allowing for
non-commutative base algebras \cite{BalSzl:fin.Gal}, \cite{Kad:note}. On
the other hand, as the study of Hopf algebroids has a quite 
short past, there are many aspects of Hopf algebras that have not yet been
investigated how to extend to Hopf algebroids. It has to be admitted that
almost nothing has been done yet e.g. towards a classification and structure
theory of Hopf algebroids.  

What does it mean that the base algebra $\R$ of a Hopf algebroid is
non-commutative? Recall that a bialgebra over a commutative base ring $k$ is a
$k$-module, with compatible algebra and coalgebra structures. By analogy, in a 
bialgebroid the coalgebra structure is replaced by a coring over any not
necessarily commutative 
$k$-algebra $\R$. Also the algebra structure is replaced by a ring over a 
non-commutative base algebra. However, in order to formulate the compatibility
between  
the ring and coring structures, the base algebra of the ring has to be not
$\R$ but $\R\sstac k \R^{op}$. A Hopf algebra is a bialgebra with an
additional antipode map. In the Hopf algebroid case, the antipode relates two
different bialgebroid structures, over the base algebras $\R$ and $\R^{op}$,
respectively. 

In these notes we arrive at the notion of a Hopf algebroid after considering
all constituent structures. In Section \ref{sec:ring-coring}
$\R$-rings and $\R$-corings are introduced. They are seen to generalise
algebras 
and coalgebras, respectively. Emphasis is put on their duality. Section
\ref{sec:bgd} is devoted to a study of bialgebroids, generalising bialgebras. 
Several equivalent descriptions are given and examples are collected. In
particular, constructions of new bialgebroids from known ones are 
presented. Some of them change the base algebra of a bialgebroid, so they have
no counterparts for bialgebras. Although bialgebroid axioms are not manifestly
self-dual, duals of finitely generated and projective bialgebroids are shown
to be bialgebroids. Key properties of a bialgebroid are monoidality
of the categories of modules and comodules. This is explained in some
detail. Section \ref{sec:bgd} is closed by a most important and most
successful application, Galois theory of bialgebroids. Hopf
algebroids are the subject of Section \ref{sec:hgd}. After presenting the
definition, listing some examples and deriving some immediate consequences of
the axioms, we discuss the theory of comodules. Since in a Hopf algebroid
there are two bialgebroids (hence corings) present, comodules of the Hopf
algebroid comprise comodule structures of both. The relation between the
categories of comodules of a Hopf algebroid, and comodules of the constituent
bialgebroids, is investigated. The category of comodules of a Hopf algebroid is
proven to be monoidal, what is essential from the point of view of Galois
theory. 
Next we turn to a study of the theory of integrals. It is a good
example of results that are obtained by using some new ideas, but that extend
analogous results for Hopf algebras in a reassuring way. The structure of
Galois extensions by Hopf algebroids is investigated. Useful theorems are
presented about situations when surjectivity of a canonical map implies Galois
property. They extend known results about Hopf Galois extensions. While there
seems to be an accord in the literature that the right generalisation of a
bialgebra to a non-commutative base is a bialgebroid, there is some discussion
about the right generalisation of a Hopf algebra. We close these notes by
collecting and comparing notions suggested by various authors. 

In order to keep the list of references perspicuous, we do not refer in these
notes to papers containing classical results about Hopf algebras, which are
generalised hereby. We believe it is more useful to give here a detailed
bibliography of those papers which deal with structures over non-commutative
base. A very good and detailed bibliography of the literature of Hopf algebras
can be found e.g. in Chapter ``Hopf Algebras'' of Handbook of Algebra
\cite{CohGelWes:HabAlg}. 

{\bf Notations and conventions.} Throughout $k$ is an associative and
commutative unital ring. All algebras are associative and unital
$k$-algebras. A $k$-algebra is denoted by ${\bf A}$ and the underlying
$k$-module is denoted by $A$. On elements of $A$, multiplication is denoted by
juxtaposition. Unit element is denoted by $1_{\bf A}$. For an algebra ${\bf
A}$, with multiplication $(a,a')\mapsto aa'$, ${\bf A}^{op}$ denotes the
opposite of ${\bf A}$. As a $k$-module it is equal to $A$ and multiplication
is $(a,a')\mapsto a'a$. The category of right (resp.left) modules of an
algebra ${\bf A}$ is denoted by $\cM_{\bf A}$ (resp. ${}_{\bf A}\cM$). Hom
sets are denoted by $\mathrm{Hom}_{\bf A} (-,-)$ (resp. ${}_{\bf
A}\mathrm{Hom}(-,-)$). The category of ${\bf A}$-bimodules is denoted by
${}_{\bf A}\cM_{\bf A}$ and its hom sets by ${}_{\bf A}\mathrm{Hom}_{\bf
A}(-,-)$. Often we identify left ${\bf A}$-modules with right 
${\bf A}^{op}$-modules, but in every such case it is explicitly said. 
Action on a (say right) module $M$ of a $k$-algebra ${\bf A}$, if evaluated
on elements $m\in M$ and $a\in A$, is denoted by $\varrho_M:m\otimes a \mapsto
m \cdot a$. 

For coproducts in a coalgebra and, more generally, in a coring $C$, Sweedler's
index notation is used. That is, for an element $c\in C$, we write $c\mapsto
c_{(1)}\otimes c_{(2)}$ (or sometimes $c\mapsto c^{(1)}\otimes c^{(2)}$) for
the coproduct, where implicit summation is understood. Similarly, for a (say
right) coaction on a comodule $M$ of a coring, evaluated on an element $m\in
M$, the notation $\varrho^M:m\mapsto m_{[0]}\otimes m_{[1]}$ (or $m\mapsto
m^{[0]}\otimes m^{[1]}$) is used, where implicit summation is understood.
The category of right (resp. left) comodules of a coring $C$
is denoted by $\cM^C$ (resp. ${}^C \cM$). Hom sets are denoted by
$\mathrm{Hom}^C(-,-)$ (resp. ${}^C \mathrm{Hom}(-,-)$). 

In any category ${\mathcal A}$, the identity morphism at an object $A$ is
denoted by the same symbol $A$. Hom sets in ${\mathcal A}$ are denoted by
$\Hom_{\mathcal A}(-,-)$. 

In a monoidal category $(\cM,\otimes,U)$ we allow for non-trivial coherence 
isomorphisms $(-\otimes -)\otimes - \cong -\otimes (-\otimes -)$ and $-\otimes
U \cong - \cong U\otimes -$, but do not denote them explicitly. (Such monoidal
categories are called in the literature sometimes {\em lax} monoidal.) The
opposite of a monoidal category $(\cM,\otimes,U)$, denoted by
$(\cM,\otimes,U)^{op}$, means the same category $\cM$ with opposite monoidal
product. A functor ${\mathbb F}$ between monoidal categories $(\cM,\otimes,U)$
and $(\cM',\otimes',U')$ is said to be monoidal if there exist natural
transformations ${\mathbb F}^2:{\mathbb F}(-)\otimes' {\mathbb F}(-) \to
{\mathbb F}(-\otimes -)$ and ${\mathbb F}^0:U'\to {\mathbb F}(U)$, satisfying
usual compatibility conditions. ${\mathbb F}$ is said to be op-monoidal 
if there exist compatible natural transformations ${\mathbb F}_2:{\mathbb
F}(-\otimes -)\to {\mathbb F}(-)\otimes' {\mathbb F}(-)$ and ${\mathbb
F}_0:{\mathbb F}(U)\to U'$. A monoidal functor $({\mathbb F},{\mathbb
F}^2,{\mathbb F}^0)$ is strong monoidal if ${\mathbb F}^2$ and ${\mathbb F}^0$
are isomorphisms, and it is strict monoidal if ${\mathbb F}^2$ and ${\mathbb
F}^0$ are identity morphisms. 

\section{$\R$-rings and $\R$-corings} \label{sec:ring-coring}

A {\em monoid} in a monoidal category $(\cM,\otimes,U)$ is a triple
$(A,\mu,\eta)$. Here $A$ is an object and $\mu:A\otimes A \to A$ and
$\eta:U\to A$ are morphisms in $\cM$, satisfying associativity and unitality
conditions
\begin{equation}\label{eq:cat.associ}
\mu\circ(\mu\otimes A)=\mu\circ (A\otimes \mu)\qquad \textrm{and}\qquad 
\mu\circ(\eta \otimes A)=A=\mu\circ (A\otimes \eta).
\end{equation}
The morphism $\mu$ is called a {\em multiplication} (or {\em product}) and
$\eta$ is called a {\em unit}.
An algebra over a commutative ring $k$ can be described as a monoid in
the monoidal category $(\cM_k,\sstac k, k)$ of $k$-modules. 

A {\em right module} of a monoid $(A,\mu,\eta)$ is a pair $(V,\nu)$, where $V$
is an object and $\nu:V \ot A \to V$ is a morphism in $\cM$, such that
$$
\nu\circ(V \ot \mu)= \nu \circ (\nu \ot A)
\qquad \textrm{and}\qquad
\nu \circ (V\ot \eta) =V.
$$
Left modules are defined symmetrically.

Reversing all arrows in the definition of a monoid, we arrive at the dual
notion of a comonoid. A {\em comonoid} in a monoidal category
$(\cM,\otimes,U)$ is a triple $(C,\Delta,\epsilon)$. Here $C$ is an object and
$\Delta:C\to C\otimes C$ and $\epsilon:C\to U$ are morphisms in $\cM$,
satisfying coassociativity and counitality conditions
\begin{equation}\label{eq:cat.coassoci}
(\Delta\otimes C)\circ \Delta=(C\otimes \Delta)\circ \Delta
\qquad \textrm{and}\qquad 
(\epsilon \otimes C)\circ \Delta=C=(C\otimes \epsilon)\circ \Delta.
\end{equation}
The morphism $\Delta$ is called a {\em comultiplication} (or {\em coproduct})
and $\epsilon$ is called a {\em counit}.
A coalgebra over a commutative ring $k$ can be described as a comonoid in
the monoidal category $(\cM_k,\sstac k, k)$ of $k$-modules. 
Dualising the definition of a module of a monoid, one arrives at the notion of
a comodule of a comonoid.

Many aspects of the theory of algebras and their modules, or coalgebras and
their comodules, can be extended to monoids or comonoids in general monoidal
categories. Here we are interested in monoids and
comonoids in a monoidal category $({}_\R\cM_\R,\sstac \R, \R)$ of bimodules
over a $k$-algebra $\R$. These monoids and comonoids are called {\em
$\R$-rings} and {\em $\R$-corings}, respectively. 

\subsection{$\R$-rings} Generalising algebras over commutative rings, we study
monoids in bimodule categories. 
\begin{definition}\label{def:R-ring}
For an algebra $\R$ over a commutative ring $k$, an {\em $\R$-ring} is a triple
$(A,\mu,\eta)$. Here $A$ is an $\R$-bimodule and $\mu:A\sstac \R A \to A$ and
$\eta:R\to A$ are $\R$-bimodule maps, satisfying the associativity and unit
conditions \eqref{eq:cat.associ}. 
A {\em morphism of $\R$-rings} $f:(A,\mu,\eta)\to (A',\mu',\eta')$ is an
$\R$-bimodule map $f:A \to A'$, such that 
$f\circ \mu=\mu'\circ(f\sstac \R f)$ and $f\circ \eta=\eta'$.
\end{definition}
For an $\R$-ring $(A,\mu,\eta)$, the {\em opposite} means the $\R^{op}$-ring
$(A^{op},\mu^{op},\eta)$. Here $A^{op}$ is the same $k$-module $A$. It is
understood to be a left (resp. right) $\R^{op}$-module via the right
(resp. left) $\R$-action.
Multiplication is $\mu^{op}(a\sstac {\R^{op}} a'):=\mu(a'\sstac \R a)$ and unit
is $\eta$. 

A most handy characterisation of $\R$-rings comes from the following
observation. 
\begin{lemma}\label{lem:char.R-rings}
There is a bijective correspondence between $\R$-rings $(A,\mu,\eta)$ and
$k$-algebra homomorphisms $\eta:\R\to \A$.
\end{lemma}
Indeed, starting with an $\R$-ring $(A,\mu,\eta)$, a multiplication map
$A\sstac k A \to A$ is obtained by composing the canonical epimorphism
$A\sstac k A \to A\sstac \R A$ with $\mu$. 
Conversely, starting with an algebra map $\eta:\R\to \A$, an $\R$-bilinear
associative multiplication $A\sstac \R A \to A$ is obtained by using the
universality of the coequaliser $A\sstac k A \to A\sstac \R A$. 

An $\R$-ring $(A,\mu,\eta)$ determines monads on the categories of right and
left $\R$-modules (i.e. monoids in the monoidal categories of endofunctors on 
$\cM_\R$ and ${}_\R \cM$, respectively). They are given by $-\sstac \R
A:\cM_\R\to \cM_\R$ and $A\sstac \R -:{}_\R \cM\to {}_\R \cM$, respectively.
\begin{definition}\label{def:mod.R-ring}
A {\em right module} for an $\R$-ring $(A,\mu,\eta)$ is an algebra for the
monad $-\sstac \R A$ on the category $\cM_\R$. A {\em right module morphism} is
a morphism of algebras for the monad $-\sstac \R A$.

A {\em left module} for an $\R$-ring $(A,\mu,\eta)$ is an algebra for the monad
$A\sstac \R \, -$ on the category ${}_\R \cM$. A {\em left module morphism} is
a morphism of algebras for the monad $A\sstac \R\, -$. 
\end{definition}
Left modules of an $\R$-ring are canonically identified with right modules for
the opposite $\R^{op}$-ring.
Analogously to Lemma \ref{lem:char.R-rings}, modules for $\R$-rings can be
characterised as follows. 
\begin{lemma}\label{lem:char.mod.R-rings}
A $k$-module $M$ is a (left or right) module of an $\R$-ring $(A,\mu,\eta)$ if
and only if it is a (left or right) module of the corresponding $k$-algebra
$\A$ in Lemma \ref{lem:char.R-rings}. Furthermore, a $k$-module map $f:M\to M'$
is a morphism of (left or right) modules of an $\R$-ring $(A,\mu,\eta)$ if
and only if it is a morphism of (left or right) modules of the corresponding
$k$-algebra $\A$ in Lemma \ref{lem:char.R-rings}.
\end{lemma}

The situation when the (left or right) regular $\R$-module extends to a (left
or right) module of an $\R$-ring $(A,\mu,\eta)$ is of particular interest.
\begin{lemma}\label{lem:character}
The right regular module of a $k$-algebra $\R$ extends to a right module of an
$\R$-ring $(A,\mu,\eta)$ if and only if there exists a $k$-module map
$\chi:A\to R$, obeying the following properties.
\begin{itemize}
\item[(i)] $\chi(a\eta(r))=\chi(a)r$, for $a\in A$ and $r\in R$ (right
  $\R$-linearity),
\item[(ii)] $\chi(aa')=\chi((\eta\circ \chi)(a)a')$, for $a,a'\in A$
  (associativity),
\item[(iii)] $\chi(1_\A)=1_\R$ (unitality). 
\end{itemize}
The map $\chi$ obeying these properties is called a {\em right character} on
the $\R$-ring $(A,\mu,\eta)$.
\end{lemma}
In terms of a right character $\chi$, a right $\A$-action on $R$ is given by
$r\cdot a:=\chi(\eta(r)a)$. Conversely, in terms of a right $\A$-action on
$R$, a right character is constructed as $\chi(a):=1_\R \cdot a$.
Symmetrically, one can define a {\em left character} on an $\R$-ring
$(A,\mu,\eta)$ via the requirement that the left regular $\R$-module extends to
a left module for $(A,\mu,\eta)$.
\begin{definition}\label{def:Gal.ring}
Let $(A,\mu,\eta)$ be an $\R$-ring possessing a right character $\chi:A\to
R$. The {\em invariants} of a right module $(M,\varrho_M)$  with respect to
$\chi$ are the elements of the $k$-submodule
$$
M_\chi:=\{\  m\in M \ |\ \varrho_M(m\stac \R a)=\varrho_M(m\stac \R (\eta\circ
\chi)(a)),\quad \forall a\in A\ \}\cong \mathrm{Hom}_\A(R,M),
$$
where the isomorphism $M_\chi\to \mathrm{Hom}_\A(R,M)$ is given by $m\mapsto
(r\mapsto m\cdot \eta(r))$.
In particular, the invariants of $R$ are the elements of the subalgebra 
$$
\B:=R_\chi =\{\  b\in R \ |\ \chi(\eta(b)a)=b \chi(a),\quad \forall a\in A\ \}.
$$
Associated to a character $\chi$, there is a canonical map
\begin{equation}\label{eq:can.R-ring}
A \to {}_\B\mathrm{End}(R), \qquad 
a\mapsto \big(\ r\mapsto \chi(\eta(r)a)\ ).
\end{equation}
The $\R$-ring $(A,\mu,\eta)$ is said to be a {\em Galois} $\R$-ring (with
respect to the character $\chi$) provided that the canonical map
\eqref{eq:can.R-ring} is bijective.
\end{definition}

\subsection{$\R$-corings}\label{sec:coring}
The theory of $\R$-corings is dual to that of $\R$-rings. A detailed study
can be found in the monograph \cite{BrzWis:coring}.
\begin{definition}\label{def:R-coring}
For an algebra $\R$ over a commutative ring $k$, an {\em $\R$-coring} is a
triple $(C,\Delta,\epsilon)$. Here $C$ is an $\R$-bimodule and $\Delta:C\to
C\sstac \R C$ and $\epsilon: C\to R$ are $\R$-bimodule maps, satisfying the
coassociativity and counit conditions \eqref{eq:cat.coassoci}. 
A {\em morphism of $\R$-corings} $f:(C,\Delta,\epsilon)\to
(C',\Delta',\epsilon')$ is an $\R$-bimodule map $f:C \to C'$, such that 
$\Delta'\circ f=(f\sstac \R f)\circ \Delta$ and
$\epsilon' \circ f =\epsilon$.
\end{definition}
For an $\R$-coring $(C,\Delta,\epsilon)$, the {\em co-opposite} means the
$\R^{op}$-coring $(C_{cop},\Delta_{cop},\epsilon)$. Here $C_{cop}$ is the same
$k$-module $C$. It is understood to be a left (resp. right) $\R^{op}$-module
via the right (resp. left) $\R$-action.
Comultiplication is $\Delta_{cop}(c):=c_{(2)}\sstac {\R^{op}} c_{(1)}$
and counit is $\epsilon$. 

An $\R$-coring $(C,\Delta,\epsilon)$ determines comonads on the categories of
right and left $\R$-modules (i.e. comonoids in the monoidal categories of
endofunctors on  $\cM_\R$ and ${}_\R \cM$, respectively). They are given by
$-\sstac \R C:\cM_\R\to \cM_\R$ and $C\sstac \R -:{}_\R \cM\to {}_\R \cM$,
respectively. 
\begin{definition}\label{def:comod.R-coring}
A {\em right comodule} for an $\R$-coring $(C,\Delta,\epsilon)$ is a coalgebra
for the comonad $-\sstac \R C$ on the category $\cM_\R$. That is, a pair
$(M,\varrho^M)$, where $M$ is a right $\R$-module and $\varrho^M:M\to M\sstac
\R C$ is a right $\R$-module map satisfying the coassociativity and counit
conditions 
\begin{equation}\label{eq:comodule}
(\varrho^M\stac \R C)\circ \varrho^M=(M\stac \R \Delta)\circ \varrho^M\qquad
\textrm{and}\qquad  
(M\stac \R \epsilon)\circ \varrho^M=M.
\end{equation}
A {\em right comodule morphism} $f:(M,\varrho^M)\to (M',\varrho^{M'})$ is a
morphism of coalgebras for the comonad $-\sstac \R C$. That is, a right
$\R$-module map $f:M\to M'$, satisfying $\varrho^{M'}\circ f = (f\sstac \R C)
\circ \varrho^M$.

Symmetrically, a {\em left comodule} is a coalgebra for the comonad $C\sstac \R
\, -$ on the category ${}_\R \cM$. A {\em left comodule morphism} is a morphism
of coalgebras for the comonad $C \sstac \R\, -$. 
\end{definition}
Left comodules of an $\R$-coring are canonically identified with right
comodules for the co-opposite $\R^{op}$-coring. 

The situation when the (left or right) regular $\R$-module extends to a (left
or right) comodule of an $\R$-coring $(C,\Delta,\epsilon)$ is of particular
interest, see \cite[Lemma 5.1]{Brz:Str}. 
\begin{lemma}\label{lem:grouplike}
The (left or right) regular $\R$-module extends to a (left or right) comodule
of an $\R$-coring $(C,\Delta,\epsilon)$ if and only if there exists an element
$g\in C$, obeying the following properties.
\begin{itemize}
\item[(i)] $\Delta(g)=g\sstac \R g$. 
\item[(ii)] $\epsilon(g)=1_\R$.
\end{itemize}
The element $g$ obeying these properties is called a {\em grouplike element}
in the $\R$-coring $(C,\Delta,\epsilon)$.
\end{lemma}
Having a grouplike element $g$ in an $\R$-coring $(C,\Delta,\epsilon)$, a right
coaction on $R$ is constructed as a map $R\to C$, $r\mapsto g\cdot
r$. Conversely, a right coaction $\varrho^R:R\to C$ determines a grouplike
element $\varrho^R(1_\R)$. 
\begin{definition}\label{def:Gal.coring}
Let $(C,\Delta,\epsilon)$ be an $\R$-coring possessing a grouplike element
$g\in C$. The {\em coinvariants} of a right comodule $(M,\varrho^M)$ with
respect to $g$ are the elements of the $k$-submodule
$$
M^g:=\{\  m\in M \ |\ \varrho^M(m)=m\stac \R g\ \}\cong \mathrm{Hom}^C(R,M),
$$
where the isomorphism $M^g\to \mathrm{Hom}^C(R,M)$ is given by $m\mapsto
(r\mapsto m\cdot r)$.
In particular, the coinvariants of $R$ are the elements of the subalgebra 
$$
\B:=R^g =\{\  b\in R \ |\ b\cdot g =g\cdot b \ \}.
$$
Associated to a grouplike element $g$, there is a canonical map
\begin{equation}\label{eq:can.R-coring}
R\stac \B R \to C, \qquad r\stac \B r'\mapsto r\cdot g\cdot r'.
\end{equation}
The $\R$-coring $(C,\Delta,\epsilon)$ is said to be a {\em Galois} $\R$-coring
(with respect to the grouplike element $g$) provided that the canonical map
\eqref{eq:can.R-coring} is bijective.
\end{definition}
Let $(C,\Delta,\epsilon)$ be an $\R$-coring possessing a grouplike element
$g$. Put $\B:=R^g$.  
For any right $C$-comodule $M$, $M^g$ is a right $\B$-module. Furthermore, any
right $C$-comodule map $M\to M'$ restricts to a right $\B$-module map $M^g\to
M^{\prime g}$. There is an adjoint pair of functors
\begin{equation}\label{eq:comod.adj}
-\stac \B R: \cM_\B \to \cM^C
\qquad \textrm{and} \qquad 
(-)^g: \cM^C \to \cM_\B . 
\end{equation}
If  $(C,\Delta,\epsilon)$ is a Galois coring (with respect to $g$), then
$\cM^C$ is equivalent to the category of descent data for the extension
$\B\subseteq \R$. Hence the situation, when the functors \eqref{eq:comod.adj}
establish an equivalence, is interesting from the descent theory point of
view. 

\subsection{Duality}\label{sec:dual.ring-coring}
Beyond the formal duality between algebras and coalgebras, it is well known
that the $k$-dual of a coalgebra over a commutative ring $k$ possesses a
canonical algebra structure. The converse is true whenever a $k$-algebra is
finitely generated and projective as a $k$-module. In what follows we recall
analogues of these facts for rings and corings over an arbitrary algebra $\R$.
\begin{proposition}\label{prop:dual.ring-coring}
Let $\R$ be an algebra over a commutative ring $k$.
\begin{itemize}
\item[(1)] For an $\R$-coring $(C,\Delta,\epsilon)$, the left dual ${}^*
  C:={}_\R\mathrm{Hom}(C,R)$ possesses a canonical $\R$-ring
  structure. Multiplication is given by $(\phi\psi)(c):=\psi\big(c_{(1)} \cdot
  \phi(c_{(2)})\big)$, for $\phi,\psi\in {}^*C$ and $c\in C$.
Unit map is $\R\to {}^* C$, $r\mapsto \epsilon(-)r$.
\item[(2)] For an $\R$-ring $(A,\mu,\eta)$, which is a finitely generated and
  projective right $\R$-module, the right dual $A^*:=\mathrm{Hom}_\R(A,R)$
  possesses a canonical $\R$-coring structure. In terms of a dual basis $(\{
  a_i\in A\},\{\alpha_i\in A^*\})$, comultiplication is given by 
$\xi\mapsto \sum_i \xi(a_i -)\sstac \R \alpha_i$,
which is independent of the choice of a dual basis.
Counit is $A^*\to R$, $\xi\mapsto \xi(1_\A)$.
\item[(3)] For an $\R$-coring $(C,\Delta,\epsilon)$, which is a finitely
  generated and projective left $\R$-module, the second dual $({}^* C)^*$ is
  isomorphic to $C$ as an $\R$-coring.
\item[(4)] For an $\R$-ring $(A,\mu,\eta)$, which is a finitely generated and
  projective right $\R$-module, the second dual ${}^*(A^*)$ is isomorphic to
  $A$ as an $\R$-ring.
\end{itemize}
\end{proposition}
Applying Proposition \ref{prop:dual.ring-coring} to the co-opposite coring and 
the opposite ring, analogous correspondences are found between right duals of
corings and left duals of rings.
\begin{proposition}\label{prop:dual.mod.comod}
Let $C$ be a coring over an algebra $\R$.
\begin{itemize}
\item[(1)] Any right $C$-comodule $(M,\varrho^M)$ possesses a right module
    structure for the $\R$-ring ${}^*C$,
\begin{equation}\label{eq:dual.ac}
m\cdot \phi:= m_{[0]} \cdot \phi(m_{[1]}), \qquad \textrm{for }m\in M,\
\phi\in {}^*  C.
\end{equation}
Any right $C$-comodule map becomes a ${}^* C$-module map. That is, there is a
  faithful functor $\cM^C\to \cM_{\,{}^*   \!C}$.
\item[(2)] The functor $\cM^C\to \cM_{\,{}^*\! C}$ is an equivalence if and
    only if  $C$ is a finitely generated and projective left $\R$-module.
\end{itemize}
\end{proposition}
There is a duality between Galois rings and Galois corings too.
\begin{proposition}\label{prop:dual.Galois.ring-coring}
Let $C$ be an $\R$-coring which is a finitely generated and
projective left $\R$-module. For an element $g\in C$, introduce the map
$\chi_g:{}^*C \to R$, $\phi \mapsto \phi(g)$.
The following statements hold.
\begin{itemize}
\item[(1)] The element $g\in C$ is grouplike if and only if $\chi_g$ is a
  right character on the $\R$-ring ${}^* C$.
\item[(2)] An element $b\in R$ is a coinvariant of the right $C$-comodule $R$
  (with coaction induced by a grouplike element $g$) if and only if $b$ is an
  invariant of the right ${}^*C$-module $R$ (with respect to the right
  character $\chi_g$). 
\item[(3)] The $\R$-coring $C$ is a Galois coring (with respect to a
  grouplike element $g$) if and only if the $\R$-ring ${}^* C$ is a Galois
  ring (with respect to the right character $\chi_g$).
\end{itemize}
\end{proposition}

\section{Bialgebroids}\label{sec:bgd}
In Section \ref{sec:ring-coring} we generalised algebras and coalgebras over
commutative rings to
monoids and comonoids in bimodule categories. We could easily do so, the
category of bimodules over any $k$-algebra $\R$ is a monoidal category, just as
the category of $k$-modules. If we try to generalise bialgebras to a
non-commutative base algebra $\R$, however, we encounter difficulties. Recall
that a $k$-bialgebra consists of an algebra $(B,\mu,\eta)$, and a coalgebra
$(B,\Delta,\epsilon)$ defined on the same $k$-module $B$. They are subject to
compatibility conditions. Unit and multiplication must be coalgebra maps
or, equivalently, counit and comultiplication must be algebra maps. This means
in particular that, for any elements $b$ and $b'$ in $B$, multiplication and
comultiplication satisfy the condition
\begin{equation}\label{eq:ba.comp}
(bb')_{(1)}\stac k (bb')_{(2)}= b_{(1)} b'_{(1)}\stac k b_{(2)} b'_{(2)}.
\end{equation}
Note that \eqref{eq:ba.comp} is formulated in terms of the symmetry
$\tw$ in $\cM_k$. For any $k$-modules $M$ and $N$, the twist map
$\tw_{M,N}:M\sstac k N\to N\sstac k M$ maps $m\sstac k n$ to $ n\sstac
k m$. Precisely, \eqref{eq:ba.comp} is equivalent to 
$$
\Delta\circ \mu = (\mu\stac k \mu)\circ (B\stac k \tw_{B,B} \stac k
B)\circ (\Delta\stac k \Delta).
$$
In the literature one can find generalisations when $\tw$ is replaced
by a braiding \cite{Tak:braid}, \cite{Maj:braid}, \cite{Scha:braid}. (For an
approach when $\tw$ is replaced by a {\em mixed 
distributive law}, see \cite{MesWis:Hopf.mon}.) However, general bimodule
categories are neither symmetric nor 
braided. There is no natural way to formulate an analogue of
\eqref{eq:ba.comp} in a bimodule category. In fact, more sophisticated ideas
are needed.

The notion which is known today as a (left) bialgebroid, was introduced
(independently) by several authors. The first definition is due to Takeuchi,
who used the name {\em $\times_R$-bialgebra} in \cite{Tak:bgd}. Some twenty
years later, with motivation coming from Poisson geometry, in \cite{Lu:hgd} Lu
proposed an equivalent definition. The term {\em bialgebroid} is due to her. A
third equivalent set of axioms was invented by Xu in \cite{Xu:QGr}. The
equivalence of the listed definitions is far from obvious. It was proven by
Brzezi\'nski and Militaru in \cite{BrzeMil:bgd}. 
Symmetrical notions of left and right bialgebroids were formulated and studied
by Kadison and Szlach\'anyi in \cite{KadSzl:D2bgd}. 
The definition presented here is a slightly reformulated version of the one in
\cite{Lu:hgd} or \cite{KadSzl:D2bgd}.

\subsection{Right and left bialgebroids}\label{sec:bgd.def}
In contrast to the definition of a bialgebra, in this section a
bialgebroid is not described as a compatible monoid and comonoid in some
monoidal category (of bimodules). The ring and coring structures of a
bialgebroid are defined over different base algebras: they are a monoid and a
comonoid in different monoidal categories. Recall from Lemma 
\ref{lem:char.R-rings} that an $\R\sstac k \R^{op}$-ring $A$ (for some algebra
$\R$ over a commutative ring $k$) is described by a $k$-algebra map
$\eta:\R\sstac k \R^{op}\to \A$. Equivalently, instead of $\eta$, we can
consider its restrictions 
\begin{equation}\label{eq:s-t}
s:=\eta(-\stac k 1_\R):\R\to \A\qquad \textrm{and}\qquad 
t:=\eta(1_\R\stac k\, -):\R^{op}\to \A,
\end{equation}
which are $k$-algebra maps with commuting ranges in $\A$. The maps $s$ and $t$
in \eqref{eq:s-t} are called the {\em source} and {\em target} maps of an
$\R\sstac k \R^{op}$-ring $A$, respectively. In what follows, an $\R\sstac k
\R^{op}$-ring will be given by a triple $(\A,s,t)$, where $\A$ is a
$k$-algebra (with underlying $k$-module $A$) and $s$ and $t$ are algebra maps
with commuting ranges as in \eqref{eq:s-t}. 
\begin{definition}\label{def:right.bgd}
Let $\R$ be an algebra over a commutative ring $k$.
A {\em right $\R$- bialgebroid} $\cB$ consists of an $\R\sstac k \R^{op}$-ring
$(\B,s,t)$ and an $\R$-coring $(B,\Delta,\epsilon)$ on the same $k$-module
$B$. They are subject to the following compatibility axioms. 
\begin{itemize}
\item[(i)] The bimodule structure in the  $\R$-coring $(B,\Delta,\epsilon)$ is
  related to the $\R\sstac k \R^{op}$-ring $(\B,s,t)$ via
\begin{equation}\label{eq:rbgd.bimod}
r\cdot b\cdot r':= b s(r') t(r),\quad \textrm{for }r,r'\in R,\ b\in B.
\end{equation}
\item[(ii)] Considering $B$ as an $\R$-bimodule as in \eqref{eq:rbgd.bimod},
  the coproduct $\Delta$ corestricts to a $k$-algebra map from $\B$ to 
\begin{equation}\label{eq:Tak.prod}
B\times_R B :=\{\ \sum_i b_i\stac \R b'_i\ |\ \sum_i s(r) b_i \stac \R b'_i= 
\sum_i b_i \stac \R t(r) b'_i,\quad \forall r\in R\ \},
\end{equation}
where $B\times_R B$ is an algebra via factorwise multiplication.
\item[(iii)] The counit $\epsilon$ is a right character on the $\R$-ring
  $(\B,s)$. 
\end{itemize}
\end{definition}
\begin{remarks}
The bialgebroid axioms in Definition \ref{def:right.bgd} have some immediate
consequences.
\begin{itemize}
\item[(1)] Note that the $k$-submodule $B\times_R B$ of $B\sstac \R B$ is
  defined in such a way that factorwise multiplication is well defined on
  it. $B\times_R B$ is called the {\em Takeuchi product}. In fact, it has more
  structure than that of a $k$-algebra: it is 
  an $\R\sstac k \R^{op}$-ring with unit map $\R\sstac k \R^{op}\to B\times_R
  B$, $r\sstac k r'\mapsto t(r')\sstac \R s(r)$. The (corestriction of the)
  coproduct is easily checked to be an $\R\sstac k \R^{op}$-bimodule map $B\to
  B\times_R B$. 
\item[(2)] Axiom (iii) is equivalent to the requirement that the counit
  $\epsilon$ is a right character on the $\R^{op}$-ring $(\B,t)$.
\item[(3)] Yet another equivalent formulation of axiom (iii) is the
  following. The map 
$\theta: \B \to \mathrm{End}_k(R)^{op}$, $b\mapsto \big(\ r\mapsto
\epsilon(s(r)b)\ \big)$
is a $k$-algebra map, where $\mathrm{End}_k(R)^{op}$ is an algebra via
opposite composition of endomorphisms. The map $\theta$ is called an {\em
  anchor map} in \cite{Xu:QGr}. 
\end{itemize}
\end{remarks}
Recall that in a bialgebra over a commutative ring, replacing the algebra with
the opposite one, or replacing the coalgebra with the co-opposite one, one
arrives at bialgebras again. Analogously, the {\em co-opposite} of a right
$\R$-bialgebroid $\cB$ with structure maps denoted as in Definition
\ref{def:right.bgd}, 
is the following right $\R^{op}$-bialgebroid $\cB_{cop}$. $\R^{op}\sstac k
\R$-ring structure is $(\B,t,s)$, $\R^{op}$-coring structure is
$(B_{cop},\Delta_{cop}, \epsilon)$. However, the 
$\R\sstac k \R^{op}$-ring $(\B^{op},t,s)$ and the $\R$-coring
$(B,\Delta,\epsilon)$ do not satisfy the same axioms in Definition
\ref{def:right.bgd}. Instead, they are subject to a symmetrical version
of Definition \ref{def:right.bgd}.
\begin{definition}\label{def:left.bgd}
Let $\R$ be an algebra over a commutative ring $k$.
A {\em left $\R$- bialgebroid} $\cB$ consists of an $\R\sstac k \R^{op}$-ring
$(\B,s,t)$ and an $\R$-coring $(B,\Delta,\epsilon)$ on the same $k$-module
$B$. They are subject to the following compatibility axioms. 
\begin{itemize}
\item[(i)] The bimodule structure in the  $\R$-coring $(B,\Delta,\epsilon)$ is
  related to the $\R\sstac k \R^{op}$-ring $(\B,s,t)$ via
\begin{equation}\label{eq:lbgd.bimod}
r\cdot b\cdot r':= s(r) t(r') b,\quad \textrm{for }r,r'\in R,\ b\in B.
\end{equation}
\item[(ii)] Considering $B$ as an $\R$-bimodule as in \eqref{eq:lbgd.bimod},
  the coproduct $\Delta$ corestricts to a $k$-algebra map from $\B$ to 
\begin{equation}\label{eq:l.Tak.prod}
B\, {}_R\!\!\times B :=\{\ \sum_i b_i\stac \R b'_i\ |\ \sum_i b_i t(r) \stac \R
b'_i=  
\sum_i b_i \stac \R b'_i s(r),\quad \forall r\in R\ \},
\end{equation}
where $B\, {}_R\!\!\times B$ is an algebra via factorwise multiplication.
\item[(iii)] The counit $\epsilon$ is a left character on the $\R$-ring
  $(\B,s)$. 
\end{itemize}
\end{definition}
Since in this note left and right bialgebroids are considered
simultaneously, we use two versions of Sweedler's index notation. In a
left bialgebroid we use lower indices to denote components of the coproduct,
i.e. we write $b\mapsto b_{(1)}\sstac \R b_{(2)}$. In a right bialgebroid we
use upper indices to denote components of the coproduct, i.e. we write
$b\mapsto b^{(1)}\sstac \R b^{(2)}$. In both cases implicit summation is
understood. 

Recall that coalgebras over a commutative ring $k$ form a monoidal category
with respect to the $k$-module tensor product. Bialgebras over $k$ can be
described as monoids in the monoidal category of $k$-coalgebras. 
In \cite{Tak:bgd} Takeuchi defined bialgebroids ({\em $\times_R$-bialgebras}
in his terminology) as monoids in a monoidal category of certain corings, too.
By \cite[Definition 3.5]{Tak:sqrMor}, for two $k$-algebras $\R$
and ${\bf S}$, an {\em ${\bf S}|\R$-coring} is an ${\bf S}\sstac k
\R$-bimodule $C$, together with an $\R$-coring structure
$(C,\Delta,\epsilon)$, such that the following identities hold.
\begin{eqnarray*}
&&\Delta((s\stac k 1_\R)\cdot c\cdot (s'\stac k 1_\R))=
c^{(1)}\cdot (s'\stac k 1_\R)\ \stac \R\ (s\stac k 1_\R)\cdot c^{(2)}\qquad
\textrm{and} 
\\
&&(s\stac k 1_\R)\cdot c^{(1)}\ \stac \R\  c^{(2)}=
c^{(1)}\ \stac \R\  c^{(2)}\cdot (s\stac k 1_\R),
\qquad\qquad\qquad\qquad \textrm{for }s,s'\in S, c\in C.
\end{eqnarray*}
Morphisms of ${\bf S}|\R$-corings are morphisms of $\R$-corings which are in
addition ${\bf S}$-bimodule maps.
In particular, it can be shown by using the same methods
as in \cite{BrzeMil:bgd}, that the notion of an $\R|\R$-coring is equivalent
to a {\em  $\times_R$-coalgebra} in \cite[Definition 4.1]{Tak:bgd}. Part (1)
of following Theorem \ref{thm:R|R.corings} is thus a reformulation of
\cite[Proposition 4.7]{Tak:bgd}. Part (2) states an equivalence of a right
$\R$-bialgebroid in Definition \ref{def:right.bgd}, and a symmetrical version
of a $\times_R$-bialgebra in \cite[Definition 4.5]{Tak:bgd}.
\begin{theorem}\label{thm:R|R.corings}
For an algebra $\R$ over a commutative ring $k$, the following statements hold.
\begin{itemize}
\item[(1)] $\R|\R$-corings form a monoidal category. Monoidal product of two
  objects $(C,\Delta,\epsilon)$ and $(C',\Delta',\epsilon')$ is the $\R\sstac
  k \R^{op}$-module tensor product
$$
C\odot C':= C\stac k C'/
\{\ (1_\R\stac k r)\cdot c\cdot (1_\R\stac k r')\stac k c' - 
c\stac k (r'\stac k 1_\R)\cdot c'\cdot (r\stac k 1_\R)\ |\ r,r'\in R\ \}.
$$
$C\odot C'$ is an $\R\sstac k\R$-bimodule, via the actions
$$
(r_1\stac k r_2)\cdot (c\odot c')\cdot (r_1'\stac k r_2'):= 
(r_1\stac k 1_\R)\cdot c\cdot (r'_1\stac k 1_\R)\odot (1_\R \stac k r_2)\cdot
c'\cdot (1_\R \stac k r_2'). 
$$
Coproduct and counit in $C\odot C'$ are 
$$
c\odot c'\mapsto (c^{(1)}\odot c^{\prime(1)})\stac \R 
(c^{(2)}\odot c^{\prime(2)})\qquad \textrm{and}\qquad 
c\odot c'\mapsto \epsilon'\big((\epsilon(c)\stac k 1_\R)\cdot c'\big).
$$
Monoidal unit is $R\sstac k R$, with $\R|\R$-coring structure described in
Section \ref{sec:RotR-op} below.
Monoidal product of morphisms $\alpha_i:(C_i,\Delta_i,\epsilon_i)\to 
(C'_i,\Delta'_i,\epsilon'_i)$, for $i=1,2$,
is given by 
$$
(\alpha_1\odot \alpha_2)(c_1\odot c_2):=\alpha_1(c_1)\odot \alpha_2(c_2).
$$
\item[(2)] Monoids in the monoidal category of $\R|\R$-corings are the same as
  right $\R$-bialgebroids.
\item[(3)] Monoidal morphisms in the monoidal category of $\R|\R$-corings are
  the same as maps of $\R\sstac k \R^{op}$-rings as well as of $\R$-corings.
\end{itemize}
\end{theorem}
Considering a right $\R$-bialgebroid as an $\R|\R$-coring, the
  $\cdot|\R$ -bimodule structure is given by right multiplications by the
  source and target maps. The $\R|\cdot$ -bimodule structure is given by left
  multiplications by the source and target maps.

Theorem \ref{thm:R|R.corings} was extended by Szlach\'anyi in
  \cite{Szl:MonMor}. He has shown that ${\bf S}|\R$-corings form a bicategory,
  monads in which are the same as right bialgebroids. This makes it possible
  to define ({\em base changing}) morphisms of bialgebroids as bimodules for
  the corresponding monads. 
  The constructions described in Section \ref{sec:Dri.twist}, 
and also in Definition \ref{def:R-matrix}, provide
  examples of bialgebroid morphisms in this sense. For more details we refer
  to  
\cite{Szl:MonMor}.  

An ${\bf S}|\R$-coring $C$ can be looked at as an ${\bf S}\sstac k {\bf
  S}^{op}$-$\R\sstac k \R^{op}$ bimodule in a canonical way. Identifying
  ${\bf S}$-bimodules with 
  right ${\bf S}\sstac k {\bf S}^{op}$-modules, and $\R$-bimodules with right
  $\R\sstac k \R^{op}$-modules, there is an op-monoidal left adjoint functor
  $- \sstac {{\bf S}\sstac k {\bf S}^{op}} C: {}_{\bf S}\cM_{\bf S}\to {}_\R
  \cM_\R$. As a matter of fact, this correspondence establishes a
  bifunctor from the bicategory of ${\bf S}|\R$-corings to the
  2-category of op-monoidal left adjoint functors between bimodule
  categories. For any two $0$-cells (i.e. algebras ${\bf S}$ and $\R$) it
  gives an equivalence of the vertical subcategories. So, in addition to a 
  characterisation of bialgebroids as monads in the bicategory of ${\bf
  S}|\R$-corings, they can be described as monads in the bicategory of
  op-monoidal left adjoint functors between bimodule categories. Monads in the
  2-category of op-monoidal functors, (i.e. op-monoidal monads such that
  multiplication and unit natural transformations of the monad are compatible
  with the op-monoidal structure) were termed {\em Hopf monads} in
  \cite{Moe:monad.tens.cat} and {\em bimonads} in
  \cite{Szlach:mon.Eil.Moore}. Using the latter terminology, the following
  characterisation of 
  bialgebroids in \cite[Theorem  4.5]{Szlach:mon.Eil.Moore} is obtained.
\begin{theorem}\label{thm:bimonad}
For an algebra $\R$, any right $\R$-bialgebroid 
induces a bimonad on ${}_\R\cM_\R$ which possesses a right adjoint. Conversely,
every bimonad on ${}_\R\cM_\R$ possessing a right adjoint is naturally
equivalent to a bimonad induced by a right $\R$-bialgebroid.
\end{theorem}
Another aspect of the equivalence in Theorem \ref{thm:bimonad} is
explained in Section \ref{sec:bgd.mod}.  

In the paper \cite{DayStr:Qcat} by Day and Street, op-monoidal monads were
studied in the more general framework of pseudo-monoids in monoidal
bicategories. Based on Theorem \ref{thm:bimonad}, a description of
bialgebroids as {\em strong monoidal morphisms between pseudo-monoids} in the
monoidal bicategory of bimodules was obtained \cite[Proposition
  3.3]{DayStr:Qcat}.  

\subsection{Examples}\label{sec:ex.bgd}
In order to make the reader more familiar with the notion, in this section we
list some examples of bialgebroids. 

\subsubsection{Bialgebras}
Obviously, a bialgebra $B$ over a commutative ring $k$ determines a
(left or right) $k$-bialgebroid in which both the source and target maps are
equal to the unit map $k\to \B$. Note, however, that there are
$k$-bialgebroids in which the source and target maps are different, or
possibly they are equal but their range is not central in the total algebra,
hence they are not bialgebras, see e.g. Section \ref{sec:q.torus}.

\subsubsection{Weak bialgebras}\label{sec:WBA}
A weak bialgebra over a commutative ring $k$ consists of an algebra and 
a coalgebra structure on the same $k$-module $B$, subject to compatibility
axioms generalising the axioms of a bialgebra \cite{Nill:WBA},
\cite{BNSz:WHAI}. Explicitly, the coproduct $\Delta$ is required to be
multiplicative in the sense of \eqref{eq:ba.comp}. Unitality of the coproduct
$\Delta$ and multiplicativity of the counit $\epsilon$ are weakened to the
identities 
\begin{eqnarray*}
&&(\Delta(1_\B)\stac k 1_\B)(1_\B \stac k \Delta(1_\B))=(\Delta\stac k B)\circ
  \Delta(1_\B)= 
(1_\B \stac k \Delta(1_\B))(\Delta(1_\B)\stac k 1_\B), \qquad \textrm{and}\\
&&\epsilon(b1_{(1)}) \epsilon(1_{(2)}b')=\epsilon (bb') =
\epsilon(b1_{(2)}) \epsilon(1_{(1)}b'),\qquad \qquad\qquad \qquad\qquad \qquad
\quad \ \ \textrm{for }b,b'\in B, 
\end{eqnarray*}
respectively. Here $1_{(1)} \sstac k 1_{(2)}$ denotes $\Delta(1_\B)$ (which
may differ from $1_\B \sstac k 1_\B$). The map
$$
\sqcap^R:B \to B,\qquad b\mapsto 1_{(1)} \epsilon(b1_{(2)})
$$
is checked to be idempotent. Its range is a subalgebra $\R$ of $\B$. $B$ is an
$\R\sstac k \R^{op}$-ring, with source map given by the inclusion $\R\to \B$
and target map given by the restriction to $R$ of the map 
$B \to B$, $b\mapsto \epsilon(b1_{(1)}) 1_{(2)}$.
Consider $B$ as an $\R$-bimodule via right multiplication by these source
and target maps. 
Coproduct in an $\R$-coring $B$ is obtained by composing $\Delta:B \to B\sstac
k B$ with the canonical epimorphism $B \sstac k B\to B \sstac \R B$. It has a
counit $\sqcap^R$. The $\R\sstac k \R^{op}$-ring and $\R$-coring structures
constructed on $B$ in this way constitute a right $\R$-bialgebroid. 

A left $\R^{op}$-bialgebroid structure in a weak bialgebra $B$ is constructed
symmetrically. Its source map is the inclusion map into $B$ of the range
subalgebra of the idempotent map
$$
\sqcap^L:B\to B,\qquad b\mapsto \epsilon(1_{(1)} b) 1_{(2)}.
$$
Coproduct is obtained by composing the weak coproduct $\Delta:B\to B\sstac k
B$ with an appropriate canonical epimorphism to an $\R^{op}$-module tensor
product, too. 

As it was observed by Schauenburg in \cite{Scha:WHA&bgd} (see also
\cite[Sections 1.2 and 1.3]{Szlach:FinQG}),
the base algebra $\R$ of a weak bialgebra $B$ is {\em Frobenius
  separable}. This means the existence of a $k$-module map $R\to k$ (given
  by the counit $\epsilon$), possessing a dual basis $\sum_i e_i\sstac k f_i
  \in R\sstac k R$, such that $\sum_i e_i f_i = 1_\R$. The dual basis property
  means
$
\sum_i e_i \epsilon(f_i r)=r=\sum_i\ \epsilon(r e_i) f_i$, for all $r\in R
$.
In a weak bialgebra a dual basis is given by $1_{(1)}\sstac k
\sqcap^R(1_{(2)})\in R\sstac k R$.  
In \cite{Scha:WHA&bgd}, \cite{KadSzl:D2bgd} and \cite{Szlach:FinQG} also the
converse is proven: a $k$-module map  
$\epsilon:R\to k$ on the base algebra $\R$ of a bialgebroid $\cB$, with
normalised dual basis $\sum_i e_i\sstac k f_i \in R\sstac k R$, determines a
weak bialgebra structure on the underlying $k$-module $B$.
In \cite[Proposition 9.3]{KadSzl:D2bgd} any separable algebra over a field was
proven to be Frobenius separable.

Consider a small category ${\mathcal C}$ with finitely many objects.
For a commutative ring $k$, the free $k$-module generated by the morphisms in
${\mathcal C}$ carries a weak bialgebra
structure. The product of two morphisms is equal to the composite morphism if
they are composable, and zero otherwise. Unit element is a (finite) sum
of the identity morphisms for all objects. Extending the product $k$-linearly,
we obtain a $k$-algebra. Coproduct is diagonal on all morphisms,
i.e. $\Delta(f)=f\sstac k f$. Counit maps every morphism to $1_k$. Extending
the coproduct and the counit $k$-linearly, we obtain a $k$-coalgebra. The
algebra and coalgebra structures constructed in this way constitute a weak
bialgebra. 

\subsubsection{The bialgebroid $\R\otimes \R^{op}$}\label{sec:RotR-op}
For any algebra $\R$ over a commutative ring $k$, a simplest possible
right $\R$-bialgebroid is constructed on the algebra $\R\sstac k
\R^{op}$. Source and target maps are given by the inclusions
$$
\R\to \R\stac k \R^{op},\quad r\mapsto r\stac k 1_\R\qquad \textrm{and}\qquad 
\R^{op}\to \R\stac k \R^{op},\quad r\mapsto 1_\R\stac k r,
$$
respectively. Coproduct is
$$
\R\stac k \R^{op}\to (\R\stac k \R^{op})\stac \R (\R\stac k \R^{op}),\qquad
r\stac k r'\mapsto (1_\R\stac k r')\stac \R (r\stac k 1_\R).
$$
Counit is 
$\R\sstac k \R^{op}\to R$, $r\sstac k r'\mapsto r'r$. The corresponding
$\R|\R$-coring occurred in Theorem \ref{thm:R|R.corings} (1).
The opposite co-opposite
of the above construction yields a left $\R^{op}$-bialgebroid structure on
$\R^{op}\sstac k \R\cong \R\sstac k \R^{op}$. 

\subsection{Duality}\label{sec:dual.bgd}
In Section \ref{sec:dual.ring-coring} the duality between $\R$-rings and
$\R$-corings has been studied. Now we shall see how it leads to a duality of
bialgebroids. Recall that the axioms of a bialgebra over a commutative ring
$k$ are self-dual. That is, the diagrams (in the category $\cM_k$ of
$k$-modules), expressing the bialgebra axioms, remain unchanged if all
arrows are reversed. As a consequence, if a bialgebra $B$ is finitely
generated and projective as a $k$-module (hence possesses a dual in $\cM_k$),
then the dual has a bialgebra structure too, which is the transpose of the
bialgebra structure of $B$. In contrast to bialgebras, axioms of a bialgebroid
are not self dual in the same sense. Although it follows by the considerations
in 
Section \ref{sec:dual.ring-coring} that the $\R$-dual of a finitely generated
and projective bialgebroid possesses an $\R$-ring, and an $\R$-coring
structure, it is not obvious that these dual structures constitute a
bialgebroid. The fact that they do indeed, was shown first by Schauenburg in 
\cite{Scha:dual.double.bgd}. A detailed study can be found also in the paper
\cite{KadSzl:D2bgd} by Kadison and Szlach\'anyi. Our presentation here is
closer to \cite[Proposition 2.5]{KadSzl:D2bgd}.
\begin{proposition}\label{prop:dual.bgd}
Let $\R$ be an algebra. Consider a left bialgebroid $\cB$ over $\R$, 
which is a finitely generated and projective left $\R$-module (via left
multiplication by the source map). Then the left $\R$-dual ${}^* B:={}_\R
\mathrm{Hom}(B,R)$ possesses a (canonical) right $\R$-bialgebroid structure.
\end{proposition}
Applying part (1) of Proposition \ref{prop:dual.ring-coring} to the
$\R$-coring $(B,\Delta,\epsilon)$ underlying $\cB$, we conclude on the
existence of an $\R$-ring structure on ${}^* B$. Its unit map is
${}^* s: \R\to {}^* B$, $r\mapsto \epsilon(-)r$.
Multiplication is given by 
$$
(\beta \beta')(b)=\beta'\big(\ t(\beta(b_{(2)}))b_{(1)}\ \big),\qquad
\textrm{for }\beta,\beta'\in B,
$$
where we used that the right $\R$-module structure in $B$ is given via the
target map $t$ in $\cB$. 
Applying (a symmetrical version of) part (2) of Proposition
\ref{prop:dual.ring-coring} to the $\R$-ring 
$(\B,s)$ underlying $\cB$, we conclude on the existence of an $\R$-coring
structure in ${}^* B$. It has a bimodule structure 
$$
r\cdot \beta\cdot r'= \beta\big( - s(r) \big) r'\qquad \textrm{for }r,r'\in
R,\ \beta\in {}^* B.
$$
In particular, $\epsilon\cdot r ={}^*s(r)$, for $r\in R$, as expected. A
to-be-target-map is defined as ${}^* t(r):=r\cdot \epsilon$, for $r\in R$.
Counit in the $\R$-coring ${}^* B$ is ${}^* \epsilon:{}^* B\to R$, $\beta
\mapsto \beta(1_\B)$. Coproduct is given in terms of a dual basis $(\{\ a_i\in
B\ \},\{\ \alpha_i\in {}^* B \ \})$ as 
$$
{}^* \Delta:{}^* B\to {}^* B \stac \R {}^* B,\qquad 
\beta\mapsto \sum_i \alpha_i \stac \R \beta (- a_i).
$$
Right bialgebroid axioms are verified by direct computations.

One can apply Proposition \ref{prop:dual.bgd} to the co-opposite left
bialgebroid, which is a finitely generated and projective right
$\R^{op}$-module via left multiplication by the target map.
In this way one verifies that the right dual $B^*:=\mathrm{Hom}_\R(B,R)$
of a left $\R$-bialgebroid $\cB$, which is a finitely generated and projective
right $\R$-module, possesses a right $\R$-bialgebroid structure
$\cB^* := ({}^*(\cB_{cop}))_{cop}$. Note that (conventionally), multiplication
is chosen in such a way which results in right bialgebroid structures on both
duals of a left bialgebroid. Applying the constructions to the opposite
bialgebroid, left and right duals of a right bialgebroid, which is a finitely
generated and projective $\R$-module on the appropriate side, are concluded to
be left bialgebroids. Our convention is to choose $(\cB^{op})^*:=(\cB^*)^{op}$
and  
${}^*(\cB^{op}):= ({}^* \cB)^{op}$.
 
\subsection{Construction of new bialgebroids from known ones}
In addition to the examples in Section \ref{sec:ex.bgd},
further (somewhat implicit) examples of bialgebroids are provided by various
constructions starting with given bialgebroids.

\subsubsection{Drinfel'd twist}\label{sec:Dri.twist}
A Drinfel'd twist of a bialgebra $B$ over a commutative ring $k$ is a
bialgebra with the same algebra structure in $B$, and coproduct deformed (or
{\em  twisted}) by an invertible normalised 2-cocycle in $B$ (the so called
{\em Drinfel'd element}). In this section we recall analogous Drinfel'd twists
of bialgebroids from \cite[Section 6.3]{Szl:MonMor}. More general twists,
which do not correspond to invertible Drinfel'd elements, are studied in
\cite{Xu:QGr}. Such generalised twists will not be considered here.
\begin{definition}
For an algebra $\R$, consider a right $\R$-bialgebroid $\cB$, 
with structure maps denoted as in Definition \ref{def:right.bgd}.
An (invertible) element $J$ of the Takeuchi product $B\times_R B$ is called an
{\em (invertible) normalised 2-cocycle in $\cB$} provided it satisfies the
following conditions. 
\begin{itemize}
\item[(i)]
$(t(r)\sstac \R s(r'))J=J(t(r)\sstac \R s(r'))$, for $r,r'\in R$,
  ({\em bilinearity}), 
\item[(ii)] $(J\sstac \R 1_\B)(\Delta\sstac \R B)(J)=(1_B\sstac \R J)(B\sstac
  \R \Delta)(J)$ ({\em cocycle condition}), 
\item[(iii)] $(\epsilon \sstac \R B)(J)=1_\B=(B\sstac \R \epsilon)(J)$
  ({\em normalisation}). 
\end{itemize}
\end{definition}
\begin{proposition}
Let $J$ be an invertible normalised 2-cocycle in a right $\R$-bialgebroid
$\cB$. The $\R\sstac k \R^{op}$-ring $(\B,s,t)$ in $\cB$, the counit
$\epsilon$ of $\cB$ and the {\em twisted form} $\Delta_J:= J\Delta(-)
J^{-1}$ of the coproduct $\Delta$ in $\cB$, constitute a right $\R$-bialgebroid
$\cB_J$. 
\end{proposition}

\subsubsection{Cocycle double twist}\label{sec:dual.twist}
Dually to the construction in Section \ref{sec:Dri.twist}, one can leave the
coproduct in a bialgebroid $\cB$ unchanged and {\em twist} multiplication by an
invertible normalised 2-cocycle {\em on} $\cB$. For an algebra $\R$ over a
commutative ring $k$, consider a left $\R$-bialgebroid $\cB$, with structure
maps denoted as in Definition \ref{def:left.bgd}. Recall from Theorem
\ref{thm:R|R.corings} that the $\R\sstac k \R^{op}$-module tensor product
$B\sstac {\R\sstac k \R^{op}} B$ (with respect to the right (resp. left)
actions given by right (resp. left) multiplications by $s$ and $t$) is an
$R$-coring. It has a bimodule structure $r\cdot(b\otimes b')\cdot r':=s(r)
t(r')b\otimes b'$, coproduct $b\otimes b'\mapsto (b_{(1)}\otimes
b'_{(1)})\sstac \R (b_{(2)}\otimes b'_{(2)})$ and counit $b\otimes b'\mapsto
\epsilon (bb')$. Hence there is a corresponding convolution algebra ${}_\R
\mathrm{Hom}_\R(B\sstac {\R\sstac k \R^{op}} B, R)$ with multiplication
$(f\diamond g)(b\otimes b'):= f(b_{(1)}\otimes b'_{(1)}) g(b_{(2)}\otimes
b'_{(2)})$. 
\begin{definition}
Let $\R$ be an algebra over a commutative ring $k$ and 
let $\cB$ be a left $\R$-bialgebroid, with 
structure maps denoted as in Definition \ref{def:left.bgd}.
An (invertible) element of the convolution algebra ${}_\R
\mathrm{Hom}_\R(B\sstac {\R\sstac k \R^{op}} B, R)$ 
is called an {\em (invertible) normalised 2-cocycle on $\cB$} provided it
  satisfies the following conditions, for 
  $b,b',b^{\prime\prime}\in B$ and $r,r'\in R$. 
\begin{itemize}
\item[(i)] $\sigma(s(r)t(r')b,b')=r\sigma(b,b')r'$ ({\em bilinearity}),
\item[(ii)] $\sigma(b,s(\sigma({b'}_{(1)},{b^{\prime\prime}}_{(1)}))
  {b'}_{(2)} {b^{\prime\prime}}_{(2)})
  = \sigma(s(\sigma(b_{(1)},{b'}_{(1)}))b_{(2)}{b'}_{(2)}, b^{\prime\prime})$
  ({\em cocycle condition}), 
\item[(iii)]$\sigma(1_\B,b)=\epsilon(b)=\sigma(b,1_\B)$ ({\em normalisation}).
\end{itemize}
\end{definition}
\begin{proposition}
Let $\sigma$ be an invertible normalised 2-cocycle on a left $\R$-bialgebroid
$\cB$, with inverse ${\widetilde \sigma}$. 
The source and target maps $s$ and $t$ in $\cB$, the $\R$-coring
$(B,\Delta,\epsilon)$ in $\cB$, and the {\em twisted product} 
$
b\cdot_\sigma b':= 
s\big(\sigma(b_{(1)}, b'_{(1)})\big)
t\big({\widetilde \sigma}(b_{(3)}, b'_{(3)})\big) b_{(2)} b'_{(2)}$, 
for $b,b'\in B$,
constitute a left $\R$-bialgebroid $\cB^\sigma$.
\end{proposition}

\subsubsection{Duality} 
The constructions in Section \ref{sec:Dri.twist} and
Section \ref{sec:dual.twist} are dual of each other, in the following sense.
\begin{proposition} 
For an algebra $\R$, let $\cB$ be a left $\R$-bialgebroid which is a finitely
generated and projective right $\R$-module via left multiplication by the
target  map $t$, and consider the right dual right $\R$-bialgebroid
$\cB^*$. The following statements hold. 

(1) An element $J=\sum_k \xi_k\sstac \R \zeta_k\in B^*\times_\R B^*$ is
    an invertible normalised 2-cocycle in $\cB^*$ if and only if
\begin{equation}\label{eq:J>sigma}
\sigma_J:B\stac {\R\otimes \R^{op}} B \to R,\qquad 
\sigma_J(b,b'):=\sum_{k} \xi_k\big(b t(\zeta_k(b'))\big)
\end{equation}
is an invertible normalised 2-cocycle on $\cB$. 

(2) Assume the equivalent properties in part (1).
    The right bialgebroid $(\cB^*)_J$, obtained by twisting the coproduct of
    $\cB^*$ by the cocycle $J$, is right dual of the left bialgebroid
    $\cB^{\sigma_J}$, obtained by twisting the product in $\cB$ by the
    cocycle $\sigma_J$ in \eqref{eq:J>sigma}. 
\end{proposition}
The inverse of the construction in \eqref{eq:J>sigma} is given by 
associating to an invertible normalised 2-cocycle $\sigma$ on $\cB$ an
invertible normalised 2-cocycle in $\cB^*$:
$J_\sigma:= \sum_i \sigma(-,a_i)\sstac \R \alpha_i$, where $(\{\
a_i\in B\ \}$, $\{\ \alpha_i\in B^*\ \})$ is a dual basis.

\subsubsection{Drinfel'd double}

For a Hopf algebra $B$ over a commutative ring $k$, such that $B$ is a
finitely generated and projective $k$-module, the $k$-module ${\mathcal
  D}(B):=B\sstac k B^*$ has a bialgebra (in fact Hopf algebra) structure. It
is known as the {\em Drinfel'd double} of $B$. The category of ${\mathcal
  D}(B)$-modules is isomorphic (as a monoidal category) to the category of
Yetter-Drinfel'd modules for $B$ and also to the monoidal center
of the category of $B$-modules. These results were extended to certain
bialgebroids by Schauenburg in \cite{Scha:dual.double.bgd}.

In this section let $\cB$ be a left bialgebroid over a $k$-algebra $\R$,
finitely generated and projective as a right $\R$-module. Assume in addition 
that the following map is bijective.
\begin{equation}\label{eq:x-Hopf.alg.map}
\vartheta: B\stac {\R^{op}} B \to B \stac \R B,\qquad 
b\stac {\R^{op}} b' \mapsto b_{(1)}\stac \R b_{(2)} b',
\end{equation}
where in the domain of $\vartheta$ module structures are given by right and
left multiplications by the target map, and in the codomain module structures
are given by left multiplications by the  target and source maps. Left
bialgebroids, for which the map \eqref{eq:x-Hopf.alg.map} is bijective, were
named {\em (left) $\times_R$-Hopf algebras} in 
\cite{Scha:dual.double.bgd} and they are discussed in more detail in
Section \ref{sec:x_R.hopf.alg}. In following Proposition
\ref{prop:double.bgd}, Sweedler's index notation is used for the coproducts,
and also the index notation $\vartheta^{-1}(b\sstac \R 1_\B)=b^{\langle
1\rangle}\sstac {\R^{op}}b^{\langle 2\rangle}$, where implicit summation is
understood.  
\begin{proposition}\label{prop:double.bgd}
Let $\R$ be an algebra over a commutative ring $k$. Let $\cB$ be a left
$\times_R$-Hopf algebra which is a finitely generated and projective right
$\R$-module via left multiplication by the target map. 
Denote the structure maps in $\cB$ as in Definition \ref{def:left.bgd}.
Consider $B$ as a right
$\R\sstac k \R^{op}$-module via right multiplications by the source and target
maps $s$ and $t$ in the left bialgebroid $\cB$. Consider the right dual $B^*$
as a left $\R\sstac k \R^{op}$-module via right multiplications by the target
and source maps $t^*$ and $s^*$ in the right bialgebroid $\cB^*$. The tensor
product ${\mathcal D}(\cB):= B\sstac {\R\otimes \R^{op}} B^*$ has a left
$\R$-bialgebroid structure, as follows. Multiplication is given by
$$
(b\bowtie \beta)(b' \bowtie \beta'):=
b s\big(\beta^{(1)}(b'_{(1)})\big) {b'_{(2)}}^{\langle 1\rangle}
\bowtie
\beta' \beta^{(2)} s^*\big(\beta^{(3)}( {b'_{(2)}}^{\langle 2\rangle})\big),
\qquad \textrm{for } b\bowtie \beta,\ b' \bowtie\beta'\in  {\mathcal D}(\cB).
$$
Source and target maps are 
$$
\R\to {\mathcal D}(\cB),\qquad r\mapsto 1_\B \bowtie t^*(r)
\qquad \textrm{and}\qquad 
\R^{op}\to {\mathcal D}(\cB),\qquad r\mapsto 1_\B \bowtie
s^*(r),
$$
respectively.
Coproduct is
$$
{\mathcal D}(\cB)\to {\mathcal D}(\cB)\stac \R {\mathcal D}(\cB),\qquad 
b\bowtie \beta \mapsto 
(b_{(1)}\bowtie \beta^{(1)})\stac \R
(b_{(2)}\bowtie \beta^{(2)}),
$$
and counit is
${\mathcal D}(\cB)\to R$,
$b\bowtie \beta \mapsto \epsilon\big( b s(\beta(1_\B))\big)$.
\end{proposition}
It is not known if the Drinfel'd double (or the dual) of a $\times_R$-Hopf
algebra is a $\times_R$-Hopf algebra too.

\subsubsection{Morita base change}\label{sec:Mor.base.ch}
In contrast to the previous sections, constructions in this section and in the
forthcoming ones {\em change the base algebra of a bialgebroid}.

Let $\R$ be an algebra over a commutative ring $k$ and let $\cB$ be a left
$\R$-bialgebroid. Denote the structure maps of $\cB$ as in Definition
\ref{def:left.bgd}. Let ${\widetilde \R}$ be a $k$-algebra which is Morita
equivalent to $\R$. Fix a strict Morita context
$
( \R\ ,\ {\widetilde \R}\ ,\ P\ ,\ Q\ ,\ \bullet\ ,\ \circ\ ).
$
Denote the inverse image of $1_\R$ under the map $\bullet$ by $\sum_i
q_i\sstac {\widetilde \R} p^i \in Q\sstac {\widetilde \R} P$ and denote the
inverse 
image of $1_{\widetilde \R}$ under the map $\circ$ by $\sum_j p_j\sstac \R q^j
\in P\sstac \R Q$. The ${\widetilde \R}$-$\R$ bimodule $P$ determines a
canonical $\R^{op}$-${\widetilde \R}^{op}$ bimodule $\overline P$. Similarly,
there is a canonical ${\widetilde \R}^{op}$-$\R^{op}$ bimodule $\overline Q$.
In \cite[Section 5]{Scha:Mor.ch}, a left ${\widetilde \R}$-bialgebroid
structure was constructed on the $k$-module 
$
{\widetilde B}:=(P\sstac k {\overline Q})\otimes B \otimes (Q\sstac k
{\overline P}), 
$
where unadorned tensor product is meant over $\R\sstac k
\R^{op}$. Multiplication is given by 
$$
[(p_1\otimes q_1)\otimes b \otimes (q_2\otimes p_2)]
[(p'_1\otimes q'_1)\otimes b' \otimes (q'_2\otimes p'_2)]:=
(p_1\otimes q_1)\otimes b s(q_2 \bullet p'_1) t(q'_1\bullet p_2) b'
\otimes (q'_2\otimes p'_2).
$$
Source and target maps in ${\widetilde \cB}$ are, for
${\widetilde r} \in {\widetilde R}$, 
$$
{\widetilde r}\mapsto \sum_{j,j'} 
({\widetilde r}\cdot p_{j'}\otimes q^j)\otimes 1_\B \otimes (q^{j'}\otimes
p_j) \qquad \textrm{and} \qquad
\qquad {\widetilde r}\mapsto \sum_{j,j'} (p_{j'}\otimes q^j\cdot {\widetilde
  r})\otimes 1_\B \otimes (q^{j'}\otimes p_j), 
$$
respectively. Coproduct and counit are  given on an element $ (p_1\otimes
q_1)\otimes b \otimes (q_2\otimes p_2) \in {\widetilde B}$ as
\begin{eqnarray*}
&& (p_1\otimes q_1)\otimes b \otimes (q_2\otimes p_2) \mapsto \sum_{i,j}
[(p_1\otimes q_i)\otimes b_{(1)} \otimes (q_2\otimes p_j)]\stac {\widetilde \R}
[(p^i\otimes q_1)\otimes b_{(2)} \otimes (q^j\otimes p_2)],\\
&& (p_1\otimes q_1)\otimes b \otimes (q_2\otimes p_2) \mapsto
p_1 \cdot \epsilon\big(b s(q_2\bullet p_2)\big)\circ q_1,
\end{eqnarray*}
respectively.
Generalisation, and a more conceptual background of Morita base change in
bialgebroids, are presented in Section \ref{sec:bgd.mod}.

\subsubsection{Connes-Moscovici's bialgebroids}\label{sec:Con.Mos.ex}
The following bialgebroid was constructed in \cite{ConMos:trans.geo} in the
framework of transverse geometry. 
Let $H$ be a Hopf algebra over a commutative ring $k$, with coproduct
$\delta:h \mapsto h_{(1)}\sstac k h_{(2)}$
and counit $\varepsilon$. Let $\R$ be a left $H$-module algebra. Consider the
$k$-module $B:=R\sstac k H \sstac k R$. It can be equipped with an associative 
multiplication
$$
(r_1\stac k h \stac k r_2)(r'_1\stac k h' \stac k r'_2) :=
r_1 (h_{(1)}\cdot r'_1)\stac k h_{(2)} h'\stac k (h_{(3)}\cdot r'_2) r_2,
$$
with unit $1_\R\sstac k 1_{\bf H}\sstac k 1_\R$.
The algebra $\B$ can be made a left $\R$-bialgebroid with source and target
maps  
$$
\R\to \B,\quad r\mapsto r\stac k 1_{\bf H} \stac k 1_\R\qquad
\textrm{and}\qquad  
\R^{op}\to \B,\quad r\mapsto 1_\R\stac k 1_{\bf H} \stac k r,
$$
respectively. Coproduct and counit are
$$
r\stac k h \stac k r' \mapsto (r\stac k h_{(1)}\stac k 1_\R) \stac \R
(1_\R\stac k h_{(2)}\stac k r')\qquad \textrm{and}\qquad 
r\stac k h \stac k r' \mapsto r\varepsilon(h)r'.
$$
  
\subsubsection{Scalar extension}\label{sec:scalar.ext}
Let $B$ be a bialgebra over a commutative ring $k$ and let $A$ be an algebra
in the category of right-right Yetter-Drinfel'd modules of $B$. Recall that
this means that $A$ is a right $B$-module algebra and a right $B$-comodule
algebra such that the following compatibility condition holds, for $a\in A$
and $b\in B$.
$$
(a\cdot b_{(2)})_{[0]}\stac k b_{(1)} (a\cdot b_{(2)})_{[1]} = 
a_{[0]} \cdot b_{(1)} \stac k a_{[1]} b_{(2)}.
$$
The category of right-right Yetter-Drinfel'd modules is pre-braided. Assume
that $A$ is braided commutative, i.e., for $a,a'\in A$, the identity
$
a'_{[0]}(a\cdot a'_{[1]}) =aa'
$
holds.
Under these assumptions, it follows by a symmetrical version of \cite[Theorem
  4.1]{BrzeMil:bgd} that the smash product algebra has a 
right $\A$-bialgebroid structure. Recall that
the smash product algebra is the $k$-module $A\# B:=A\sstac k B$, with
multiplication
$(a\# b)(a'\# b'):=a'(a\cdot b'_{(1)})\stac k b b'_{(2)}$.
Source and target maps are
$$
s:\A\to A\# B, \qquad a\mapsto a_{[0]}\# a_{[1]}\qquad \textrm{and}\qquad
t:\A^{op}\to A\# B, \qquad a\mapsto a\# 1_\B,
$$
respectively. Coproduct is
$$
\Delta:A\# B\to(A\# B)\stac \A (A\# B),\qquad a\# b\mapsto (a\# b_{(1)})\stac
\A (1_\A\# b_{(2)})
$$
and counit is given in terms of the counit $\varepsilon$ in $B$ as 
$
A\# B \to A$, $a\# b \mapsto a\varepsilon (b)
$.
The name {\em `scalar extension'} comes from the feature that the base algebra
$k$ of $B$ (the subalgebra of `scalars') becomes replaced by the base
algebra ${\bf A}$ of $A\# B$.

A solution $\mathcal{R}$ of the quantum Yang-Baxter equation on a finite
dimensional vector space determines a bialgebra $B(\mathcal{R})$ via the so
called FRT construction. In \cite[Proposition 4.3]{BrzeMil:bgd} a braided
commutative algebra in the category of Yetter-Drinfel'd modules for
$B(\mathcal{R})$ was constructed, thus a bialgebroid was associated to a
finite dimensional solution $\mathcal{R}$ of the quantum Yang-Baxter equation.

Above construction of a scalar extension was extended in \cite[Theorem
4.6]{BalSzl:fin.Gal}. Following it, a smash product of a
right bialgebroid $\cB$ over an algebra $\R$, with a braided commutative
algebra $A$ in the category of right-right Yetter-Drinfel'd modules for
$\cB$, is shown to possess a right ${\bf A}$-bialgebroid structure. The
fundamental importance of scalar extensions from the point of view 
of Galois extensions by bialgebroids is discussed in Section
\ref{sec:bgd.gal}. 

\subsection{The monoidal category of modules}\label{sec:bgd.mod}
An algebra $\B$ over a commutative ring $k$ is known to have a $k$-bialgebra
structure if and only if the category of (left or right) $\B$-modules is
monoidal such that the forgetful functor to $\cM_k$ is strict
monoidal \cite{Par:nc}. Generalisation \cite[Theorem 5.1]{Scha:ba.nc.base} of
this fact to 
bialgebroids is due to Schauenburg. Recall that any (right) module of an
$\R\sstac k \R^{op}$-ring $(\B,s,t)$ is an $\R$-bimodule via the actions by
$t(r)$ and $s(r)$, for $r\in R$.
\begin{theorem}\label{thm:bgd.mod.cat}
For an algebra $\R$ over a commutative ring $k$, the following data are
equivalent on an $\R\sstac k \R^{op}$-ring $(\B,s,t)$.

(1) A right bialgebroid structure on $(\B,s,t)$;

(2) A monoidal structure on the category $\cM_\B$ of right $\B$-modules, such
that the forgetful functor $\cM_\B\to {}_\R \cM_\R$ is strict monoidal.
\end{theorem}
Applying Theorem \ref{thm:bgd.mod.cat} to the opposite $\R\sstac k
\R^{op}$-ring, an analogous equivalence is obtained between left bialgebroid
structures and monoidal structures on the category of left modules. At the
heart of Theorem \ref{thm:bgd.mod.cat} lies the fact that the right regular
$\B$-module is a generator in the category $\cM_\B$. Hence the $\B$-module 
structure on the $\R$-module tensor product of any two $\B$-modules can be
expressed in terms of the action on $1_\B\sstac \R 1_\B\in B\sstac \R B$,
which defines a coproduct. In terms of a coproduct $b\mapsto b^{(1)}\sstac \R
b^{(2)}:=(1_\B\sstac \R 1_\B)\cdot b$, the 
right $\B$-action in the $\R$-module tensor product of two right $\B$-modules
$M$ and $N$ is written as
\begin{equation}\label{eq:tens.mod}
(m\stac \R n)\cdot b=m\cdot b^{(1)}\stac \R n\cdot b^{(2)}, \qquad \textrm{for
}m\stac \R n\in M\stac \R N,\ b\in B.
\end{equation}
A $\B$-module structure on the monoidal unit $R$ is equivalent to a right
character $\epsilon:B\to R$ by Lemma \ref{lem:character}.

Recall that, for an $\R\sstac k \R^{op}$-ring $(\B,s,t)$, the category
$\cM_\B$ is isomorphic to the category of algebras for the monad $-\sstac
{\R\otimes \R^{op}} B$ on ${}_\R \cM_\R$ (where $\R$-bimodules are considered
as right $\R\sstac k \R^{op}$-modules). In light of this fact,
Theorem \ref{thm:bgd.mod.cat} is a particular case of a question
discussed by Moerdijk in \cite{Moe:monad.tens.cat}: Having a monad
$\mathbb{B}$ on a monoidal category $\cM$, the monoidal structure of $\cM$
lifts to a monoidal structure on the category $\cM_{\mathbb B}$ of ${\mathbb
B}$-algebras (in the sense that the forgetful functor $\cM_{\mathbb B}\to \cM$
is strict monoidal, as in part (2) in Theorem \ref{thm:bgd.mod.cat}) if and
only if ${\mathbb B}$ is a monoid in the category of op-monoidal endofunctors
on $\cM$, i.e. a {\em Hopf monad} in the terminology of
\cite{Moe:monad.tens.cat} (called a {\em  bimonad} in
\cite{Szlach:mon.Eil.Moore}). Comparing this result with Theorem
\ref{thm:bgd.mod.cat}, we obtain another evidence for a characterisation of 
bialgebroids as bimonads in Theorem \ref{thm:bimonad}.

In the paper \cite{Szl:MonMor} by Szlach\'anyi it was investigated what
bialgebroids possess monoidally equivalent module categories. That is, a {\em
  monoidal Morita theory} of bialgebroids was developed. Based on Theorem
\ref{thm:bimonad}, one main result in  \cite{Szl:MonMor} can be reformulated
as in Theorem \ref{thm:mon.Mor} below. 
Recall that a bimodule, for two monads $\mathbb{B}:\cM\to \cM$ and
$\mathbb{B}':\cM'\to \cM'$, is a functor $\mathbb{M}:\cM \to \cM'$ together
with natural transformations $\varrho:\mathbb{M}\mathbb{B}\to \mathbb{M}$ and
$\lambda:\mathbb{B}' \mathbb{M} \to \mathbb{M}$, satisfying the usual
compatibility conditions for the right and left actions in a bimodule. If any
pair of parallel 
morphisms in $\cM'$ possesses a coequaliser, and $\mathbb{B}'$ preserves
coequalisers, then the bimodule $\mathbb{M}$ induces a functor
$\overline{\mathbb{M}}:\cM_\mathbb{B}\to \cM'_{\mathbb{B}'}$, between the
categories of algebras for $\mathbb{B}$ and $\mathbb{B}'$, respectively, 
with object map
$(V,v)\mapsto\mathrm{Coeq}(\mathbb{M}(v),\varrho_V)$. If both categories
$\cM$ and $\cM'$ possess coequalisers and both $\mathbb{B}$ and $\mathbb{B}'$
preserve them, then one can define the inverse of a bimodule as in Morita's
theory. A $\mathbb{B}'$-$\mathbb{B}$ bimodule $\mathbb{M}'$ is said to be the
inverse of $\mathbb{M}$ provided that $\overline{\mathbb{M}'}\,
\overline{\mathbb{M}}$ is naturally equivalent to the identity functor on
$\cM_\mathbb{B}$ and $\overline{\mathbb{M}}\, \overline{\mathbb{M}'}$ is
naturally equivalent to the identity functor on $\cM_{\mathbb{B}'}$.
\begin{theorem}\label{thm:mon.Mor}
For two right bialgebroids $\cB$ and $\cB'$, over respective base algebras
$\R$ and $\R'$, the right module categories $\cM_\B$ and $\cM_{\B'}$ are
monoidally equivalent if and only if there exists an invertible bimodule 
in the 2-category of op-monoidal left adjoint functors, 
for the monads $-\sstac {\R\otimes \R^{op}} B:{}_\R \cM_\R\to {}_\R \cM_\R$
and  $-\sstac {\R'\otimes \R^{\prime op}} B':{}_{\R'} \cM_{\R'}\to
{}_{\R'} \cM_{\R'}$.
\end{theorem}
By standard Morita theory, an equivalence $\overline{\mathbb{M}}:\cM_\B \to
\cM_{\B'}$ is of the form $\overline{\mathbb{M}}= -\sstac \B M$, for some
invertible $\B$-$\B'$ bimodule $M$. In \cite{Szl:MonMor} monoidality of the
equivalence is translated to properties of the Morita equivalence bimodule
$M$. 

In \cite[Definition 2.1]{Tak:sqrMor} two algebras $\R$ and $\widetilde \R$
over a commutative ring $k$ were said to be {\em
$\sqrt{\textit{Morita}}$-equivalent} whenever the bimodule categories ${}_\R
\cM_\R$ and ${}_{{\widetilde \R}}\cM_{{\widetilde \R}}$ are strictly
equivalent as $k$-linear monoidal categories. This property implies that the
algebras $\R\sstac k \R^{op}$ and ${\widetilde \R}\sstac k {\widetilde
  \R}^{op}$ are  $\sqrt{\textrm{Morita}}$-equivalent (but not 
conversely). In this situation, any $\R\sstac k \R^{op}$-ring $B$ (i.e. monoid
in the category of  $\R\sstac k \R^{op}$-bimodules) 
determines an  ${\widetilde \R}\sstac k {\widetilde \R}^{op}$-ring $\widetilde
B$, with underlying $k$-algebra ${\widetilde \B}$ Morita equivalent to
$\B$. If $\cB$ is a right $\R$-bialgebroid then the forgetful functor $\cM_\B
\to 
{}_\R\cM_\R$ is strict monoidal by Theorem \ref{thm:bgd.mod.cat}. Hence the
equivalence $\cM_\B\cong \cM_{\widetilde \B}$ can be used to induce a monoidal
structure on $\cM_{\widetilde \B}$ such that the forgetful functor
$\cM_{\widetilde \B}\to {}_{{\widetilde \R}} \cM_{{\widetilde \R}}$ is strict
monoidal. By Theorem \ref{thm:bgd.mod.cat} we conclude that there is a right
${\widetilde \R}$-bialgebroid structure on ${\widetilde B}$. In
\cite{Scha:Mor.ch}, the bialgebroid ${\widetilde \cB}$ was said to be obtained
from $\cB$ via {\em $\sqrt{\textit{Morita}}$-base change}. Since Morita
equivalent algebras are also $\sqrt{\textrm{Morita}}$-equivalent (but not
conversely), the construction in Section \ref{sec:Mor.base.ch} is a special
instance of a $\sqrt{\textrm{Morita}}$-base change. 

For a right $\R$-bialgebroid $\cB$, with structure maps denoted as in
Definition \ref{def:right.bgd}, consider $B$ as an $\R$-bimodule (or
$\R^{op}$-bimodule) via right multiplications by the source map $s$ and the
target map $t$. Both $B\sstac \R B$ and $B\sstac {\R^{op}} B$ are left modules
for $\B\sstac k \B$ via the regular actions on the two factors.
Associated to $\cB$, we construct a category ${\mathcal
C}(\cB)$, with two objects $\circ$ and $\bullet$. Morphisms with source
$\circ$ are elements of $B\sstac \R B$ and morphisms with source
$\bullet$ are elements of $B\sstac {\R^{op}} B$. 
Morphisms $F$, with target $\circ$ and $\bullet$, are required to satisfy the 
following ($\R$-centralising) conditions $(T\circ)$ and $(T\bullet)$,
respectively, for all $r\in R$. 
\begin{itemize}
\item[{$(T\circ)$}]$\quad \qquad \qquad \qquad \qquad (s(r) \stac k 1_\B)\cdot
  F =(1_\B \stac k t(r))\cdot F$, 
\item[{$(T\bullet)$}]$\quad  \qquad \qquad \qquad \qquad (t(r) \stac k
  1_\B)\cdot F =(1_\B \stac k s(r))\cdot F$.
\end{itemize}
Via composition given by factorwise multiplication, ${\mathcal C}(\cB)$ is a
category. Unit morphisms at the objects $\circ$ and $\bullet$ are $1_\B\sstac
\R 1_\B$ and $1_\B\sstac {\R^{op}} 1_\B$, respectively. The range of the
coproduct $\Delta$ lies in $\mathrm{Hom}(\circ,\circ)=B\times_R B$ and the
range of the co-opposite coproduct $\Delta_{cop}$ lies in
$\mathrm{Hom}(\bullet,\bullet)=B_{cop}\times_{R^{op}} B_{cop}$.  
In terms of the category ${\mathcal C}(\cB)$, definition of a quasi-triangular
bialgebroid as formulated in \cite[Proposition 3.13]{DonMud:Qgr} can be
described as follows.
\begin{definition}\label{def:R-matrix}
For a right $\R$-bialgebroid $\cB$, with structure maps denoted as in
Definition \ref{def:right.bgd}, let ${\mathcal C}(\cB)$ be the category
constructed above. An invertible morphism $\mathcal {R}=\mathcal {R}^1\sstac
{\R^{op}} \mathcal {R}^2\in \mathrm{Hom}(\bullet, \circ)$ (where implicit
summation is understood) is a {\em universal R-matrix} provided that for any
$b\in B$ the following identity holds in $\mathrm{Hom}(\bullet, \circ)$
$$
\Delta(b) \mathcal {R} = \mathcal {R} \Delta_{cop}(b)
$$
and 
$$
(\Delta_{cop} \,\stac{\R^{op}}\, B)(\mathcal {R})= \mathcal {R} \triangleleft
\mathcal {R} \qquad \textrm{and}\qquad 
(B\,\stac {\R^{op}}\, \Delta_{cop})(\mathcal {R})= \mathcal {R} \triangleright
\mathcal {R},
$$
where the (well defined) maps 
\begin{eqnarray*}
&-\ \triangleleft \ \mathcal {R}: \ B\stac {\R^{op}} B \to B\stac {\R^{op}}
  B\stac {\R^{op}} B,\qquad  
&b\,\stac {\R^{op}}\, b'\mapsto b\, \stac {\R^{op}}\, \mathcal {R}^1\, \stac
  {\R^{op}}\, b' \mathcal   {R}^2\qquad \textrm{and}\\
&\mathcal {R} \ \triangleright \ - :\ B\stac {\R^{op}} B \to B\stac {\R^{op}}
  B\stac {\R^{op}} B,\qquad 
&b\,\stac {\R^{op}}\ b'\mapsto b \mathcal {R}^1\,\stac {\R^{op}}\, \mathcal
  {R}^2 \,\stac {\R^{op}}\, b' 
\end{eqnarray*}
are used. A right bialgebroid $\cB$ with a given universal R-matrix
$\mathcal{R}$ is called a {\em quasi-triangular} bialgebroid. 
\end{definition}
Following Theorem \ref{thm:braid} was obtained in \cite[Theorem
  3.15]{DonMud:Qgr}, as a generalisation of an analogous result for 
quasi-triangular bialgebras. 
\begin{theorem}\label{thm:braid}
Consider a quasi-triangular right bialgebroid $(\cB,\mathcal{R})$ over a base
algebra $\R$. The monoidal category of right $B$-modules is braided, with
braiding natural isomorphism
$$ 
M\stac \R M' \to M'\stac \R M,\qquad m\stac \R m'\mapsto m'\cdot \mathcal{R}^2
\stac \R m\cdot \mathcal{R}^1.
$$
\end{theorem}

\subsection{The monoidal category of comodules}\label{sec:bgd_comod}
For a bialgebra over a commutative ring $k$, not only the category of
modules, but also the category of (left or right) comodules has a monoidal
structure, such that the forgetful functor to $\cM_k$ is strict monoidal. In
trying to prove an analogue of this result for bialgebroids, the first question
is to find a forgetful functor. A right, say, comodule of (the constituent
coring in) an $\R$-bialgebroid is by definition only a right $\R$-module. In
order to obtain a forgetful functor to the {\em monoidal} category of
$\R$-{\em bimodules}, following \cite[Lemma 1.4.1]{PHH:TK} is needed.
\begin{lemma}
Let $\R$ and ${\bf S}$ be two algebras over a commutative ring $k$ and let $C$
be an
${\bf S}|\R$-coring. Any right comodule $(M,\varrho^M)$ of the $\R$-coring $C$
can be equipped with a unique left ${\bf S}$-module structure such that
$\varrho^M(m)$ belongs to the center of the ${\bf S}$-bimodule $M\sstac \R C$,
for every $m\in M$. This unique left ${\bf S}$-action makes $M$ an ${\bf
S}$-$\R$ bimodule. Every $C$-comodule map becomes an ${\bf S}$-$\R$ bimodule
map. That is, there is a forgetful functor $\cM^C \to {}_{\bf S}\cM_\R$.
\end{lemma}
The left ${\bf S}$-action on a right comodule $(M,\varrho^M)$ of the
$\R$-coring $C$ is constructed as 
\begin{equation}\label{eq:S|R.coring.comod}
s\cdot m:= m^{[0]}\cdot \epsilon(m^{[1]} \cdot (s\stac k 1_\R)), \qquad
  \textrm{for   } s\in {\bf   S}, \ m\in M.
\end{equation}
In particular, \eqref{eq:S|R.coring.comod} can be used to equip a right
comodule of a right $\R$-bialgebroid with an $\R$-$\R$ bimodule structure. 
Applying the construction to co-opposite and opposite
bialgebroids, forgetful functors are obtained from categories of left and
right comodules of left and right $\R$-bialgebroids to ${}_\R\cM_\R$.
\begin{theorem}\label{thm:comod.cat.mon}
Let $\R$ be an algebra and let $\cB$ be a right
$\R$-bialgebroid. The category $\cM^B$ of right $\cB$-comodules is monoidal,
such that the forgetful functor $\cM^B \to {}_{\bf R}\cM_\R$ is strict
monoidal.  
\end{theorem}
The monoidal unit $R$ in ${}_\R\cM_\R$ is a right $\cB$-comodule, via a
coaction provided by the source map. One has to verify that, for any two right
$\cB$-comodules $M$ and $N$, the diagonal coaction 
\begin{equation}\label{eq:diag.coac}
M\stac \R N \to M\stac \R N \stac \R B, \qquad m\stac \R n \mapsto
m^{[0]}\stac \R n^{[0]}\stac \R m^{[1]} n^{[1]}
\end{equation}
is well defined and that the coherence natural isomorphisms in ${}_\R\cM_\R$
are $\cB$-comodule maps. 

We do not know about a converse of Theorem \ref{thm:comod.cat.mon}, i.e. an
analogue of the correspondence (2)$\Rightarrow$(1) in Theorem
\ref{thm:bgd.mod.cat} for the category of comodules. A reason  
for this is that (in contrast to modules of an $\R$-bialgebroid $\cB$, which
are algebras for the monad $-\sstac {\R\otimes \R^{op}} B$ on the monoidal
category 
of $\R$-bimodules), it is not known if comodules can be described as coalgebras
of a comonad on ${}_\R \cM_\R$.

However, the definition of a bialgebroid can be dualised in the sense of
reversing all arrows in the diagrams in $\cM_k$, expressing the axioms of a
bialgebroid over a $k$-algebra. For a flat $k$-coalgebra $C$, (e.g. when $k$
is a field), $C$-bicomodules constitute a monoidal category
${}^C\cM^C$. Monoidal structure is given by {\em cotensor products} -- a
notion dual to a module tensor product. That is, for $C$-bicomodules $M$ and
$N$, the cotensor product $M\Box_C N$ is the equaliser of the maps $\varrho^M
\ot_k N$ and $M\ot_k {}^N\varrho$, where $\varrho^M$ is the right coaction on
$M$ and ${}^N \varrho$ is the left coaction on $N$. Flatness of the $k$-module
$C$ implies that $M\Box_C N$ is a $C$-bicomodule via the left $C$-coaction on
$M$ and right $C$-coaction on $N$. In this case, dualisation of the
bialgebroid axioms leads to the notion of a {\em bicoalgebroid} over the
$k$-coalgebra $C$, see \cite{BrzeMil:bgd}. The relation
between $C$-bicoalgebroid structures on a comonoid in ${}^C\cM^C$, and strict
monoidal structures on the forgetful functor from its comodule category to
${}^C\cM^C$, is studied in \cite{Bimi:bicoa} and \cite{Bimi:PhD}. 

Applying Theorem \ref{thm:comod.cat.mon} to the co-opposite bialgebroid, we
conclude on the strict monoidality of the forgetful functor ${}^B \cM\to
{}_{\R^{op}}\cM _{\R^{op}}$, for a right $\R$-bialgebroid $\cB$. Applying
Theorem \ref{thm:comod.cat.mon} to opposite bialgebroids, it follows that the
forgetful functors ${}^B \cM \to {}_{\bf R}\cM_\R$ and $\cM^B\to
{}_{\R^{op}}\cM _{\R^{op}}$ are strict monoidal, for a left $\R$-bialgebroid
$\cB$. 
Note that, for an $\R$-bialgebroid $\cB$ which is a finitely generated and
projective $\R$-module on the appropriate side, the equivalence in Proposition
\ref{prop:dual.mod.comod} between the categories of comodules for $\cB$ and
modules for its dual is strict (anti-)monoidal.

The reader should be warned that in the paper \cite{BrugVire:Hopf.monad} a
different notion of a comodule is used. For an ${\bf S}|\R$-coring $C$, 
the coproduct and the counit of the $\R$-coring $C$ project to a coproduct and
a counit on the quotient $\R$-bimodule $C/ \{\ (s\sstac k 1_\R)\cdot c-c\cdot
(s\sstac k 1_\R)\ |\ c\in C,\ s\in S\ \}$. Applying the definition of a
comodule of a bimonad in \cite[Section 4.1]{BrugVire:Hopf.monad} to a bimonad
induced by a right $\R$-bialgebroid $\cB$, the resulting notion is a comodule
for 
the corresponding quotient coring $B/\{\ s(r)b - t(r)b\ |\ b\in B,\ r\in R\ \}$
(where $s$ and $t$ are the source and target map of $\cB$, respectively). The
category of such comodules is not known to be monoidal. 

\subsection{Algebra extensions by bialgebroids. Galois extensions}
\label{sec:bgd.gal}   
In analogy with bialgebra extensions,
there are two symmetrical notions of an algebra extension by a bialgebroid
$\cB$. In the {\em action picture} one deals with a $\cB$-module algebra $M$
and its invariant subalgebra (with respect to a character defined by the 
counit). In this picture, {\em Galois property} means Galois property of an
associated ${\bf M}$-ring. Dually, in the {\em coaction picture} one deals
with a $\cB$-comodule algebra $M$ and its coinvariant subalgebra (with
respect to a grouplike element defined by the unit). In this picture, {\em 
Galois property} means Galois property of an associated ${\bf
M}$-coring. Although the two approaches are symmetric (and equivalent for
finitely generated and projective bialgebroids), the coaction picture is  
more popular and more developed. We present it in more detail but, for
the sake of completeness, we shortly describe the action picture as
well. 

\subsubsection{The action and coaction pictures}\label{sec:ac.coac.pic}
By Theorem \ref{thm:bgd.mod.cat}, the category $\cM_\B$ of right modules of a
right 
$\R$-bialgebroid $\cB$ is monoidal. By definition, a {\em right $\cB$-module
algebra} is a monoid $M$ in $\cM_\B$. 
Denote the structure maps of $\cB$ as in Definition \ref{def:right.bgd}.
In view of Lemma \ref{lem:char.R-rings}, a right
$\cB$-module algebra is the same as an algebra and right $\B$-module ${\bf M}$,
such that the multiplication in ${\bf M}$ is $\R$-balanced and 
$$
(mm')\cdot b = (m\cdot b^{(1)})(m'\cdot b^{(2)})\qquad \textrm{and} \qquad 
1_{\bf M}\cdot b =1_{\bf M}\cdot s\big(\epsilon(b)\big),\qquad \textrm{for
}m,m'\in M, \ b\in B, 
$$
cf. \eqref{eq:tens.mod}. 
Note in passing that, by strict monoidality of the forgetful functor
$\cM_\B\to {}_\R\cM_\R$, a right $\cB$-module algebra $M$ has a canonical
$\R$-ring structure. Its unit is the map $\R\to {\bf M}$, $r\mapsto 1_{\bf
  M}\cdot s(r)=1_{\bf M}\cdot t(r)$. 

For a right $\cB$-module algebra $M$, $B\sstac \R M$
has an ${\bf M}$-ring structure. It is called a {\em smash product}, with
multiplication
$
(b\sstac \R m)(b'\sstac \R m')=bb^{\prime (1)} \sstac \R (m\cdot b^{\prime
  (2)}) m',
$
and unit $m\mapsto 1_\B\sstac \R m$. The right
character $\epsilon$ on the $\R$-ring $(\B,s)$ determines a right character
$\epsilon\sstac \R M$ on the ${\bf M}$-ring $B\sstac \R M$. Hence we can
consider the invariant subalgebra of the base algebra ${\bf M}$, with respect
to the right character $\epsilon\sstac \R M$. It coincides with the
$\epsilon$-invariants of the $(\B,s)$-module $M$, 
$$
N:=M_{\epsilon} =\{\ n\in M\ |\ n\cdot b = n \cdot s\left(\epsilon
(b)\right),\quad \forall b\in B\ \}.
$$
In the {\em action picture} the algebra ${\bf M}$ is said to be a {\em
right $\cB$-Galois extension} of the invariant subalgebra ${\bf N}$ provided
that $B\sstac \R M$ is a 
  Galois ${\bf M}$-ring with respect to the right character $\epsilon\sstac \R
  M$. That is, the {\em canonical map}
$$
B\stac \R M \to {}_{\bf N} \mathrm{End}(M),\qquad 
b\stac \R m \mapsto \big(\ m'\mapsto (m'\cdot b)m\ \big)
$$
is bijective. Left Galois extensions by a left bialgabroid $\cB$ are defined
symmetrically, referring to a left $\cB$-module algebra and its
$\epsilon$-invariant subalgebra. 

By Theorem \ref{thm:comod.cat.mon}, also the category $\cM^B$ of right
comodules of a right $\R$-bialgebroid $\cB$ is monoidal. By definition, a {\em
right $\cB$-comodule algebra} is a monoid in $\cM^B$. 
In view of Lemma \ref{lem:char.R-rings}, a right 
$\cB$-comodule algebra is the same as an algebra and right $\cB$-comodule ${\bf
  M}$, with coaction $m\mapsto m^{[0]}\sstac \R m^{[1]}$, such that the
multiplication in ${\bf M}$ is $\R$-balanced and, for $m,m'\in M$, 
\begin{equation}\label{eq:com.alg}
(mm')^{[0]}\stac \R (mm')^{[1]}=m^{[0]} m^{\prime [0]}\stac \R m^{[1]}
m^{\prime [1]}\qquad \textrm{and}\qquad {1_{\bf M}}^{[0]}\stac \R {1_{\bf
  M}}^{[1]} = 1_{\bf M}\stac \R 1_\B.
\end{equation}
Note in passing that, by strict monoidality of the forgetful functor
$\cM^B\to {}_\R\cM_\R$, a right $\cB$-comodule algebra $M$ has a canonical
$\R$-ring structure. Its unit is the map $\R\to {\bf M}$, $r\mapsto 1_{\bf
  M}\cdot r=r\cdot 1_{\bf M}$. 

For a right $\cB$-comodule algebra $M$, $M\sstac \R B$ has an ${\bf M}$-coring
structure with left and right ${\bf M}$-actions
\begin{equation}\label{eq:comod.alg.coring}
m_1\cdot(m\stac \R b)\cdot m_2 = m_1 m {m_2}^{[0]} \stac \R b {m_2}^{[1]},
\qquad \textrm{for } m_1,m_2\in M, \ m\stac \R b\in M\stac \R B,
\end{equation}
comultiplication
$
m\sstac \R b \mapsto (m\sstac \R b^{(1)})\sstac {\bf M} (1_{\bf M}\sstac \R
b^{(2)})
$
and counit $m\sstac \R b \mapsto m\cdot \epsilon(b)$. 
The grouplike element $1_\B$ in the $\R$-coring $(B,\Delta,\epsilon)$
determines a grouplike element $1_{\bf M}\sstac \R 1_\B$ in the ${\bf
M}$-coring $M\sstac \R B$. Hence we can consider the coinvariant subalgebra
of the base algebra ${\bf M}$, with respect to the grouplike element $1_{\bf
  M}\sstac \R 1_\B$. It coincides with the $1_\B$-coinvariants of the
$\cB$-comodule $M$,  
$$
N:=M^{1_\B}=\{\ n\in M\ |\ n^{[0]}\stac \R n^{[1]}=n\stac \R 1_\B\ \}.
$$
Note that, by right $\R$-linearity of the $\cB$-coaction on $M$ and
\eqref{eq:com.alg}, for $n\in N$ and $r\in R$, 
$$
(n\cdot r)^{[0]}\stac \R (n\cdot r)^{[1]} = n\stac \R s(r) =
(r\cdot n)^{[0]}\stac \R (r\cdot n)^{[1]}.
$$
Hence, for $n\in N$ and $r\in R$, 
\begin{equation}\label{eq:n-r_comm}
n(1_{\bf M}\cdot r)=n \cdot r=r\cdot n =(1_{\bf M}\cdot r)n.
\end{equation}

In the {\em coaction picture} the algebra ${\bf M}$ is said to be a {\em right 
$\cB$-Galois extension} of the coinvariant subalgebra ${\bf N}$ provided that
$M\sstac \R B$ is a Galois ${\bf M}$-coring with respect to the grouplike
element $1_{\bf M}\sstac \R 1_\B$. That is, the {\em canonical map} 
\begin{equation}\label{eq:can}
\mathrm{can}:M\stac {\bf N} M \to M\stac \R B, \qquad 
m\stac {\bf N} m'\mapsto m m^{\prime [0]} \stac \R m^{\prime [1]}
\end{equation}
is bijective. Since for a right $\R$-bialgebroid $\cB$ also the category of
left comodules is monoidal, there is a symmetrical notion of a left
$\cB$-Galois extension ${\bf N}\subseteq {\bf M}$. It is a left $\cB$-comodule
algebra $M$, with coinvariant subalgebra ${\bf N}$, such that an 
associated ${\bf M}$-coring $B\sstac \R M$ is a Galois coring, with respect to
the grouplike element $1_{\bf B}\sstac \R 1_{\bf M}$. Left and right Galois
extensions by left bialgebroids are treated symmetrically.

For a right comodule algebra $M$ of a right $\R$-bialgebroid $\cB$, a {\em
right-right relative Hopf module} is a right $M$-module {\em in} $\cM^B$.
The category of right-right relative Hopf modules is denoted by $\cM^B_M$
and it turns out to be isomorphic to the category of right comodules for the
${\bf M}$-coring \eqref{eq:comod.alg.coring}.   
Hence the grouplike element $1_{\bf M}\sstac \R 1_\B\in M\sstac \R B$
determines an adjunction as in \eqref{eq:comod.adj} between $\cM^B_M$ and 
the category $\cM_{\bf N}$ of right modules for the coinvariant subalgebra
${\bf N}$ of $M$. It will be denoted as
\begin{equation}\label{eq:R.adj}
-\stac {\bf N} M:\cM_{\bf N} \to \cM^B_M
\qquad \textrm{and}\qquad 
(-)^{coB}: \cM^B_M \to \cM_{\bf N}. 
\end{equation}
Recall from Section \ref{sec:coring} that, for a $\cB$-Galois extension ${\bf
N}\subseteq {\bf M}$, this adjunction is interesting from the descent theory
point of view. 

For a finitely generated and projective bialgebroid, the action and coaction
pictures are equivalent in the sense of Proposition
\ref{prop:dual.ac.coac.pic}. This equivalence was observed (in a slightly more
restricted context) in \cite[Theorem \& Definition 3.3]{BalSzl:fin.Gal}.
\begin{proposition}\label{prop:dual.ac.coac.pic}
Let $\cB$ be a right $\R$-bialgebroid which is a finitely generated and
projective right $\R$-module via right multiplication by the source map.

(1) There is a bijective correspondence between right $\cB$-module algebra
    structures and right $(\cB^*)^{op}$-comodule algebra structures on a given
    algebra ${\bf M}$.  

(2) The invariant subalgebra ${\bf N}$ of a right $\cB$-module algebra ${M}$
    (with respect to the right character given by the counit) is the same as
    the coinvariant subalgebra of the corresponding right
    $(\cB^*)^{op}$-comodule algebra  $M$ (with respect to the grouplike element
    given by the unit). 

(3) A right $\cB$-module algebra ${\bf M}$ is a $\cB$-Galois extension of its
    invariant subalgebra ${\bf N}$ in the action picture if and only if ${\bf
    M}$ is a $(\cB^*)^{op}$-Galois extension of ${\bf N}$ in the coaction
    picture. 
\end{proposition}
Part (1) of Proposition \ref{prop:dual.ac.coac.pic} follows by the strict
monoidal equivalence $\cM_\B\cong \cM^{(B^*)^{op}}$. Parts (2) and (3) follow
by Proposition \ref{prop:dual.Galois.ring-coring}, since the ${\bf M}$-ring
$B\sstac \R M$, associated to a right $\cB$-module algebra $M$, is the left
${\bf M}$-dual of the ${\bf M}$-coring $M\sstac \R (B^*)^{op}$, associated to
the right $(\cB^*)^{op}$-comodule algebra $M$.

Consider a right $\R$-bialgebroid $\cB$, which is a finitely generated and
projective right $\R$-module via the source map.
Then, for a right $\cB$-module algebra $M$, the category 
of right-right $(M,(\cB^*)^{op})$ relative
Hopf modules is equivalent also to the category of right modules for the smash
product algebra $B\sstac \R M$.

In the rest of these notes {\em only coaction picture} of Galois extensions
will be used. 

\subsubsection{Quantum torsors and bi-Galois extensions}
Following the work of Grunspan and Schauenburg \cite{Grun:Qtor},
\cite{Scha:Qtor}, \cite{Scha:biGal}, \cite{Scha:HbiGal}, for a bialgebra $B$
over a 
commutative ring $k$, a faithfully flat right $B$-Galois extension ${\bf T}$
of $k$ can be described without explicit mention of the
bialgebra $B$. Instead, a {\em quantum torsor} structure is introduced on
${\bf T}$, from which $B$ can be reconstructed uniquely. What is more, a
quantum torsor determines a second $k$-bialgebra $B'$, for which ${\bf T}$ is
a left 
$B'$-Galois extension of $k$. It is said that any faithfully flat Galois
extension of $k$ by a $k$-bialgebra is in fact a {\em bi-Galois
extension}. The categories of (left) comodules for the bialgebras $B$ and
$B'$ are monoidally equivalent.
Such a description of faithfully flat Galois extensions by bialgebroids was
developed in the PhD thesis of Hobst \cite{Hobst:phd} and in the paper
\cite{BohmBrz:torsor}.
\begin{definition}\label{def:torsor}
For two algebras ${\bf R}$ and ${\bf S}$ over a commutative ring $k$, an
$\R$-${\bf S}$ torsor is a pair $(T,\tau)$. Here $T$ is an $\R\sstac k {\bf
  S}$-ring with underlying $k$-algebra ${\bf T}$ and unit maps $\alpha:\R\to
{\bf T}$ and $\beta:{\bf S}\to {\bf T}$ (with commuting ranges in
${\bf T}$). Considering $T$ as an $\R$-${\bf S}$ bimodule and as an ${\bf
  S}$-$\R$ bimodule via the maps $\alpha$ and $\beta$, $\tau$ is an ${\bf
  S}$-$\R$ bimodule map $T\to T\sstac \R T\sstac {\bf S} T$, 
$t\mapsto t^{\langle 1\rangle} \sstac \R t^{\langle 2\rangle} \sstac {\bf S}
t^{\langle 3\rangle}$ (where implicit summation is understood), satisfying the
following axioms, for $t,t'\in T$, $r\in R$ and $s\in S$.
\begin{itemize}
\item[(i)] $(\tau\sstac\R T\sstac{\bf S}T)\circ \tau = (T \sstac\R T\sstac
  {\bf S} \tau )\circ \tau$ ({\em coassociativity}),
\item[(ii)] $(\mu_\R\sstac{\bf S} T)\circ \tau =\beta \sstac{\bf S} T$
and $(T\sstac\R \mu_{\bf S})\circ \tau = T\sstac\R \alpha$ ({\em left and
  right counitality}),
\item[(iii)] $\tau(1_{\bf T})=1_{\bf T} \sstac\R 1_{\bf T}\sstac{\bf S}
  1_{\bf  T}$ ({\em unitality}), 
\item[(iv)] $\alpha(r) t^{\langle 1\rangle} \sstac\R t^{\langle 2\rangle}
  \sstac{\bf S} t^{\langle 3\rangle} = 
 t^{\langle 1\rangle} \sstac\R t^{\langle 2\rangle} \alpha(r) \sstac
  {\bf S} t^{\langle 3\rangle}$ and \\
$t^{\langle 1\rangle} \sstac\R \beta(s)t^{\langle 2\rangle} \sstac
  {\bf S} t^{\langle 3\rangle}=
t^{\langle 1\rangle} \sstac\R t^{\langle
    2\rangle} \sstac{\bf S} t^{\langle 3\rangle}\beta(s)$, ({\em centrality
  conditions}) 
\item[(v)] $\tau(tt')= t^{\langle 1\rangle}t^{\prime \langle 1\rangle} \sstac
    \R t^{\prime\langle 2\rangle} t^{\langle 2\rangle}\sstac{\bf S}
     t^{\langle 3\rangle}t^{\prime \langle 3\rangle}$ ({\em multiplicativity}),
\end{itemize}
where $\mu_\R$ and $\mu_{\bf S}$ denote multiplication in the $\R$-ring $({\bf
T},\alpha)$ and the ${\bf S}$-ring $({\bf T},\beta)$, respectively.

An $\R$-${\bf S}$ torsor $(T,\tau)$ is said to be {\em faithfully flat} if $T$
is a faithfully flat right $\R$-module and a faithfully flat left ${\bf
  S}$-module. 
\end{definition}
Note that axiom (iv) in Definition \ref{def:torsor} is needed in order for
the multiplication in axiom (v) to be well defined.
\begin{theorem}\label{thm:torsor}
For two $k$-algebras ${\bf R}$ and ${\bf S}$, there is a bijective
correspondence between the following sets of data.
\begin{itemize}
\item[(i)] Faithfully flat $\R$-${\bf S}$ torsors $(T,\tau)$.
\item[(ii)] Right $\R$-bialgebroids $\cB$ and left faithfully flat right
  $\cB$-Galois extensions ${\bf S}\subseteq {\bf T}$, such that $T$ is
  a right faithfully flat $\R$-ring.
\item[(iii)] Left ${\bf S}$-bialgebroids $\cB'$ and right faithfully flat left 
  $\cB'$-Galois extensions ${\bf R}\subseteq {\bf T}$, such that $T$ is
  a left faithfully flat ${\bf S}$-ring.
\end{itemize}
Furthermore, a faithfully flat $\R$-${\bf S}$ torsor $T$ is a $\cB'$-$\cB$
bicomodule, i.e. the left $\cB'$, and right $\cB$-coactions on $T$ do commute.
\end{theorem}
Starting with the data in part (ii) of Theorem \ref{thm:torsor}, a torsor map
on $T$ is constructed in terms of the $\cB$-coaction $\varrho^T:T\to T\sstac \R
B$, and the inverse of the canonical map \eqref{eq:can} (with the role of the
comodule algebra $M$ in \eqref{eq:can} played by $T$), as $\tau:= (T\sstac \R
\mathrm{can}^{-1}(1_{\bf T}\sstac \R -))\circ \varrho^T$. Conversely, to a
faithfully flat $\R$-${\bf S}$ torsor $(T,\tau)$ (with multiplication
$\mu_{\bf S}$ in the ${\bf S}$-ring $({\bf T},\beta)$) one associates a right
$\R$-bialgebroid $\cB$, defined on the $\R$-$\R$ bimodule given by the
equaliser 
of the maps $(\mu_{\bf S}\sstac \R T \sstac {\bf S} T)\circ (T\sstac {\bf S}
\tau)$ and $\alpha\sstac \R T \sstac {\bf S} T:T\sstac {\bf S} T \to T\sstac
\R T \sstac {\bf S} T$.
\begin{theorem}\label{thm:torsor.eq}
For two $k$-algebras ${\bf R}$ and ${\bf S}$, consider a faithfully flat
$\R$-${\bf S}$ torsor $(T,\tau)$. Let $\cB$ and $\cB'$ be the associated
bialgebroids in Theorem \ref{thm:torsor}. Assume that $T$ is a faithfully flat
right ${\bf S}$-module and $B'$ is a flat right ${\bf S}$-module. Then the
categories of left $\cB$-, and $\cB'$-comodules are monoidally equivalent. 
\end{theorem}
Note that the assumptions made about the right ${\bf S}$-modules $T$ and $B'$
in Theorem \ref{thm:torsor.eq} become redundant if working with one
commutative base ring $\R={\bf S}$ and equal unit maps $\alpha=\beta$.
The equivalence in Theorem \ref{thm:torsor.eq} is given by $T\Box_\cB -: {}^B
\cM \to {}^{B'}\cM$, a cotensor product with the $\cB'$-$\cB$ bicomodule
$T$. (Recall that the notion of a {\em cotensor product} is dual to the one of
module tensor product. That is, for a right $\cB$-comodule $(M,\varrho^M)$ and
a left $\cB$-comodule $(N,{}^N \varrho)$, $M\Box_\cB N$ is the equaliser of
the maps $\varrho^M \sstac \R N$ and $M\sstac \R {}^N\varrho$.)

\subsubsection{Galois extensions by finitely generated and projective
  bialgebroids}
The bialgebra $B$, for which a given algebra extension is $B$-Galois, is
non-unique. Obviously, there is even more possibility for a choice of
$\cB$ if it is allowed to be a bialgebroid. Still, as a main advantage of
studying Galois extensions by bialgebroids, in an appropriately finite case
all possible bialgebroids $\cB$ can be related to a canonical one. Following
Theorem \ref{thm:Gal.bgds} is a mild generalisation of \cite[Proposition
4.12]{BalSzl:fin.Gal}.
\begin{definition}\label{def:Yet.Dri.mod}
Let $\cB$ be a right bialgebroid over an algebra $\R$. A {\em right-right
  Yetter-Drinfel'd module} for $\cB$ is a right $\cB$-module and right
  $\cB$-comodule $M$ (with one and the same underlying $\R$-bimodule
  structure), such that the following compatibility condition holds.
$$
(m\cdot b^{(2)})^{[0]} \stac \R b^{(1)} (m\cdot b^{(2)})^{[1]}= m^{[0]}\cdot 
  b^{(1)}\stac \R m^{[1]} b^{(2)},\qquad \textrm{for }m\in  M,\ b\in B.
$$
\end{definition}
It follows by a symmetrical version of \cite[Proposition
  4.4]{Scha:dual.double.bgd} that the category of right-right Yetter-Drinfel'd
modules of a right bialgebroid $\cB$ is isomorphic to the weak center of the
monoidal category of right $\B$-modules. Hence it is monoidal and
pre-braided. Following the paper \cite{BalSzl:fin.Gal}, the 
construction in Section \ref{sec:scalar.ext} can be extended to a braided
commutative algebra $A$ in the category of right-right Yetter-Drinfel'd
modules of a right $\R$-bialgebroid $\cB$. That is, the smash product algebra
$A \# B$ can be proven to carry the structure of a right ${\bf
  A}$-bialgebroid, called a {\em scalar extension} of $\cB$ by $A$. 

In following Theorem \ref{thm:Gal.bgds}, the center of a bimodule $M$ of an
algebra $\R$ is
denoted by $M^\R$.
\begin{theorem}\label{thm:Gal.bgds}
For an algebra $\R$ consider a right $\R$-bialgebroid $\cB$ which is a
finitely  
generated and projective left $\R$-module via right multiplication by the
target map. Let ${\bf N}\subseteq {\bf M}$ be a right $\cB$-Galois
extension. Then ${\bf N}\subseteq {\bf M}$ is a right Galois extension by a
right bialgebroid $(M\sstac {\bf N} M)^{\bf N}$ over the base algebra $M^{\bf
N}$. What is more, the bialgebroid $(M\sstac {\bf N} M)^{\bf N}$ is isomorphic
to a scalar extension of $\cB$.
\end{theorem}
In proving Theorem \ref{thm:Gal.bgds} the following key ideas are used. First
a braided commutative algebra structure, in the category of right-right
Yetter-Drinfel'd modules for $\cB$, is constructed on $M^{\bf N}$. The
$\cB$-coaction on $M^{\bf N}$ is given by restriction of the $\cB$-coaction on
$M$. The $\B$-action on $M^{\bf N}$ is of the Miyashita-Ulbrich type, i.e. it
is given in terms of the inverse of the canonical map
\eqref{eq:can}. Introducing an index notation $\mathrm{can}^{-1}(1_{\bf M}
\sstac \R b)= b^{\{1\}}\sstac {\bf N} b^{\{2\}}$, for $b\in B$ (implicit
summation is understood), the right $\B$-action on $M^{\bf N}$ is $a\cdot b:=
b^{\{1\}} a b^{\{2\}}$, for  
$b\in B$ and $a\in M^{\bf N}$. Since in this way $M^{\bf N}$ is a braided
commutative algebra in the category of right-right Yetter-Drinfel'd modules
for $\cB$, there exists a right $M^{\bf N}$-bialgebroid $M^{\bf N}\sstac \R B$
(cf. Section \ref{sec:scalar.ext}). Restriction
of the $\cB$-canonical map \eqref{eq:can} establishes a bijection
$(M\sstac {\bf N} M)^{\bf N}\to M^{\bf N}\sstac \R B$. Hence it induces an 
$M^{\bf N}$-bialgebroid structure on  $(M\sstac {\bf N} M)^{\bf N}$, and also
an $(M\sstac {\bf N} M)^{\bf N}$-comodule algebra structure on $M$. After
checking that coinvariants of the $(M\sstac {\bf N} M)^{\bf N}$-comodule $M$
are precisely the elements of ${\bf N}$, the $(M\sstac {\bf N} M)^{\bf
  N}$-Galois property of the extension ${\bf N}\subseteq {\bf M}$ becomes
obvious: the $(M\sstac {\bf N} M)^{\bf N}$-canonical map differs from the
$\cB$-canonical map \eqref{eq:can} by an isomorphism. 

\subsubsection{Depth two algebra extensions}
Classical finitary Galois extensions of fields 
can be characterised inherently, by normality and separability properties,
without referring to the Galois group $G$. That is, a (unique upto isomorphism)
finite Galois group $G:=\mathrm{Aut}_K(F)$ is {\em determined} by any normal
and separable field extension $F$ of $K$. While no such inherent
characterisation of Galois extensions by 
(finitely generated and projective) bialgebras is known, a most important
achievement in the Galois theory of bialgebroids is a characterisation of
Galois extension by finitely generated and projective bialgebroids. A first
result in this direction was \cite[Theorem 3.7]{BalSzl:fin.Gal}. At the level
of generality presented here, it was proven in \cite[Theorem
  2.1]{Kad:note}.  

Following definition in \cite[Definition 3.1]{KadSzl:D2bgd} was abstracted 
from depth 2 extensions of $C^*$-algebras.
\begin{definition}
Consider an extension ${\bf N}\subseteq {\bf M}$ of algebras.
It is said to satisfy the {\em right (resp. left) depth 2} condition if 
the ${\bf M}$-${\bf N}$ bimodule (resp. ${\bf N}$-${\bf M}$ bimodule) $M\sstac
{\bf N} M$ is a direct summand in a 
finite direct sum of copies of $M$.
\end{definition}
Note that the right depth 2 property of an algebra extension ${\bf
  N}\subseteq {\bf M}$ is equivalent to the existence of finitely many
  elements $\gamma_k\in {}_{\bf N}\mathrm{End}_{\bf N}(M)\cong
  {}_{\bf M}\mathrm{Hom}_{\bf N}(M\sstac {\bf N} M,M)$ and $c_k \in (M\sstac
  {\bf N} M)^{\bf N}\cong {}_{\bf M}\mathrm{Hom}_{\bf N}(M,M\sstac {\bf N}
  M)$, the so called {\em right depth 2 quasi-basis}, satisfying the identity
\begin{equation}\label{eq:D2.quasi.basis}
\sum_k m\gamma_k(m')c_k = m\stac {\bf N} m'\qquad \textrm{for }m,m'\in M.
\end{equation}
\begin{definition}
An extension ${\bf N}\subseteq {\bf M}$ of algebras 
is {\em balanced} if all endomorphisms of $M$, as a left module for the
algebra ${\mathcal E}:=\mathrm{End}_{\bf N}(M)$, are given by right
multiplication by some element of ${\bf N}$. 
\end{definition}
\begin{theorem}\label{thm:bgd.Gal.ext}
For an algebra extension ${\bf N}\subseteq {\bf M}$, the following properties
are equivalent.
\begin{itemize}
\item[{(i)}] ${\bf N}\subseteq {\bf M}$ is a right Galois extension by
some right $\R$-bialgebroid $\cB$, which is a finitely generated and
projective left $\R$-module via right multiplication by the target map.
\item[{(ii)}] The algebra extension ${\bf N}\subseteq {\bf M}$ is balanced and
  satisfies the right depth 2 condition.  
\end{itemize}
\end{theorem}
If ${\bf N}\subseteq {\bf M}$ is a right Galois extension by a right
$\R$-bialgebroid $\cB$, then $M\sstac {\bf N} M\cong M\sstac \R B$ as ${\bf
  M}$-${\bf N}$ bimodules. Hence the
right depth two condition follows by finitely generated projectivity of the
left $\R$-module $B$. A left ${\mathcal E}$-module endomorphism of $M$ is
given by right multiplication by an element $x\in M^\R$, by 
the right ${\bf N}$-linearity of the maps, given by left 
multiplication by an element $m\in M$, and right multiplication by
$r\in R$. Since by \eqref{eq:n-r_comm} also the right action
\eqref{eq:dual.ac} on $M$ by $\phi\in {}^* B$ is a right ${\bf N}$-module map,
it follows that $x^{[0]} 
\phi(x^{[1]})=x\phi(1_\B)$, for all $\phi\in {}^* B$. Together with the
finitely generated projectivity of the left $\R$-module $B$, this implies that
$x$ belongs to the coinvariant subalgebra ${\bf N}$, hence the extension ${\bf
  N}\subseteq {\bf M}$ is balanced. 

By Theorem \ref{thm:Gal.bgds}, if ${\bf N}\subseteq {\bf M}$ is a right Galois
extension by some finitely generated projective right $\R$-bialgebroid $\cB$,
then it is a Galois extension by the (finitely generated projective) right
$M^{\bf N}$-bialgebroid $(M\sstac {\bf N} M)^{\bf N}$. Hence in the converse
direction a balanced algebra extension ${\bf N}\subseteq {\bf M}$, satisfying
the right depth 2 condition, is shown to be a right Galois extension by a
right $M^{\bf N}$-bialgebroid $(M\sstac {\bf N} M)^{\bf N}$, constructed in
\cite[Section 5]{KadSzl:D2bgd}. (In fact, in \cite{KadSzl:D2bgd} both the left
and right depth two properties are assumed. It is proven in \cite{Kad:note}
that the construction works for one sided depth two extensions as well.)
The coproduct in $(M\sstac {\bf N} M)^{\bf N}$ and its
coaction on $M$ are constructed in terms of the right depth 2 quasi-basis
\eqref{eq:D2.quasi.basis}. Let us mention that the only point in the proof,
where the balanced property is used, is to show that the  $(M\sstac {\bf N}
M)^{\bf N}$-coinvariants in $M$ are precisely the elements of $N$.

\section{Hopf algebroids} \label{sec:hgd}
A Hopf algebra is a bialgebra $H$ equipped with an additional {\em antipode}
map $H\to H$.  The antipode is known to be a bialgebra map from $H$ to the
opposite co-opposite of $H$. It does not seem to be possible to define a Hopf
algebroid based on this analogy. Starting with a, say left, bialgebroid $\hH$,
its opposite co-opposite $\hH^{op}_{cop}$ is a right bialgebroid. There is no
sensible notion of a bialgebroid map $\hH\to \hH^{op}_{cop}$. If we choose as a
guiding principle the antipode of a Hopf algebroid $\hH$ to be a
bialgebroid map $\hH\to \hH^{op}_{cop}$, then $\hH$ and $\hH^{op}_{cop}$ need
to carry the same, say left, bialgebroid structure. This means that the
underlying 
algebra ${\bf H}$ must be equipped both with a left, and a right bialgebroid
structure. The first definition fulfilling this requirement was proposed in
\cite[Definition 4.1]{BohmSzl:hgdax}, where however the antipode was defined
to be bijective. Bijectivity of the antipode was relaxed in \cite[Definition
  2.2]{Bohm:hgdint}. Here we present a set of axioms which is equivalent to
\cite[Definition 2.2]{Bohm:hgdint}, as it was formulated in
\cite[Remark 2.1]{BohmBrz:hgdcleft}.
\begin{definition}\label{def:hgd}
For two algebras $\R$ and $\L$ over a commutative ring $k$,
a {\em Hopf algebroid} over the base algebras $\R$ and $\L$
is a triple $\hH=(\hH_L,\hH_R,S)$. Here $\hH_L$ is a left $\L$-bialgebroid and
$\hH_R$ is a right $\R$-bialgebroid, 
such that their underlying $k$-algebra ${\bf H}$ is the {\em same}. The {\em
  antipode} $S$ is a $k$-module map $H\to H$. 
Denote the $\R\sstac k \R^{op}$-ring structure of
$\hH_R$ by $({\bf H},s_R,t_R)$ and its $\R$-coring structure by
$(H,\Delta_R,\epsilon_R)$. Similarly, denote the $\L\sstac k \L^{op}$-ring
structure of $\hH_L$ by $({\bf H},s_L,t_L)$ and its $\L$-coring structure by
$(H,\Delta_L,\epsilon_L)$. 
Denote the multiplication in the $\R$-ring $(\bH,s_R)$ by $\mu_R$ and denote
  the multiplication in the $\L$-ring $(\bH,s_L)$ by $\mu_L$. 
These structures are subject to the following compatibility axioms.
\begin{itemize}
\item[(i)] $s_L\circ \epsilon_L \circ t_R=t_R$, $\quad t_L\circ \epsilon_L
  \circ s_R=s_R$, $\quad s_R\circ \epsilon_R \circ t_L=t_L$ and $t_R\circ
  \epsilon_R \circ s_L=s_L$.
\item[(ii)] $(\Delta_L\sstac \R H)\circ \Delta_R = (H\sstac \L \Delta_R)\circ
  \Delta_L$ and $(\Delta_R\sstac \L H)\circ \Delta_L = (H\sstac \R
  \Delta_L)\circ \Delta_R$.
\item[(iii)] For $l\in L$, $r\in R$ and $h\in H$, $S(t_L(l)ht_R(r))=s_R(r)
  S(h) s_L(l)$.
\item[(iv)] $\mu_L\circ (S\sstac \L H)\circ \Delta_L = s_R\circ \epsilon_R$
  and $\mu_R\circ (H\sstac \R S)\circ \Delta_R = s_L\circ \epsilon_L$.
\end{itemize}
\end{definition}
\begin{remarks}\label{rem:hgd.def}
Hopf algebroid axioms in Definition \ref{def:hgd} require some interpretation.
\begin{itemize}
\item[{(1)}] By the bialgebroid axioms, all maps $s_L\circ \epsilon_L$,
  $t_L\circ 
\epsilon_L$, $s_R\circ \epsilon_R$ and $t_R\circ \epsilon_R$ are idempotent
maps $H\to H$. Hence the message of axiom (i) is that the ranges of $s_L$ and
$t_R$, and also the ranges of $s_R$ and $t_L$, are coinciding subalgebras of
$\bH$. These axioms imply that the coproduct $\Delta_L$ in $\hH_L$ is not only
an $\L$-bimodule map, but also an $\R$-bimodule map. Symmetrically, $\Delta_R$
is an $\L$-bimodule map, so that axiom (ii) makes sense.
\item[{(2)}] The $k$-module $H$ underlying a (left or right) bialgebroid is a
  left and 
right comodule via the coproduct. Hence a $k$-module $H$ underlying a Hopf
algebroid is a left and right comodule for the constituent left and right
bialgebroids $\hH_L$ and $\hH_R$, via the two coproducts $\Delta_L$ and
$\Delta_R$. Axiom (ii) expresses  a property that these regular coactions do
commute, i.e. $H$ is an $\hH_L$-$\hH_R$ bicomodule and also an $\hH_R$-$\hH_L$
bicomodule. 

Alternatively, considering $H$ and $H\ot_\L H$ as right $\hH_R$-comodules via
the respective coactions $\Delta_R$ and $H\ot_\L \Delta_R$, the first axiom in
(ii) expresses that $\Delta_L$ is a right $\hH_R$-comodule map. Symmetrically,
this condition can be interpreted as a left $\hH_L$-comodule map property of
$\Delta_R$. Similarly, the second axiom in (ii) can be read as right
$\hH_L$-colinearity of $\Delta_R$ or left $\hH_R$-colinearity of $\Delta_L$.
\item[{(3)}] Axiom (iii) formulates the $\R$-$\L$ bimodule map property of the 
  antipode, needed in order for axiom (iv) to make sense.
\item[{(4)}] Analogously to the Hopf algebra axioms, axiom (iv) tells us that
  the antipode is convolution inverse of the identity map $H$, in some
  generalised sense. The notion of convolution products in the case of two
  different base 
algebras $\L$ and $\R$ is discussed in Section \ref{sec:hgd.cleft}.
\end{itemize}
\end{remarks}

Since in a Hopf algebroid $\hH$ there are two constituent bialgebroids $\hH_L$
and $\hH_R$ present, in these notes we use two versions of Sweedler's index
notation in parallel, to denote components of the coproducts $\Delta_L$ and
$\Delta_R$. We will use lower indices in the case of a left
bialgebroid $\hH_L$, i.e. we write $\Delta_L(h)=h_{(1)}\sstac \L h_{(2)}$, and
upper indices in the case of a right bialgebroid $\hH_R$, i.e. we write
$\Delta_R(h)=h^{(1)}\sstac \R h^{(2)}$, for $h\in H$, where implicit summation
is understood in both cases. Analogously, we use upper indices to denote
components of a coaction by $\hH_R$, and lower indices to denote components of
a coaction by $\hH_L$. 

\subsection{Examples and constructions}\label{sec:hgd_ex}
Before turning to a study of the structure of Hopf algebroids, let us see some
examples.

\subsubsection{Hopf algebras} 
A Hopf algebra $H$ over a commutative ring $k$ is an example of Hopf algebroids
over base algebras $\R=k=\L$. Both bialgebroids $\hH_L$ and $\hH_R$ are equal
to the $k$-bialgebra $H$ and the Hopf algebra antipode of $H$ satisfies the
Hopf algebroid axioms. 

Certainly, not every Hopf algebroid over base algebras $\R=k=\L$ is a
Hopf algebra (see e.g. Section \ref{sec:q.torus}). Examples of this kind have
been constructed in \cite{ConnMos:Cyc.cohom}, as follows. 
Let $H$ be a Hopf algebra over $k$, with coproduct 
$\Delta_L: h\mapsto h_{(1)}\sstac k h_{(2)}$, counit $\epsilon_L$ and antipode
$S$. 
Let $\chi$ be a character on $H$, i.e. a $k$-algebra map $H\to k$. The
coproduct 
$\Delta_R: h\mapsto h_{(1)}\sstac k \chi\left(S(h_{(2)})\right)h_{(3)}$ and
the counit $\epsilon_R:= \chi$ define a second bialgebra structure on the
$k$-algebra $\bH$. 
Looking at these two bialgebras as left and right $k$-bialgebroids
respectively, we obtain a Hopf algebroid with a {\em twisted}  
antipode $h\mapsto \chi(h_{(1)})S(h_{(2)})$. This construction was
extended in \cite[Theorem 4.2]{Bohm:alternative}, where new
examples of Hopf algebroids were constructed by twisting a (bijective)
antipode of a given Hopf algebroid. 

\subsubsection{Weak Hopf algebras}
A {\em weak Hopf algebra} over a commutative ring $k$ is a weak bialgebra
$H$ equipped with a $k$-linear antipode map $S:H\to H$, subject to the
following axioms \cite{BNSz:WHAI}. For $h\in H$,
$$
h_{(1)} S(h_{(2)}) =\sqcap^L(h),\qquad 
S(h_{(1)}) h_{(2)} =\sqcap^R(h),\qquad 
S(h_{(1)}) h_{(2)} S(h_{(3)}) = S(h),
$$
where the maps $\sqcap^L$ and $\sqcap^R$
were introduced in Section \ref{sec:WBA}.

The right $\R$-bialgebroid and the left $\R^{op}$-bialgebroid,
constructed for a weak Hopf algebra $H$ in Section \ref{sec:WBA}, together
with the antipode $S$, satisfy the Hopf algebroid axioms.

In particular, consider a small groupoid with finitely many objects. By
Section \ref{sec:WBA}, the free $k$-module spanned by its morphisms is a weak
$k$-bialgebra. It can be equipped with an antipode by putting $S(f):=f^{-1}$
for every morphism $f$, and extending it $k$-linearly. Motivated by this
example, weak Hopf algebras, and sometimes also Hopf algebroids, are called
{\em quantum groupoids} in the literature.

Weak Hopf algebras have a nice and well understood representation
theory. The category of finite dimensional modules of a finite dimensional
semisimple weak Hopf algebra $H$ over a field $k$ is a $k$-linear semisimple
category with finitely many inequivalent irreducible objects, with all finite
dimensional hom spaces. It is a monoidal category with left and right duals.
A category with the listed properties is termed a {\em fusion category}.
Conversely, based on Tannaka-Krein type reconstruction theorems in
\cite{Szlach:FinQG} and \cite{Hay:TK.rec}, it was proven in
\cite{EtiNikOst:Fusicat} that any fusion category is
equivalent to the category of finite dimensional modules of a (non-unique)
finite dimensional semisimple weak Hopf algebra. 
 
\subsubsection{$\R\otimes \R^{op}$} 
The left and right bialgebroids on an algebra of the form $\R\sstac k \R^{op}$,
constructed for any $k$-algebra $\R$ in Section \ref{sec:RotR-op}, form a Hopf
algebroid together with the antipode $r\sstac k r'\mapsto r'\sstac k r$.

\subsubsection{The algebraic quantum torus} \label{sec:q.torus}
Consider an algebra ${\bf T}_q$ over a commutative ring $k$, generated by two
invertible elements $U$ and $V$, subject to the relation $UV=qVU$, where $q$
is an invertible element in $k$. ${\bf T}_q$ possesses a right bialgebroid
structure over the commutative subalgebra $\R$ generated by $U$. 
Both the source and 
target maps are given by the inclusion $\R \to {\bf T}_q$. 
Coproduct and counit are defined by $\Delta_R:V^m U^n \mapsto V^m U^n\sstac \R
V^m$ and $\epsilon_R:V^m U^n \mapsto U^n$, respectively. 
Symmetrically, there is a left $\R$-bialgebroid structure given by the
coproduct $\Delta_L:U^n V^m \mapsto U^n V^m \sstac \R V^m$ and counit
$\epsilon_L:U^n V^m \mapsto U^n$.  Together with the antipode $S:U^n
V^m\mapsto V^{-m} U^n$ they constitute a Hopf algebroid. 

\subsubsection{Scalar extension}
Consider a Hopf algebra $H$ and a braided commutative algebra $A$ in the
category of right-right Yetter-Drinfel'd modules of $H$. As it was seen in
Section \ref{sec:scalar.ext}, the smash product algebra $A\# H$ carries a
right ${\bf A}$-bialgebroid structure. If the antipode $S$ of $H$ is
bijective, then the ${\bf A}$-bialgebroid structure of $A\# H$ extends to a
Hopf algebroid. Indeed, $A\# H$ is a left ${\bf A}^{op}$-bialgebroid via the
source map $a \mapsto a_{[0]}\cdot S(a_{[1]})\# 1_\B$, target map
$a \mapsto a_{[0]}\# a_{[1]}$, coproduct 
$
a\# h \mapsto (a\# h_{(1)})\sstac {{\bf A}^{op}} (1_{\bf A}\# h_{(2)})
$
and counit $a\# h \mapsto a_{[0]}\cdot S^{-1}\big(h S^{-1}(a_{[1]})\big)$. The
(bijective) antipode is given by
$a\# h \mapsto a_{[0]}\cdot S(h_{(2)})\# a_{[1]} S(h_{(1)})$.

\subsection{Basic properties of Hopf algebroids}
In this section some consequences of Definition \ref{def:hgd} of a Hopf
algebroid will be recalled from \cite[Section 2]{Bohm:hgdint}.

The opposite co-opposite of a Hopf algebra $H$ is a Hopf algebra, with the
same antipode $S$ of $H$. If $S$ is bijective, then also the opposite and the
co-opposite of $H$ are Hopf 
algebras, with antipode $S^{-1}$. A generalisation of these facts to Hopf
algebroids is given below. 
\begin{proposition}
Consider a Hopf algebroid $(\hH_L,\hH_R,S)$ over base 
algebras $\L$ and $\R$. The following hold true.
\begin{itemize}
\item[{(1)}] The triple $((\hH_R)^{op}_{cop}, (\hH_L)^{op}_{cop},S)$ is a Hopf
  algebroid over the base algebras $\R^{op}$ and $\L^{op}$.
\item[{(2)}] If the antipode $S$ is bijective then $((\hH_R)^{op},
  (\hH_L)^{op},S^{-1})$ is a Hopf algebroid over the base algebras $\R$ and
  $\L$, and $((\hH_L)_{cop}, (\hH_R)_{cop},S^{-1})$ is a Hopf algebroid over
  the base algebras $\L^{op}$ and $\R^{op}$.  
\end{itemize}
\end{proposition}
Proposition \ref{prop:antipode} states the expected behaviour of the
antipode of a Hopf algebroid with respect to the underlying ring and coring
structures. Consider a Hopf algebroid $\hH$
over base algebras $\L$ and
$\R$, with structure maps denoted as in Definition \ref{def:hgd}. It follows
immediately by axiom (i) in Definition \ref{def:hgd} that the
base algebras $\L$ and $\R$ are anti-isomorphic. Indeed, there are inverse  
algebra isomorphisms
\begin{equation}\label{eq:eps_L.s_R}
\epsilon_L\circ s_R:\R^{op}\to \L
\qquad \textrm{and}\qquad 
\epsilon_R\circ t_L: \L\to \R^{op}.
\end{equation}
Symmetrically, there are inverse algebra isomorphisms
\begin{equation}\label{eq:eps_R.s_L}
\epsilon_R\circ s_L:\L^{op}\to  \R
\qquad \textrm{and}\qquad 
\epsilon_L\circ t_R:\R\to \L^{op}.
\end{equation}
\begin{proposition}\label{prop:antipode}
Let $\hH$ be a Hopf algebroid over base 
algebras $\L$ and $\R$, with structure maps denoted as in Definition
\ref{def:hgd}. The following assertions hold. 
\begin{itemize}
\item[{(1)}] The antipode $S$ is a homomorphism of $\R\sstac k\R^{op}$-rings 
$$
(\ \bH\ ,\ s_R\ ,\ t_R\ )\to (\ \bH^{op}\ ,\  s_L\circ (\epsilon_L \circ s_R)\
  ,\  t_L\circ (\epsilon_L \circ s_R)\ )
$$ 
and also a homomorphism of $\L\sstac k\L^{op}$-rings 
$$
(\ \bH\ ,\ s_L\ ,\ t_L\ )\to (\ \bH^{op}\ ,\  s_R\circ (\epsilon_R \circ s_L)\
,\ t_R\circ (\epsilon_R \circ s_L)\ ).
$$
In particular, $S$ is a $k$-algebra anti-homomorphism $\bH\to \bH$.  
\item[{(2)}] The antipode $S$ is a homomorphism of $\R$-corings
$$
(\ H\ ,\ \Delta_R\ ,\ \epsilon_R\ )\to (\ H\ ,\ \Delta_L^{cop}\ ,\
  (\epsilon_R\circ s_L)   \circ \epsilon_L\ ),
$$ 
where $\Delta_L^{cop}$ is considered as a map $H\to H \sstac {\L^{op}} H \cong
H\sstac \R H$, via the isomorphism induced by the algebra isomorphism
\eqref{eq:eps_R.s_L}. Symmetrically, $S$ is a homomorphism of $\L$-corings  
$$
(\ H\ ,\ \Delta_L\ ,\ \epsilon_L\ )\to (\ H\ ,\  \Delta_R^{cop}\ ,\
(\epsilon_L\circ s_R) \circ  \epsilon_R\ ),
$$ 
where $\Delta_R^{cop}$ is considered as a map $H\to H \sstac {\R^{op}} H \cong
H\sstac \L H$, via the isomorphism induced by the algebra isomorphism
\eqref{eq:eps_L.s_R}.
\end{itemize}
\end{proposition}
A Hopf algebroid has a number of module structures over its base algebras $\L$
and $\R$. They turn out to be strongly related.

\begin{proposition}\label{prop:fgp.hgd} 
Let $\hH$ be a Hopf algebroid over base $k$-algebras $\L$ and $\R$, with
structure maps denoted as in Definition \ref{def:hgd}. 
If the antipde is bijective then the following statements hold.
\begin{itemize}
\item[{(1)}] The right $\L$-module $H$, given by left multiplication by 
$t_L$, is finitely generated and projective if and only if left $\R$-module
$H$, given by right multiplication by $t_R$, is finitely generated and
  projective. 
\item[{(2)}] The left $\L$-module $H$, given by left multiplication by 
$s_L$, is finitely generated and projective if and only if the right
$\R$-module $H$, given by right multiplication by  
$s_R$, is finitely generated and projective.
\end{itemize}
\end{proposition}

The $k$-dual $H^*$ of a finitely generated and projective Hopf algebra $H$
over a commutative ring $k$ is a Hopf algebra. The antipode in $H^*$ is the
transpose of the antipode of $H$. No generalisation of this fact for Hopf
algebroids 
is known. Although the dual $\hH^*$ of a finitely generated and projective Hopf
algebroid $\hH$ has a bialgebroid structure (cf. Section \ref{sec:dual.bgd}),
the transpose of the antipode in $\hH$ is not an endomorphism of $\hH^*$. Duals
only of {\em Frobenius} Hopf algebroids are known to be (Frobenius) Hopf
algebroids, see Section \ref{sec:Frob.hgd}.

\subsection{Comodules of Hopf algebroids}\label{sec:hgd.com}
In a Hopf algebroid, the constituent left and right bialgebroids are defined
on the {\em same} underlying algebra. Therefore, modules for the two
bialgebroids coincide. This is not the case with comodules: the two
bialgebroids have different underlying corings (over anti-isomorphic base
algebras, cf. \eqref{eq:eps_L.s_R}-\eqref{eq:eps_R.s_L}), that have a priori
different categories of (say, right) comodules. 
We take the opportunity to call here the reader's attention to a regrettable
error in the literature. Based on \cite[Theorem 2.6]{Brz:corext}, whose proof
turned out to be {\em incorrect}, the categories of (right) comodules of two
constituent bialgebroids in a Hopf algebroid were claimed to be strict
monoidally isomorphic in \cite[Theorem 2.2]{BohmBrz:hgdcleft}. Since it turned
out recently that in \cite[Theorem 2.6]{Brz:corext} there are some assumptions
missing, the derived result \cite[Theorem 2.2]{BohmBrz:hgdcleft} needs not hold
either at the stated level of generality. (There is a similar error in
\cite[Proposition 3.1]{Bohm:hgd.Gal}). Regrettably, this error influences
some results also in \cite{ArdBohmMen:Sch.type}, \cite{BohmBrz:hgdcleft}
and \cite{BohmVer:morita}. In the current section and in \ref{sec:hgd.Gal} we
present the corrected statements.

\subsubsection{Comodules of a Hopf algebroid and of its constituent
  bialgebroids} 
Since, as it is explained above, comodules of the two constituent bialgebroids
in a Hopf algebroid are in general different notions, none of them can be
expected to be a well working definition of a comodule of a Hopf
algebroid. The following definition of a comodule for a Hopf algebroid, as a
compatible comodule of both constituent bialgebroids, was suggested in
\cite[Definition 3.2]{Bohm:hgd.Gal} and \cite[Section 2.2]{BalSzl:fin.Gal}.

\begin{definition}\label{def:hgd_comod}
For a Hopf algebroid $\hH=(\hH_L,\hH_R,S)$ over base $k$-algebras $\L$ and $\R$,
denote the structure maps as in Definition \ref{def:hgd}. A {\em right
  $\hH$-comodule} is a right 
$\L$-module as well as a right $\R$-module 
$M$, together with a right
$\hH_R$-coaction $\varrho_R: M \to M\ot_\R H$ and a right $\hH_L$-coaction
$\varrho_L:M\to M\ot_\L H$, 
such that $\varrho_R$ is an $\hH_L$-comodule map and $\varrho_L$ is an
$\hH_R$-comodule map. Explicitly, 
$\varrho_R:m\mapsto m^{[0]}\ot_\R m^{[1]}$ is a right $\L$-module map in the
sense that
$$
(m\cdot l)^{[0]}\stac \R (m\cdot l)^{[1]} = m^{[0]} \stac \R t_L(l) m^{[1]},
\qquad \textrm{for } m\in M,\ l\in L,
$$
$\varrho_L:m\mapsto m_{[0]}\ot_\L m_{[1]}$ is a right $\R$-module map in the
sense that 
$$
(m\cdot r)_{[0]}\stac \L (m\cdot r)_{[1]} = m_{[0]} \stac \L m_{[1]} s_R(r),
\qquad \textrm{for } m\in M,\ r\in R,
$$
and the following compatibility conditions hold. 
$$
(\varrho_R \stac \L H)\circ \varrho_L = (M\stac \R \Delta_L)\circ \varrho_R
\qquad \textrm{and} \qquad 
(\varrho_L \stac \R H)\circ \varrho_R = (M\stac \L \Delta_R)\circ \varrho_L.
$$
Morphisms of $\hH$-comodules are $\hH_L$- as well as $\hH_R$-comodule maps.
\end{definition}

The category of right comodules of a Hopf algebroid $\hH$ is denoted by
$\cM^\hH$. 

It is not difficult to see that the right $\R$- and $\L$-actions on a right
$\hH$-comodule $M$ necessarily commute. That is, $M$ carries the structure of
a right $\L\ot_k \R$-module.

Note that by Definition \ref{def:hgd} (i) and (ii), the right $\R\ot_k
\L$-module $H$, with $\R$-action via the right source map $s_R$ and $\L$-action
via the left target map $t_L$, is a right comodule of the Hopf algebroid $\hH$,
via the two coactions given by the two coproducts $\Delta_R$ and $\Delta_L$.

Left comodules of a Hopf algebroid $\hH$ are defined symmetrically and their
category is denoted by $^\hH\cM$. 

\begin{remark}\label{rem:antip.functor}
The antipode $S$ in a Hopf algebroid $\hH$ defines a functor ${\mathcal
  M}^{\mathcal H}\to {}^\hH{\mathcal M}$. Indeed, if $M$ is a right
$\hH$-comodule, with $\hH_R$-coaction 
$m\mapsto m^{[0]}\ot_\R m^{[1]}$ and $\hH_L$-coaction $m\mapsto m_{[0]}\ot_\L
m_{[1]}$, then it is a left $\hH$-comodule with 
left $\R$-action $r\cdot m=m \cdot \epsilon_L(t_R(r))$, left $\L$-action
$l\cdot m=m\cdot \epsilon_R(t_L(l))$ (where the notations in Definition
\ref{def:hgd} are used) and respective coactions
$$
m\mapsto S(m_{[1]})\stac \R m_{[0]}
\qquad \textrm{and}\qquad 
m\mapsto S(m^{[1]})\stac \L m^{[0]}.
$$
Right $\hH$-comodule maps are also left $\hH$-comodule maps for
these coactions. 

A functor ${}^\hH{\mathcal M}\to {\mathcal M}^\hH$ is constructed
symmetrically. 
\end{remark}

Although comodules of a Hopf algebroid $\hH=(\hH_L,\hH_R,S)$ can not be
described as comodules of a coring, the free-forgetful adjunction
(cf. \cite[18.13(2)]{BrzWis:coring}), corresponding to the $L$-coring
underlying $\hH_L$, lifts to an adjunction between the categories $\cM^\hH$
and $\cM_\L$. Indeed, the forgetful functor $\cM^\hH \to \cM_\L$
has a right adjoint $-\ot_\L H: \cM_\L \to \cM^\hH$. Unit and counit of the
adjunction are given, for a right $\hH$-comodule $(M,\varrho_L,\varrho_R)$ and
a right $\L$-module ${\bf N}$, by the maps
\begin{equation}\label{eq:free.ff}
\varrho_L:M\to M\stac \L H  
\qquad \textrm{and}\qquad 
N\stac \L \epsilon_L : N\stac \L H \to N,
\end{equation}
respectively, where $\epsilon_L$ is the counit of $\hH_L$.
There is a similar adjunction between the cateories $\cM^\hH$ and $\cM_\R$.

Our next task is to look for situations when the category of comodules of a
Hopf algebroid coincides with the comodule category of any of the constituent
bialgebroids. 
Recall from Remark \ref{rem:hgd.def} (2) that for a Hopf algebroid
$\hH=(\hH_L,\hH_R,S)$, the underlying $k$-module $H$   
is an $\hH_L$-$\hH_R$ bicomodule and an $\hH_R$-$\hH_L$ bicomodule, via the
coactions given by the coproducts.  
In appropriate situations, taking {\em cotensor products} with a bicomodule
induces a functor between the comodule categories of the two corings, see
\cite[22.3]{BrzWis:coring} and its Erratum.  
In Theorem \ref{thm:pure_hgd} functors of this type are considered.

Recall from Section \ref{sec:bgd_comod} that any right comodule of a right
bialgebroid over a $k$-algebra $\R$ possesses a unique $\R$-bimodule
structure such that any comodule map is $\R$-bilinear. 
Thus if $\hH$ is a Hopf algebroid over the anti-isomorphic base $k$-algebras
$\L$ and $\R$, then any right comodule of the constituent right (or left)
bialgebroid can be regarded as a right $\L\ot_k \R$-module and the
coaction is a right $\R\ot_k \L$-module map.

\begin{theorem}\label{thm:pure_hgd}
Let $\hH=(\hH_L,\hH_R,S)$ be a Hopf algebroid, over base $k$-algebras $\L$
and $\R$, with structure maps denoted as in Definition \ref{def:hgd}. Consider
the $\R$- and $\L$-actions on $H$ that define its coring structures,
cf. Definitions \ref{def:right.bgd} and \ref{def:left.bgd}. If the equaliser
\begin{equation}\label{eq:M_eq}
\xymatrix{
M \ar[r]^-{\varrho} & 
M\stac \R H \ar@<2pt>[rr]^-{\varrho\ot_\R H} \ar@<-2pt>[rr]_-{M\ot_\R
  \Delta_R}&&  
M\stac \R H \stac \R H
}
\end{equation}
in the category ${\mathcal M}_\L$ of right $\L$-modules is preserved by both
functors $-\ot_\L H \ot_\L H$ and $- \ot_\L H \ot_\R H:{\mathcal M}_\L \to
{\mathcal M}_k$, for any right $\hH_R$-comodule $(M,\varrho)$, then the
forgetful functor 
$\cM^\hH \to \cM^{\hH_R}$ is an isomorphism.  
\end{theorem}
By standard terminology, the conditions in Theorem \ref{thm:pure_hgd} are
phrased as the equaliser \eqref{eq:M_eq} in ${\mathcal M}_\L$ is {\em $H\ot_\L
H$-pure and $H\ot_\R H$-pure}. 
Symmetrical conditions imply that the forgetful functor to the category of
right $\hH_L$-comodules is an isomorphism. 

Theorem \ref{thm:pure_hgd} is proven by constructing the inverse of the
forgetful functor, i.e. by equipping any right $\hH_R$-comodule with an
$\hH$-comodule structure. 
Any right $\hH_R$-comodule $M$ is isomorphic, as a
right $\R$-module, to the cotensor product $M \Box H$ of $M$ with the
$\hH_R$-$\hH_L$ bicomodule $H$. Since under the assumptions in Theorem
\ref{thm:pure_hgd}, $M\Box H$ is a right $\hH_L$-comodule
via $M\Box \Delta_L:M\Box H \to M\Box (H\ot_\L H)\cong (M\Box H)\ot_\L H$, the
isomorphism $M\cong M\Box H$ induces a right $\hH_L$-coaction on 
$M$. Moreover, by the assumptions in Theorem \ref{thm:pure_hgd}, 
we have isomorphisms 
$M \Box (H\ot_\L H) \cong M\ot_\L H$, 
$M\Box (H\ot_\R H)\cong M\ot_\R H$,
$M\Box (H \ot_\R H\ot_\L H)\cong M\ot_\R H\ot_\L H$ and
$M\Box (H \ot_\L H\ot_\R H)\cong M\ot_\L H\ot_\R H$.
The compatibility conditions between the $\hH_R$- and
$\hH_L$-coactions follow by these isomorphisms, the Hopf algebroid axioms
Definition \ref{def:hgd} (ii) and functoriality of the cotensor product. With
similar methods any $\hH_R$-comodule map is checked to be also
$\hH_L$-colinear. 

The purity conditions in Theorem \ref{thm:pure_hgd} are checked to hold in all
of the examples in Section \ref{sec:hgd_ex}. Moreover, if a Hopf
algebroid with a bijective antipode is finitely generated and projective in
all of the four senses in Proposition \ref{prop:fgp.hgd} then it
satisfies all purity conditions in Theorem \ref{thm:pure_hgd} since taking
tensor products with flat modules preserves any equaliser.

\subsubsection{Coinvariants in a comodule of a Hopf
  algebroid} \label{sec:coinv}  
By Definition \ref{def:hgd_comod}, a comodule $M$ of a Hopf algebroid 
$\hH=(\hH_L,\hH_R,S)$ is a comodule for both constituent bialgebroids $\hH_L$
and $\hH_R$. Since the unit element $1_H$ is grouplike for both corings
underlying $\hH_L$ and $\hH_R$, one can speak about coinvariants 
$$
M^{co\hH_R} =\{\ m\in M\ |\ \varrho_R(m) = m\stac \R 1_\bH\ \}
$$
of $M$ with respect to the $\hH_R$-coaction $\varrho_R$, or coinvariants 
$$
M^{co\hH_L} =\{\ m\in M\ |\ \varrho_L(m) = m\stac \L 1_\bH\ \}
$$
with respect to the $\hH_L$-coaction $\varrho_L$. Proposition
\ref{prop:hgd.coinv} relates these two notions. It is of crucial importance
from the point of view of Galois theory, see Section \ref{sec:hgd.Gal}.

\begin{proposition}\label{prop:hgd.coinv}
Let $\hH=(\hH_L,\hH_R,S)$ be a Hopf algebroid and $(M,\varrho_L,\varrho_R)$ be
a right $\hH$-comodule. Then any coinvariant of the $\hH_R$-comodule
$(M,\varrho_R)$ is coinvariant also for the $\hH_L$-comodule $(M,\varrho_L)$.

If moreover the antipode $S$ is bijective then coinvariants of the
$\hH_R$-comodule $(M,\varrho_R)$ and the $\hH_L$-comodule $(M,\varrho_L)$
coincide. 
\end{proposition}

For a right $\hH$-comodule $(M,\varrho_L,\varrho_R)$, consider the map 
\begin{equation}\label{eq:Phi}
\Phi_M:M\stac \R H \to M\stac \L H,\qquad m\ot_\R h \mapsto \varrho_L(m)\cdot
S(h), 
\end{equation}
where (using the notations in Definition \ref{def:hgd}) $H$ is a left
$\L$-module via the source map $s_L$ and a left $\R$-module 
via the target map $t_R$, and $M\ot_\L H$ is understood to be a right
$H$-module via the second factor. 
Since $\Phi_M(\varrho_R(m))=m\ot_\L 1_H$ and $\Phi_M(m\ot_\R 1_H)=\varrho_L(m)$,
we have the first claim in Proposition \ref{prop:hgd.coinv} proven. 
In order to prove the second assertion, note that if $S$ is an isomorphism
then so is $\Phi_M$, with inverse $\Phi_M^{-1}(m\ot_\L h) = S^{-1}(h)\cdot 
\varrho_R(m)$, where $M\ot_\R H$ is understood to be a left $H$-module via the
second factor.

\subsubsection{Comodule algebras of a Hopf algebroid.} \label{sec:hgd.com.alg}
As it was explained in Section \ref{sec:bgd.gal}, from the point of view of
Galois theory (in the coaction picture) monoidality of the category of
comodules is of central importance. Theorem \ref{thm:hgd.com.mon} replaces 
{\em unjustified} \cite[Theorem 2.2]{BohmBrz:hgdcleft} (cf. first paragraph of
Section \ref{sec:hgd.com}). 

By Definition \ref{def:hgd_comod} a right comodule of a Hopf algebroid $\hH$, 
over base $k$-algebras $\L$ and $\R$, is a right $\L\ot_k \R$-module. Since
$\L$ and $\R$ are anti-isomorphic algebras, we may regard, alternatively, any
$\hH$-comodule as an $\R$-bimodule by translating the right $\L$-action to a
left $\R$-action via the algebra anti-isomorphism \eqref{eq:eps_L.s_R}. 
 \begin{theorem}\label{thm:hgd.com.mon}
For a Hopf algebroid $(\hH_L,\hH_R,S)$ over base $k$-algebras $\L$ and $\R$, the
category ${\mathcal M}^{\mathcal H}$ of right $\hH$-comodules is
monoidal. Moreover, there are strict monoidal forgetful functors, rendering
commutative the following diagram. 
$$
\xymatrix{
{\mathcal M}^{\mathcal H} \ar[r]\ar[d]&
{\mathcal M}^{\hH_R} \ar[d]\\
{\mathcal M}^{\hH_L} \ar[r]&
{}_\R {\mathcal M}_\R.
}
$$
\end{theorem}
Commutativity of the diagram in Theorem \ref{thm:hgd.com.mon} follows by
comparing the unique $\R$-actions that make $\R$-bilinear the $\hH_R$-coaction
and the $\hH_L$-coaction in an $\hH$-comodule, respectively.
Strict monoidality of the functor on the right hand side was proven in Theorem
\ref{thm:comod.cat.mon}. Strict monoidality of the functor in the bottom row
follows by applying Theorem \ref{thm:comod.cat.mon} to the opposite of the
bialgebroid $\hH_L$ and identifying $\L^{op}$-bimodules and $\R$-bimodules via
the algebra isomorphism \eqref{eq:eps_L.s_R}. 
In order to see strict monoidality of the remaining two functors, recall that
by Theorem \ref{thm:comod.cat.mon} -- applied to $\hH_R$ and the opposite of
$\hH_L$ --, the $\R$-module tensor product of any two $\hH$-comodules is an
$\hH_R$-comodule and an $\hH_L$-comodule, via the diagonal coactions,
cf. \eqref{eq:diag.coac}. It is straightforward to check compatibility of
these coactions in the sense of Definition \ref{def:hgd_comod}. 
Similarly, $\R$($\cong \L^{op}$) is known to be
an $\hH_R$-comodule and an $\hH_L$-comodule, and compatibility of the coactions
is obvious. Finally, the $\R$-module tensor product of $\hH$-comodule maps is an
$\hH_R$-comodule map and an $\hH_L$-comodule map by Theorem
\ref{thm:comod.cat.mon}. Thus it is an $\hH$-comodule map. By Theorem
\ref{thm:comod.cat.mon} also the coherence natural transformations in ${}_\R
\cM_\R$ are $\hH_R$- and $\hH_L$-comodule maps, so $\hH$-comodule maps, what
proves Theorem \ref{thm:hgd.com.mon}.

Theorem \ref{thm:hgd.com.mon} enables us to introduce comodule algebras of
Hopf algebroids. 
\begin{definition}\label{def:hgd.comod.alg}
A {\em right comodule algebra} of a Hopf algebroid ${\mathcal H}$ is a
monoid in the monoidal category ${\mathcal M}^\hH$ of right
$\hH$-comodules. Explicitly, an $\R$-ring $(M,\mu,\eta)$, such that $M$ is a
right $\hH$-comodule and $\eta:R \to M$ and $\mu:M\ot_\R M \to M$ are right
$\hH$-comodule maps. Using the notations $m\mapsto m^{[0]}\ot_\R m^{[1]}$ and
$m\mapsto m_{[0]}\ot_\L m_{[1]}$ for the $\hH_L$- and $\hH_R$-coactions,
respectively, $\hH$-colinearity of $\eta$ and $\mu$ means the identities, for
all $m,m'\in M$, 
\begin{eqnarray*}
&{1_{\bf M}}^{[0]} \stac \R {1_{\bf M}}^{[1]} = 1_{\bf M}\stac \R 1_\bH,\qquad 
&(mm')^{[0]} \stac \R (mm')^{[1]} = m^{[0]} m^{\prime [0]} \stac \R m^{[1]}
  m^{\prime [1]}\\
&1_{{\bf M} {[0]}} \stac \L 1_{{\bf M} {[1]}} = 1_{\bf M}\stac \L 1_\bH,\qquad 
&(mm')_{[0]} \stac \L (mm')_{[1]} = m_{[0]} m'_{[0]} \stac \L m_{[1]}
  m'_{[1]} .
\end{eqnarray*}
Symmetrically, a {\em left $\hH$-comodule algebra} is a monoid in
${}^\hH\cM$. 
\end{definition}
The functors in Remark \ref{rem:antip.functor} induced by the antipode are
checked to be 
strictly anti-monoidal. Therefore, the opposite of a right $\hH$-comodule
algebra, with coactions in Remark \ref{rem:antip.functor}, is a left
${\mathcal H}$-comodule algebra and conversely. Thus there are four different
categories of modules of a comodule algebra of a Hopf algebroid.
\begin{definition}
Let $\hH$ be a Hopf algebroid and $M$ be a right $\hH$-comodule
algebra. Left and right $M$-modules in ${\mathcal M}^\hH$ are called
left-right and right-right {\em relative Hopf modules}, respectively. Their
categories are denoted by ${}_M{\mathcal M}^\hH$ and ${\mathcal M}_M^\hH$,
respectively. Left and right $M^{op}$-modules in ${}^\hH 
{\mathcal M}$ are called right-left and left-left relative Hopf modules,
respectively, and their categories are denoted by ${}^\hH{\mathcal M}_M$
and ${}^\hH_M{\mathcal M}$, respectively.

Explicitly, e.g. a right-right $(M,\hH)$-relative Hopf module is
a right module $W$ for the $\R$-ring $M$, such that the action is a right
$\hH$-comodule map $W\ot_\R M \to W$. Using index notations, with superscripts 
for the $\hH_R$-coactions and subscripts for the $\hH_L$-coactions, both on
$W$ and $M$, this means the identities, for $w\in W$ and $m\in M$,
$$
(w\cdot m)^{[0]}\stac \R (w\cdot m)^{[1]}= w^{[0]}\cdot m^{[0]}\stac \R
w^{[1]} m^{[1]}
\qquad \textrm{and} \qquad 
(w\cdot m)_{[0]}\stac \L (w\cdot m)_{[1]}= w_{[0]}\cdot m_{[0]}\stac \L
w_{[1]} m_{[1]} .
$$
\end{definition}
In contrast to relative Hopf modules of bialgebroids in Section
\ref{sec:ac.coac.pic}, relative Hopf modules of Hopf algebroids can not be
identified with comodules of a coring. Still, they determine an adjunction,
very similar to \eqref{eq:R.adj}.
Consider a right comodule algebra $M$ of a Hopf algebroid
$\hH=(\hH_L,\hH_R,S)$, over base algebras $\L$ and $\R$. Denote the
$\hH_R$-coinvariant subalgebra of $M$ by ${\bf N}$.
It follows from Proposition \ref{prop:hgd.coinv} that for any right ${\bf
N}$-module $V$, $V\ot_{\bf N} M$ is a right-right relative Hopf module via the
second 
factor. The resulting functor $-\ot_{\bf N} M: {\mathcal M}_{\bf N} \to
{\mathcal M}_M^\hH$ turns out to have a right adjoint: Any object $W$ in
$\cM_M^\hH$ can be regarded as an object in $\cM_M^{\hH_R}$, 
so we can take its $\hH_R$-coinvariants (cf. Section \ref{sec:coinv}). 
These considerations lead to an adjoint pair of functors 
\begin{equation}\label{eq:com.alg.adj}
- \stac {\bf N} M : {\mathcal M}_{\bf N} \to {\mathcal M}_M^\hH
\qquad \textrm{and}\qquad 
(-)^{co\hH_R} : {\mathcal M}_M^\hH \to {\mathcal M}_{\bf N}.
\end{equation}
The unit of the adjunction is given, for any right ${\bf N}$-module $V$, by
the map 
\begin{equation}\label{eq:adj.unit}
V \to (V\stac {\bf N} M)^{co\hH_R},\qquad 
v \mapsto v \stac {\bf N} 1_{\bf M}
\end{equation}
and the counit is given, for an $(M,\hH)$-relative Hopf module $W$, by 
\begin{equation}\label{eq:adj.counit}
W^{co\hH_R} \stac {\bf N} M \to M,\qquad  
w\stac {\bf N} m \mapsto w\cdot m.
\end{equation}
The message of this observation is that studying {\em descent
theory} of Galois extensions of a Hopf algebroid $\hH=(\hH_L,\hH_R,S)$, one
can examine both the adjunction in \eqref{eq:com.alg.adj}, corresponding
to an $\hH$-comodule algebra $M$, and also the adjunction
\eqref{eq:R.adj}, determined by $M$, regarded as an 
$\hH_R$-comodule algebra. Proposition
\ref{prop:two_adj} is obtained by observing that the units of the two
adjunctions coincide, and the counit of the adjunction in 
\eqref{eq:adj.counit} is obtained by restricting to the objects in ${\mathcal
M}^\hH_M$ the counit of the adjunction \eqref{eq:R.adj}. 
\begin{proposition}\label{prop:two_adj}
Consider a Hopf algebroid $\hH=(\hH_L,\hH_R,S)$ 
and a right $\hH$-comodule algebra $M$. Denote the $\hH_R$-coinvariant
subalgebra of $M$ by ${\bf N}$. 

(1) The  functor $-\ot_{\bf N} M:{\mathcal M}_{\bf N} \to {\mathcal
  M}^{\hH_R}_M$ is fully faithful if and only if the  functor $-\ot_{\bf
  N} M:{\mathcal M}_{\bf N} \to {\mathcal M}^\hH_M$ is fully faithful.

(2) If the functor $-\ot_{\bf N} M:{\mathcal M}_{\bf N} \to {\mathcal
  M}^{\hH_R}_M$ is an equivalence then also the  functor $-\ot_{\bf N}
M:{\mathcal M}_{\bf N} \to {\mathcal M}^\hH_M$ is an equivalence.
\end{proposition}

\subsubsection{The Fundamental Theorem of Hopf modules}\label{sec:fund.thm}
In this section we investigate the adjunction \eqref{eq:com.alg.adj} in a
special case.  

The coproducts $\Delta_L$ and $\Delta_R$ in a Hopf algebroid $\hH$ make the
underlying algebra $\bH$ a right $\hH$-comodule algebra. Corresponding
right-right relative Hopf modules are called simply Hopf modules and their
category 
is denoted by $\cM^H_H$. Coinvariants of the right $\hH_R$-comodule algebra
$H$ are the elements $t_R(r)$, for $r\in R$, where $t_R$ is the target map. 
If $\hH_R$ is the underlying bialgebroid in a Hopf algebroid $(\hH_L,\hH_R,S)$,
then $t_R:\R^{op}\to \bH$, equivalently, $s_L:\L \to \bH$, is a right
$\hH_R$-Galois extension, cf. Section \ref{sec:x_R.hopf.alg}. 
Hence Theorem \ref{thm:fund.thm}, known as the {\em Fundamental Theorem of 
Hopf modules}, can be interpreted as a Descent Theorem for this Galois
  extension $\L \cong \R^{op}\subseteq \bH$. 
\begin{theorem}\label{thm:fund.thm}
For a Hopf algebroid $\hH$,
over base algebras $\L$ and $\R$, the
functor $-\sstac {\L} H: \cM_{\L} \to \cM^H_H$ is an equivalence.
\end{theorem}
Theorem \ref{thm:fund.thm} is proven by constructing the inverses of the unit
\eqref{eq:adj.unit} and the counit \eqref{eq:adj.counit} of the relevant
adjunction. 
Use the notations for the structure maps of a Hopf algebroid in Definition
\ref{def:hgd}. 
For a right $\L$-module $V$, the inverse of \eqref{eq:adj.unit} is 
the map $(V\ot_\L H)^{co\hH_R} \to V$, $\sum_i v_i \ot_\L h_i \mapsto
\sum_i v_i \cdot s_L(\epsilon_L(h_i))$. 
For a Hopf module $W$, denote the $\hH_R$-coaction by $w\mapsto
w^{[0]}\ot_\R w^{[1]}$ and for the $\hH_L$-coaction write $w\mapsto
w_{[0]}\ot_\L w_{[1]}$. Then an epimorphism $W \to W^{co\hH_R}$ is
given by $w\mapsto w^{[0]} \cdot S(w^{[1]})$. The inverse of
\eqref{eq:adj.counit} is the map $W \to W^{co\hH_R}\ot_\L H$, $w\mapsto
{w_{[0]}}^{[0]} \cdot S({w_{[0]}}^{[1]}) \ot_\L w_{[1]}$.

\subsection{Integral theory}
In a Hopf algebra $H$, over a commutative ring $k$, integrals are invariants
of the regular module of the underlying $k$-algebra $\bH$, with respect to a
character given by the counit. The study of integrals provides a lot of
information about the structure of the $k$-algebra $\bH$. Most significantly, 
($k$-relative) semisimplicity of $\bH$ is equivalent to its separability over
$k$, and also to the existence of a normalised integral in $H$
\cite{CaeMil:Masc}, \cite{Lomp:int}. Since this extends Maschke's theorem
about group algebras, it is known as a Maschke type theorem. Its dual version
relates cosemisimplicity and coseparability of the $k$-coalgebra underlying
$H$ to the existence of normalised cointegrals \cite{CaeMil:Masc}. Another
group of results concerns Frobenius property of $\bH$, equivalent to the
existence of a non-degenerate integral in $H$ \cite{LaSw:int},
\cite{Par:Frob}. Integral theory of Hopf algebras was 
generalised to Hopf algebroids in \cite{Bohm:hgdint}. 
\begin{definition}\label{def:right.int}
For an algebra $\R$, consider a right
$\R$-bialgebroid $\cB$, with structure maps denoted as in Definition
\ref{def:right.bgd}. 
{\em Right integrals} in $\cB$ are the invariants of the right regular module
of the underlying $\R$-ring $(\B,s)$, with respect to the right character
$\epsilon$. Equivalently, invariants of the right regular module of
the $\R^{op}$-ring $(\B,t)$, with respect to $\epsilon$. That is, the elements
of   
$$
B_\epsilon=\{\ i \in B \ |\  i b = i s(\epsilon (b)),\quad \forall b\in B\ \}
=\{\ i \in B \ |\  i b = i t(\epsilon (b)),\quad \forall b\in B\ \}
\cong \mathrm{Hom}_\B(R,B).
$$
A right integral $i$ is {\em normalised} if $\epsilon(i)=1_\R$.
\end{definition}
Symmetrically, {\em left integrals} in a left $\R$-bialgebroid are defined as
invariants of the left regular module of the 
underlying $\R$-ring or $\R^{op}$-ring, with respect to a left character
defined by the counit. 
Note that a left integral $i$ in a left bialgebroid $\cB$ is a left integral
also in $\cB_{cop}$ and a right integral in the right bialgebroids $\cB^{op}$
and $\cB^{op}_{cop}$.  

Since the counit in a right (resp. left) bialgebroid is a right (resp. left)
character, there is no way to consider left (resp. right) integrals in a right
(resp. left) bialgebroid. On the contrary, the base algebra in a right
bialgebroid $\cB$ is both a right and a left $\cB$-comodule, via coactions
given by the source and target maps, respectively. Hence there are
corresponding notions of left and right cointegrals. 
\begin{definition}
For an algebra $\R$, consider a right
$\R$-bialgebroid $\cB$, with structure maps denoted as in Definition
\ref{def:right.bgd}.
A {\em right cointegral} on $\cB$ is an element of 
$$
\mathrm{Hom}^B(B,R)=\{\ \iota\in \mathrm{Hom}_\R(B,R)\ |\ (\iota\stac \R
B)\circ \Delta = s\circ \iota \}.
$$
A right cointegral $\iota$ is {\em normalised} if $\iota(1_\B)=1_\R$.

Symmetrically, a {\em left cointegral} on $\cB$ is an element of
${}^B\mathrm{Hom}(B,R)$. 
\end{definition}
Left and right cointegrals on a left $\R$-bialgebroid $\cB$ are defined
analogously, as left and right comodule maps $B\to R$.  

If a right $\R$-bialgebroid $\cB$ is finitely generated and projective as a,
say right, $\R$-module (via right multiplication by the source map), then the
isomorphism $\mathrm{Hom}_\R(R,B^*)\cong \mathrm{Hom}_\R (B,R)$ induces an
isomorphism ${}_{B^*}\mathrm{Hom}(R,B^*)\cong \mathrm{Hom}^B(B,R)$. Hence in
this case left integrals in the left $\R$-bialgebroid $\cB^*$ are the same as
right cointegrals on $\cB$. Similar statements hold for all other duals of left
and right bialgebroids.

For a Hopf algebroid $(\hH_L,\hH_R,S)$, left (resp. right) cointegrals on
$\hH_L$ and $\hH_R$ can be shown to be left (resp. right) $\hH$-comodule maps.

\subsubsection{Maschke type theorems}
Recall that an $\R$-ring $B$ is said to be {\em separable}
provided that the multiplication map $B\sstac \R B \to B$ is a split
epimorphism of $B$-bimodules. The $\R$-ring $B$ is said to be left
(resp. right) {\em semisimple} (or sometimes $\R$-relatively semisimple) if
every left (resp. right) $B$-module is $\R$-relative projective. That is,
every $B$-module epimorphism, which has an $\R$-module section, is a split
epimorphism of $B$-modules. By a classical result due to Hirata and  Sugano
\cite{HirSug:sep}, 
a separable $\R$-ring is left, and right semisimple. For Hopf algebroids also
the converse can be proven.
\begin{theorem}\label{thm:Maschke}
For a Hopf algebroid $(\hH_L,\hH_R,S)$ over base 
algebras $\L$ and $\R$, denote the structure maps as in Definition
\ref{def:hgd}. 
The following properties are equivalent.
\begin{itemize}
\item[{(i)}] The $\R$-ring $(\bH,s_R)$ underlying $\hH_R$ is separable.
\item[{(ii)}] The $\R^{op}$-ring $(\bH,t_R)$ underlying $\hH_R$ is separable.
\item[{(iii)}] The $\L$-ring $(\bH,s_L)$ underlying $\hH_L$ is separable.
\item[{(iv)}] The $\L^{op}$-ring $(\bH,t_L)$ underlying $\hH_L$ is separable.
\item[{(v)}] The $\R$-ring $(\bH,s_R)$ underlying $\hH_R$ is right semisimple.
\item[{(vi)}] The $\R^{op}$-ring $(\bH,t_R)$ underlying $\hH_R$ is right
  semisimple. 
\item[{(vii)}] The $\L$-ring $(\bH,s_L)$ underlying $\hH_L$ is left semisimple.
\item[{(viii)}] The $\L^{op}$-ring $(\bH,t_L)$ underlying $\hH_L$ is left
  semisimple. 
\item[{(ix)}] There exists a normalised right integral in $\hH_R$.
\item[{(x)}] There exists a normalised left integral in $\hH_L$.
\item[{(xi)}] The counit $\epsilon_R$ in $\hH_R$ is a split epimorphism of
  right $\bH$-modules. 
\item[{(xii)}] The counit $\epsilon_L$ in $\hH_L$ is a split epimorphism of
  left $\bH$-modules. 
\end{itemize}
\end{theorem}
Since the source map in a right $\R$-bialgebroid is a right $\R$-module
section of the counit, implication (v)$\Rightarrow$(xi) is obvious. For a
right $\bH$-module section $\nu$ of the counit, $\nu(1_\R)$ is a normalised
right integral. Thus (xi)$\Rightarrow$(ix). The antipode in a Hopf algebroid
maps a normalised right integral in $\hH_R$ to a normalised left integral in
$\hH_L$, and vice versa. So (ix)$\Leftrightarrow$(x). If $i$ is a
normalised left integral in $\hH_L$, then the map $H\to H\sstac \R H$,
$h\mapsto h i^{(1)}\sstac \R S(i^{(2)})=i^{(1)}\sstac \R S(i^{(2)})h$ is an
$\bH$-bimodule section of the multiplication in the $\R$-ring underlying
$\hH_R$ (where the index notation $\Delta_R(h)=h^{(1)}\ot_\R h^{(2)}$ is used,
for $h\in H$.). This proves
(x)$\Rightarrow$(i). The remaining equivalences follow by symmetry. 
Note that equivalences (iv)$\Leftrightarrow$(viii)$\Leftrightarrow$(x)
$\Leftrightarrow$(xii) hold also for a $\times_L$-Hopf algebra $\hH_L$
(discussed Section \ref{sec:x_R.hopf.alg}). 

As an alternative of Theorem \ref{thm:Maschke}, one can ask about properties
of the $\R\sstac k \R^{op}$-ring, underlying a right bialgebroid $\hH_R$, and
the $\L\sstac k \L^{op}$-ring, underlying a left bialgebroid $\hH_L$, in a
Hopf algebroid $(\hH_L,\hH_R,S)$. Theorem \ref{thm:alt.Maschke} is
obtained by application of \cite[Theorem 6.5]{BrugVire:Hopf.monad}. 
For a $k$-algebra $\L$, consider a left $\L$-bialgebroid $\hH_L$. Denote its
$\L\sstac k \L^{op}$-ring structure by $(\bH,s_L,t_L)$ and its $\L$-coring
structure by $(H,\Delta_L,\epsilon_L)$. Look at $L$ as a left $\L\sstac k
\L^{op}$-module, with action given by left and right multiplications. 
Look at
$H$ as a right $\L\sstac k \L^{op}$-module, with action given by right
multiplications by $s_L$ and $t_L$. Note that 
\begin{equation}\label{eq:quotient.coring}
H\sstac{\L\sstac k \L^{op}} L \cong H/\{\ hs_L(l) - ht_L(l)\ |\ h\in H, \ l\in
L\ \} 
\end{equation}
is an $\L$-coring (via quotient maps of $\Delta_L$ and $\epsilon_L$) and a
left $\bH$-module. Hence we can speak about the invariants of
$H\sstac{\L\sstac k \L^{op}} L$ with respect to $\epsilon_L$. An 
invariant of $H\sstac{\L\sstac k\L^{op}} L$
is said to be {\em normalised} if the quotient of
$\epsilon_L$ maps it to $1_\L$. 
\begin{theorem}\label{thm:alt.Maschke}
Consider a Hopf algebroid $(\hH_L,\hH_R,S)$, over base $k$-algebras $\L$ and
$\R$. The following assertions are equivalent. 
\begin{itemize}
\item[{(i)}] The $\R\sstac k \R^{op}$-ring underlying $\hH_R$ is separable.
\item[{(ii)}] The $\L\sstac k \L^{op}$-ring underlying $\hH_L$ is separable.
\item[{(iii)}] The $\R\sstac k \R^{op}$-ring underlying $\hH_R$ is right
  semisimple. 
\item[{(iv)}] The $\L\sstac k \L^{op}$-ring underlying $\hH_L$ is left
  semisimple. 
\item[{(v)}] There is a normalised invariant in the right $\bH$-module
  $R\sstac{\R\sstac k \R^{op}} H$. 
\item[{(vi)}] There is a normalised invariant in the left $\bH$-module
  $H \sstac{\L\sstac k \L^{op}} L$. 
\end{itemize}
\end{theorem}
For a $\times_L$-Hopf algebra $\hH_L$ (discussed in Section
\ref{sec:x_R.hopf.alg}), equivalences
(ii)$\Leftrightarrow$(iv)$\Leftrightarrow$(vi) in Theorem
\ref{thm:alt.Maschke} hold true.  

Recall that an $\R$-coring $B$ is said to be {\em coseparable}
provided that the comultiplication map $B\to B\sstac \R B$ is a split
monomorphism of $B$-bicomodules. The $\R$-coring $B$ is said to be left
(resp. right) {\em cosemisimple} (or sometimes $\R$-relatively cosemisimple)
if every left (resp. right) $B$-comodule is $\R$-relative injective. That is,
every $B$-comodule monomorphism, which has an $\R$-module retraction, is a
split monomorphism of $B$-comodules. A coseparable $\R$-coring is left, and
right cosemisimple. For Hopf algebroids also the converse can be proven.
\begin{theorem}\label{thm:dual.Maschke}
For a Hopf algebroid $(\hH_L,\hH_R,S)$ over base 
algebras $\L$ and $\R$, the following properties are equivalent.
\begin{itemize}
\item[{(i)}] The $\R$-coring underlying $\hH_R$ is coseparable.
\item[{(ii)}] The $\L$-coring underlying $\hH_L$ is coseparable.
\item[{(iii)}] The $\R$-coring underlying $\hH_R$ is right cosemisimple.
\item[{(iv)}] The $\R$-coring underlying $\hH_R$ is left cosemisimple.
\item[{(v)}] The $\L$-coring underlying $\hH_L$ is right cosemisimple.
\item[{(vi)}] The $\L$-coring underlying $\hH_L$ is left cosemisimple.
\item[{(vii)}] There exists a normalised right cointegral on $\hH_R$.
\item[{(viii)}] There exists a normalised left cointegral on $\hH_R$.
\item[{(ix)}] There exists a normalised right cointegral on $\hH_L$.
\item[{(x)}] There exists a normalised left cointegral on $\hH_L$.
\item[{(xi)}] The source map in $\hH_R$ is a split right $\hH_R$-comodule
  monomorphism.  
\item[{(xii)}] The target map in $\hH_R$ is a split left $\hH_R$-comodule
  monomorphism.  
\item[{(xiii)}] The source map in $\hH_L$ is a split left $\hH_L$-comodule
  monomorphism.  
\item[{(xiv)}] The target map in $\hH_L$ is a split right $\hH_L$-comodule
  monomorphism.  
\end{itemize}
\end{theorem}

\subsubsection{Frobenius Hopf algebroids} \label{sec:Frob.hgd}
It was proven by Larson and Sweedler in \cite{LaSw:int} that every finite
dimensional Hopf 
algebra over a field is a Frobenius algebra. Although this is not believed to
be true for any finitely generated projective Hopf algebra over a commutative
ring, Frobenius Hopf algebras form a distinguished class. A Hopf algebra is
known to be a Frobenius algebra if and only if it possesses a non-degenerate
integral \cite{Par:Frob}. It is a self-dual property: a non-degenerate
integral determines a non-degenerate cointegral, i.e. a non-degenerate
integral in the dual Hopf algebra. 

In a Hopf algebroid there are four algebra extensions present: the ones given
by  the source and target maps of the two constituent bialgebroids. Among Hopf
algebroids, those in which these are Frobenius extensions, play an even more
distinguished role. Although the dual of any finitely generated projective Hopf
algebroid is not known to be a Hopf algebroid, duals of Frobenius Hopf
algebroids are Frobenius Hopf algebroids. 

While every finitely generated and projective Hopf algebra over a commutative
ring $k$ was proven by Pareigis to be a quasi-Frobenius $k$-algebra in
\cite{Par:Frob}, an analogous statement 
fails to hold for Hopf algebroids. In \cite[Section 5]{Bohm:hgdint} Hopf
algebroids were constructed, which are finitely generated and projective over
their base algebras (in all the four senses in Proposition \ref{prop:fgp.hgd}),
but are not quasi-Frobenius extensions of the base algebra. 

Recall (e.g. from \cite{Kad:Frob}) that an $\R$-ring $(\bH,s)$ is said to be
{\em Frobenius} provided that 
$H$ is a finitely generated and projective left $\R$-module and
${}^*H:={}_\R\mathrm{Hom}(H,R)$ is isomorphic to $H$ as an $\bH$-$\R$
bimodule. Equivalently, if $H$ is a finitely generated and projective right
$\R$-module and $H^*:=\mathrm{Hom}_\R(H,R)$ is isomorphic to $H$ as an
$\R$-$\bH$ bimodule. These properties are equivalent also to the existence of
an $\R$-bimodule map $\psi:\bH\to \R$, the so called Frobenius functional,
possessing a dual basis $\sum_i e_i\sstac \R f_i\in H\sstac \R H$, satisfying,
for all $h\in H$, $\sum_i e_i\cdot \psi(f_i h)=h = \sum_i\psi(h e_i) \cdot f_i$.
The following characterisation of Frobenius Hopf algebroids was obtained in
\cite[Theorem 4.7]{Bohm:hgdint}, see the {\em Corrigendum}.
\begin{theorem}\label{thm:frob.hgd}
Consider a Hopf algebroid $\hH$, over base 
algebras $\L$ and $\R$, with structure maps denoted as in Definition
\ref{def:hgd}. 
Assume that $H$ is finitely generated and projective 
as a right $R$-module via right multiplication by $s_R$,  
as a left $R$-module via right multiplication by $t_R$, 
as a left $L$-module via left multiplication by $s_L$ and 
as a right $L$-module via left multiplication by $t_L$.
The following statements are equivalent.
\begin{itemize}
\item[{(i)}] The $\R$-ring $(\bH,s_R)$ is Frobenius. 
\item[{(ii)}] The $\R^{op}$-ring $(\bH^{op},t_R)$ is Frobenius. 
\item[{(ii)}] The $\L$-ring $(\bH,s_L)$ is Frobenius. 
\item[{(iv)}] The $\L^{op}$-ring $(\bH^{op},t_L)$ is Frobenius. 
\item[{(v)}] There exists a right cointegral
  $\iota$ on $\hH_R$, such that the map ${\widetilde\iota}:H\to
  \mathrm{Hom}_R(H,R)$, $h\mapsto \iota(h-)$ is bijective.  
\item[{(vi)}] There exists a left cointegral 
  $\upsilon$ on $\hH_L$, such that the map ${\widetilde\upsilon}:H\to
  {}_\L\mathrm{Hom}(H,L)$, $h\mapsto \upsilon(-h)$ is bijective.  
\item[{(vii)}] There exists a right integral $i$ in $\hH_R$, such that the maps
  ${\widetilde i}:{}_\L\mathrm{Hom}(H,L)\to H$, $\psi\mapsto
  t_L(\psi(i_{(2)}))i_{(1)}$ 
and 
is   bijective.  
\item[{(viii)}] There exists a left integral $j$ in $\hH_L$, such that the map
  ${\widetilde j}:\mathrm{Hom}_\R(H,R)\to H$, $\phi\mapsto
  j^{(2)}t_R(\phi(j^{(1)}))$ is bijective.  
\end{itemize}
A Hopf algebroid for which these equivalent conditions hold is said to be {\em
  Frobenius}, and a right (resp. left) integral obeying property (vii)
  (resp. (viii)) is said to be {\em non-degenerate}. 

In a Frobenius Hopf algebroid the antipode $S$ is  bijective. 
\end{theorem}
If $j$ is a non-degenerate left integral in $\hH_L$
then $\iota:=(\,{\widetilde j}\,)^{-1}(1_\bH)$
is a right cointegral on $\hH_R$ and a Frobenius functional for the $\R$-ring
$(\bH,s_R)$ with dual basis $j^{(1)}\sstac \R S(j^{(2)})$. Thus
(viii)$\Rightarrow$(i). If property (i) holds, then the right cointegrals on
$\hH_R$ are shown to form a free rank one left $\R$-module ${\mathcal I}$, via
the action $r\cdot \iota=\iota(t_R(r)-)$. Using finitely generated
projectivity of the right $\R$-module $H$, the dual
$H^*:=\mathrm{Hom}_\R(H,R)$ can be equipped with a Hopf module structure, with
coinvariants ${\mathcal I}$. 
Hence Theorem \ref{thm:fund.thm} implies an isomorphism $H^*\cong H\sstac \R
{\mathcal I}$. 
This isomorphism is used to show that the cyclic generator $\iota$ of the
$\R$-module ${\mathcal I}$ satisfies condition (v). If there is a right
cointegral $\iota$ as in part (v), then a non-degenerate left integral $j$ as
in part (viii) is constructed in terms of $(\,{\widetilde{\iota}}\,)^{-1}$ and 
a dual basis for the finitely generated projective right $\R$-module $H$. It
is shown to satisfy $(\,{\widetilde j}\,)^{-1}={\widetilde   \iota\circ S}$,
which implies bijectivity 
of $S$. The remaining equivalences follow by relations between the source and
target maps in $\hH_L$ and $\hH_R$, and symmetrical versions of the arguments
above.    

For a Frobenius Hopf algebroid $(\hH_L,\hH_R,S)$, over base algebras $\L$ and
$\R$, all the four duals $\mathrm{Hom}_\R(H,R)$, ${}_\R \mathrm{Hom}(H,R)$,
$\mathrm{Hom}_\L(H,L)$ and ${}_\L\mathrm{Hom}(H,L)$ possess (left or right)
bialgebroid structures. A left integral $j$ in $\hH_L$, such that the map
${\widetilde j}$  
in part (viii) of Theorem \ref{thm:frob.hgd} is bijective, determines  
further similar isomorphisms ${}_\R \mathrm{Hom}(H,R) \to H$,
$\mathrm{Hom}_\L(H,L)\to H$ and ${}_\L\mathrm{Hom}(H,L)\to H$. What is more,
putting $\iota:=(\,{\widetilde j}\,)^{-1}(1_\bH)$, there is an algebra
automorphism of $H$,
\begin{equation}\label{eq:zeta}
\zeta:H\to H,\qquad  h\mapsto h^{(2)} t_R(\iota(j\,h^{(1)})).
\end{equation}
These isomorphisms 
combine to bialgebroid (anti-) isomorphisms between the four
duals of $H$, cf. \cite[Theorem 5.16]{BohmSzl:hgdax}.

Frobenius property of a Hopf algebroid was shown to be self-dual in
\cite[Theorem 5.17 and Proposition 5.19]{BohmSzl:hgdax}, in the following
sense. 
\begin{theorem}
Consider a Frobenius Hopf algebroid $(\hH_L,\hH_R,S)$, over base algebras $\L$
and $\R$.
Let $j$ be a left integral in $\hH_L$, such that the map ${\widetilde j}$ in
part (viii) of Theorem \ref{thm:frob.hgd} is bijective. 
Then the left $\R$-bialgebroid $H^*:=\mathrm{Hom}_\R(H,R)$ extends to a
Hopf algebroid. A bijective antipode is given in terms of the map
\eqref{eq:zeta}, by $S^*:=(\,{\widetilde j}\,)^{-1}\circ S\circ \zeta\circ
{\widetilde j}$. The right bialgebroid structure is determined by the 
requirement that $S^*$ is a bialgebroid anti-isomorphism in the sense of
Proposition \ref{prop:antipode}. This dual Hopf algebroid is Frobenius, with
non-degenerate left integral $(\,{\widetilde j}\,)^{-1}(1_\bH)\in\hH^*$.
\end{theorem}
In a paper \cite{Szlach:double.alg} by Szlach\'anyi, an equivalent description
of a Frobenius Hopf algebroid $(\hH_L,\hH_R,S)$ was proposed, via so called
{\em double algebras}. In this picture the isomorphism ${\widetilde j}$ in part
(viii) of Theorem \ref{thm:frob.hgd} is used to transfer the multiplication in 
$H^*:=\mathrm{Hom}_\R(H,R)$ to a second algebra structure in $H$, with unit
$j$. In this way four Frobenius ring structures are obtained on $H$. Note that
the coproducts in $\hH_L$ and $\hH_R$ correspond canonically to the Frobenius
ring structures transferred from $H^*$. In this approach the Hopf algebroid
axioms are formulated as compatibility conditions between the two algebra
structures on $H$.

\subsection{Galois theory of Hopf algebroids}\label{sec:hgd.Gal}
Galois extensions by Hopf algebras are the same as Galois extensions by the
constituent bialgebra. Still, since the structure of a Hopf algebra is more
complex than that of a bialgebra, it allows to derive stronger results.
On the contrary, Galois theory of a Hopf algebroid is more conceptually
different from Galois theory of the constituent bialgebroids: As it is
discussed in Section \ref{sec:hgd.com}, comodules of a Hopf algebroid carry
more structure than comodules of any constituent bialgebroid. Consequently,
Galois theory of Hopf algebroids, discussed in the current section, is
significantly richer than the theory of bialgebroids. In particular, for a
comodule algebra of a Hopf algebroid (cf. Section \ref{sec:hgd.com.alg}),
several theorems concerning an equivalence between the 
category of relative Hopf modules and the category of modules of the
coinvariant subalgebra -- i.e. descent theorems -- can be proven.

By Definition \ref{def:hgd.comod.alg} and Theorem \ref{thm:hgd.com.mon}, a
right comodule algebra $M$ of a Hopf algebroid $\hH=(\hH_L,\hH_R,S)$ is
both an $\hH_L$-comodule algebra and an $\hH_R$-comodule algebra.  
Denote the $\hH_R$-coinvariant subalgebra of $M$ by ${\bf N}$.
In light of Proposition \ref{prop:hgd.coinv}, there are two corresponding
canonical maps  
\begin{eqnarray}\label{eq:two.canR}
&M\stac {\bf N} M \to M\stac \R H,\qquad &m\stac {\bf N} m'\mapsto m m^{\prime 
  [0]}\stac \R m^{\prime  [1]}\qquad \textrm{and}\\
&M\stac {\bf N} M \to M\stac \L H,\qquad &m\stac {\bf N} m'\mapsto m_{[0]} m'
  \stac \L m_{[1]},\label{eq:two.canL}
\end{eqnarray}
where $m\mapsto m^{[0]}\ot_R m^{[1]}$ and $m\mapsto m_{[0]}\ot_L m_{[1]}$
denote the $\hH_R$-coaction and the $\hH_L$-coaction on $M$, respectively.
In general, bijectivity of the two canonical maps \eqref{eq:two.canR} and
\eqref{eq:two.canL}, 
are not known to be equivalent. Only a partial result \cite[Lemma
  3.3]{Bohm:hgd.Gal} is known. 
\begin{proposition}\label{prop:bij.antip.Gal}
If the antipode $S$ in a Hopf algebroid $(\hH_L,\hH_R,S)$ is bijective, then
the 
$\hH_R$-canonical map \eqref{eq:two.canR} is bijective if and only if the
$\hH_L$-canonical map \eqref{eq:two.canL} is bijective.
\end{proposition}

This follows by noting that the two canonical maps differ by the isomorphism
$\Phi_A$ in \eqref{eq:Phi}.  

By Proposition \ref{prop:hgd.coinv} and Proposition \ref{prop:bij.antip.Gal},
a for right comodule algebra $M$ of a Hopf algebroid $(\hH_L,\hH_R,S)$ with a
bijective antipode, an
algebra extension ${\bf N}\subseteq {\bf M}$ is $\hH_L$-Galois if and only if
it is $\hH_R$-Galois.
\begin{remark}\label{rem:fgp.hgd}
Consider a Hopf algebroid $\hH=(\hH_L,\hH_R,S)$ over base algebras $\L$ and
$\R$, which is finitely generated and projective in the equivalent senses in
Proposition \ref{prop:fgp.hgd} (2). Then $H$ is in particular flat as a left
$\L$-module and as a left $\R$-module. Therefore, by Theorem
\ref{thm:pure_hgd} and Theorem \ref{thm:hgd.com.mon}, the category of right
$\hH$-comodules is strict monoidally isomorphic to the category of comodules
for any of the constituent left and right bialgebroids $\hH_L$ and $\hH_R$.  
Thus comodule algebras for $\hH$ coincide with comodule algebras for $\hH_L$
or $\hH_R$.

If furthermore the antipode is bijective, we
conclude by Proposition \ref{prop:hgd.coinv} and Proposition 
\ref{prop:bij.antip.Gal} that 
an algebra extension ${\bf N}\subseteq {\bf M}$
is a right $\hH_R$-Galois extension if and only if it is a right
$\hH_L$-Galois extension. 
\end{remark}

\subsubsection{Depth two Frobenius extensions}
An analogue of Theorem \ref{thm:bgd.Gal.ext} for Frobenius Hopf algebroids is
\cite[Theorem 3.6]{BalSzl:fin.Gal}. 
\begin{theorem}\label{thm:Frob.hgd.Gal.ext}
An algebra extension ${\bf N}\subseteq {\bf M}$ is a right Galois extension by
some 
Frobenius Hopf algebroid (i.e. by any of its constituent bialgebroids) if and
only if it is a Frobenius extension, it is balanced, and satisfies the (left
and right) depth 2 conditions.  
\end{theorem}
In a case of a Frobenius extension ${\bf N}\subseteq {\bf M}$, left and right
depth 2 properties are equivalent. 
By Remark \ref{rem:fgp.hgd}, for a Frobenius Hopf algebroid
$(\hH_L,\hH_R,S)$, $\hH_L$- and $\hH_R$-Galois properties of an extension are
equivalent. 
By Theorem \ref{thm:Gal.bgds}, a right
depth 2 and balanced  algebra extension ${\bf N}\subseteq {\bf M}$
is a Galois extension by a right bialgebroid $(M\sstac {\bf N} M)^{\bf N}$. 
If in addition ${\bf N}\subseteq {\bf M}$ is a  Frobenius extension then
$(M\sstac {\bf N} M)^{\bf N}$ is shown to be a constituent right bialgebroid
in a Frobenius Hopf algebroid. A non-degenerate (left and right) integral 
$\sum_i m_i\sstac {\bf N} m'_i\in (M\sstac {\bf N} M)^{\bf N}$ is provided by
the dual basis of a Frobenius functional $\psi:M\to N$. A non-degenerate (left
and right) integral in the dual Hopf algebroid ${}_{\bf N} \mathrm{End}_{\bf
  N}(M)$ is  
$\psi$. In the converse direction, note that a right comodule $M$ of a Hopf
algebroid $\hH$ is an $\hH^*$-module. If $\hH$ is a Frobenius Hopf algebroid
and ${\bf N}\subseteq {\bf M}$ is an $\hH$-Galois extension, a Frobenius
functional $M\to N$ is given by the action by a non-degenerate integral in
$\hH^*$. 

\subsubsection{Cleft extensions by Hopf algebroids}\label{sec:hgd.cleft}
For an algebra ${\bf M}$ and a coalgebra $C$ over a commutative ring $k$,
$\mathrm{Hom}_k(C,M)$ is a $k$-algebra via the convolution product
\begin{equation}\label{eq:conv.prod}
(f\diamond g)(c):= f(c_{(1)}) g(c_{(2)}), \qquad \textrm{for } f,g\in
\mathrm{Hom}_k(C,M), \ c\in C.
\end{equation}
A comodule algebra $M$ of a $k$-Hopf algebra $H$ is said to be a cleft
extension of its coinvariant subalgebra ${\bf N}$ provided that there exists
a convolution invertible map $j\in \mathrm{Hom}_k(H,M)$ which is an
$H$-comodule map. The relevance of cleft extensions by Hopf algebras stems
from Doi and Takeuchi's observation in \cite{DoiTak:cleft} that ${\bf
  N}\subseteq {\bf M}$ is a 
cleft extension if and only if it is a Galois extension and an additional
normal basis property holds, i.e. $M\cong N\sstac k H$ as a left ${\bf
N}$-module right $H$-comodule. What is more, (establishing an even stronger
similarity with Galois extensions of fields), ${\bf N}\subseteq {\bf M}$ is a
cleft extension if and only if ${\bf M}$ is isomorphic to a crossed product of
${\bf N}$ with $H$ with respect to an invertible 2-cocycle
\cite{DoiTak:cleft}, \cite{BCM:cleft}.

Above results have been extended to Hopf algebroids in
\cite{BohmBrz:hgdcleft}. In order to formulate the definition of a cleft
extension, as a first step, a generalised convolution product has to be
introduced. Using notations as in Definition \ref{def:hgd}, in a Hopf algebroid
$(\hH_L,\hH_R,S)$ there is an 
$\L$-coring $(H,\Delta_L,\epsilon_L)$ and an $\R$-coring 
$(H,\Delta_R,\epsilon_R)$ present. 
Consider an $\L\sstac k \R$-ring $M$, with multiplications $\mu_\L:M\sstac \L M
\to M$ and $\mu_\R: M\sstac \R M \to M$.
For these data, the convolution algebra \eqref{eq:conv.prod} 
can be generalised to a {\em convolution category}. It has two objects,
conveniently labelled by $\L$ and $\R$. For $P,Q\in \{\L,\R\}$, morphisms from
$P$ to $Q$ are $Q$-$P$ bimodule maps $H\to M$, where the bimodule structure of
the domain is determined by the ($P$- and $Q$-) coring structures of $H$ and
the bimodule structure of the codomain is determined by the ($P$- and $Q$-)
ring structures of $M$. For $P,Q,T\in \{\L,\R\}$, and 
morphisms $f:Q\to P$ and $g:T\to Q$,
composition is given by a convolution product 
\begin{equation}\label{eq:conv.cat}
f\diamond g := \mu_Q\circ (f\stac Q g)\circ \Delta_Q.
\end{equation}
Recall from Theorem \ref{thm:hgd.com.mon} that a right comodule algebra
of a 
Hopf algebroid $(\hH_L,\hH_R,S)$ has a canonical $\R$-ring structure over the
base algebra $\R$ of $\hH_R$. 
\begin{definition}
Consider a Hopf algebroid $\hH=(\hH_L,\hH_R,S)$, over base algebras $\L$ and
$\R$. A right $\hH$-comodule algebra $M$
is said to be a {\em cleft extension} of the $\hH_R$-coinvariant subalgebra
${\bf N}$ provided that the following properties hold. 
\begin{itemize}
\item[{(i)}] The canonical $\R$-ring structure of $M$ extends to an $\L\sstac
  k \R$-ring structure.
\item[{(ii)}] There exists an invertible morphism $j:\R\to \L$ in the
  convolution category \eqref{eq:conv.cat} which is a right $\hH$-comodule map.
\end{itemize}
\end{definition}
In an $\hH$-cleft extension ${\bf N}\subseteq {\bf M}$, ${\bf N}$ can be
proven to be an $\L$-subring of ${\bf M}$.

In a Hopf algebroid $\hH$, using the notations introduced in Definition
\ref{def:hgd}, an $\L\sstac k \R$-ring is given by $(\bH,s_L,s_R)$. The
identity map of $H$ is a morphism $\R\to \L$ in the corresponding convolution
category \eqref{eq:conv.cat}. It is obviously right $\hH$-colinear. What is
more, the antipode is 
its inverse by axioms (iii) and (iv) in Definition \ref{def:hgd}. Hence the
right regular 
comodule algebra of a Hopf algebroid is a cleft extension of the coinvariant 
subalgebra $t_R(\R)$. This extends a well known fact that the right regular
comodule algebra of a Hopf algebra, over a commutative ring $k$, is a cleft
extension of $k$. A further similarity between cleft extensions by Hopf
algebras and Hopf algebroids is expressed by the following theorem.
\begin{theorem}
Consider a Hopf algebroid $(\hH_L,\hH_R,S)$, over base algebras $\L$ and $\R$.
Its right comodule algebra
$M$ is a cleft extension of the $\hH_R$-coinvariant subalgebra ${\bf N}$ if
and only if the following properties hold.
\begin{itemize}\label{thm:cleft.norm.b}
\item[{(i)}] ${\bf N}\subseteq {\bf M}$ is a Galois extension by $\hH_R$;
\item[{(ii)}] the {\em normal basis condition} holds, i.e. $M\cong N\sstac \L
  H$ as left ${\bf N}$-modules right $\hH$-comodules.
\end{itemize}
\end{theorem}
Note the appearance of the two base algebras $\L$ and $\R$ in conditions (i)
and (ii) in Theorem \ref{thm:cleft.norm.b}.

Another characterisation of a cleft extension by a Hopf algebroid can be
given by using the construction of a crossed product.
\begin{definition}\label{def:measuring}
Consider a left bialgebroid $\cB$ over a base $k$-algebra $\L$. Denote its
structure maps as in Definition \ref{def:left.bgd}. 
We say that $\cB$ {\em measures} an $\L$-ring $N$ with unit map $\iota:\L\to
{\bf N}$ if there exists a $k$-module map $\cdot: B\otimes_k N\to N$, the so
called {\em measuring}, such that, for $b\in B$, $l\in L$ and $n,n'\in N$, the
following axioms are satisfied. 
\begin{itemize}
\item[{(i)}] $b\cdot 1_{\bf N}=\iota(\epsilon(b))$,
\item[{(ii)}] $(t(l)b)\cdot n =(b\cdot n)\iota(l)$ and $(s(l)b)\cdot n =
 \iota(l) (b\cdot n)$, 
\item[{(iii)}] $b\cdot (nn')=(b_{(1)}\cdot n)(b_{(2)}\cdot n')$.
\end{itemize}
\end{definition}
Note that in Definition \ref{def:measuring} condition (iii) makes sense in
view of (ii).

Consider a left bialgebroid $\cB$ over a $k$-algebra $\L$ and denote its
structure maps as in Definition \ref{def:left.bgd}. Let $N$ be an $\L$-ring,
with unit $\iota:\L\to {\bf N}$, which is measured by $\cB$. These data
determine a category ${\mathcal C}(\cB,N)$ as follows. Consider $B\sstac k B$
as an $\L$-bimodule, via left 
multiplication by $s$ and $t$ in the first factor. For an element $f$ in
${}_\L\mathrm{Hom}_\L(B\sstac k B,N)$, consider the following ($\L$-balancing)
conditions. For $a,b\in B$ and $l\in L$, 
\begin{itemize}
\item[{$(T\circ)$}]$\quad \qquad \qquad \qquad \qquad f(a\stac k
  s(l)b)=f(as(l)\stac   k b)$ 
\item[{$(S\circ)$}]$\quad \qquad \qquad \qquad \qquad f(a\stac k
  t(l)b)=f(at(l)\stac   k b)$ 
\item[{$(T\bullet)$}]$\quad  \qquad \qquad \qquad \qquad f(a\stac k
  s(l)b)=(a_{(1)}\cdot  \iota(l)) f(a_{(2)}\stac k b)$
\item[{$(S\bullet)$}]$\quad \qquad \qquad \qquad \qquad f(a\stac k
  t(l)b)=f(a_{(1)}\stac k  b) (a_{(2)}\cdot \iota(l))$.
\end{itemize}
Define a category ${\mathcal C}(\cB,N)$ of two objects $\circ$ and
$\bullet$. For two objects $X,Y\in \{\circ,\bullet\}$, morphisms $X\to Y$ are
elements of ${}_\L\mathrm{Hom}_\L(B\sstac k B,N)$, satisfying conditions
$(SX)$ and $(TY)$. Composition of morphisms $g:X\to Y$ and $f:Y\to Z$ is given
by 
\begin{equation*}
(f\diamond g)(a\stac k b):= f(a_{(1)}\stac k b_{(1)}) g(a_{(2)}\stac k
  b_{(2)}).  
\end{equation*}
Unit morphism at the object $\circ$ is the map $a\sstac k b \mapsto (ab)\cdot
1_{\bf N}=\iota(\epsilon(ab))$ and unit morphism at the object $\bullet$ is
the map $a\sstac k b \mapsto a\cdot(b\cdot 1_{\bf N})$.
\begin{definition}\label{def:cocycle}
Consider a left bialgebroid $\cB$ over a $k$-algebra $\L$ and denote its
structure maps as in Definition \ref{def:left.bgd}. Let $N$ be a
$\cB$-measured $\L$-ring, with unit $\iota:\L\to {\bf N}$.
An $N$-valued {\em 2-cocycle} on $\cB$ is a morphism $\circ\to \bullet$ in the
category ${\mathcal C}(\cB,N)$ above, such that, for $a,b,c\in B$,
\begin{itemize}
\item[{(i)}] $\sigma(1_\B,b)=\iota\big( \epsilon(b)\big)=\sigma(b,1_\B)$, 
\item[{(ii)}] $\big(a_{(1)}\cdot\sigma(b_{(1)},c_{(1)})\big)\sigma(a_{(2)},
  b_{(2)} c_{(2)}) = \sigma(a_{(1)},b_{(1)})\sigma(a_{(2)}b_{(2)}, c)$.
\end{itemize}
The $\cB$-measured $L$-ring $N$ is called a {\em $\sigma$-twisted
  $\cB$-module} if in addition, for $n\in N$ and $a,b\in B$,
\begin{itemize}
\item[{(iii)}] $1_\B \cdot n =n$,
\item[{(iv)}] $\big(a_{(1)}\cdot(b_{(1)}\cdot n)\big)\sigma(a_{(2)},b_{(2)})=
\sigma(a_{(1)},b_{(1)}) (a_{(2)}b_{(2)}\cdot n)$.
\end{itemize}
An $N$-valued 2-cocycle $\sigma$ on $\cB$ is said to be {\em invertible} if 
it is invertible as a morphism in ${\mathcal C}(\cB,N)$. 
\end{definition}
Note that in Definition \ref{def:cocycle} 
conditions (ii) and (iv) make sense in view of the module map and balanced
properties of $\sigma$.
If an $\L$-ring $N$ is measured by a left $\L$-bialgebroid $\cB$, then the
$\L^{op}$-ring $N^{op}$ is measured by the co-opposite left
$\L^{op}$-bialgebroid $\cB_{cop}$. The inverse of an $N$-valued 2-cocycle
$\sigma$ on $\cB$ turns out to be an $N^{op}$-valued 2-cocycle on $\cB_{cop}$. 

For a left bialgebroid $\cB$ over an algebra $\L$, with structure maps denoted
as in Definition \ref{def:left.bgd}, the base algebra $\L$ is measured by
$\cB$, via the left action $b\cdot l:=\epsilon(bs(l))$. For this measuring
conditions $(S\circ)$ and $(S\bullet)$ are equivalent and also conditions
$(T\circ)$ and $(T\bullet)$ are equivalent. Consequently, an
$\L$-valued 2-cocycle on $\cB$ in the sense of Definition \ref{def:cocycle} is
equivalent to a cocycle considered in Section \ref{sec:dual.twist}.
Extending cocycle double twists in Section \ref{sec:dual.twist}, one can
consider more general deformations of a Hopf algebroid $\hH$ 
(or a $\times_L$-Hopf algebra $\cB$, discussed in Section
\ref{sec:x_R.hopf.alg}) by an $N$-valued invertible 2-cocycle $\sigma$ in
Definition \ref{def:cocycle}, cf. 
\cite[Appendix]{BohmBrz:torsor}. In that construction the base algebra $\L$ of
$\hH_L$ is replaced by an $\hH_L$-measured $\L$-ring $N$. In particular, 
Connes and Moscovici's bialgebroids in Section \ref{sec:Con.Mos.ex} arise in
this way. 

A {\em crossed product} $N\#_\sigma B$
of a left $\L$-bialgebroid $\cB$ with a $\sigma$-twisted $\cB$-module $N$,
with respect to an $N$-valued 2-cocycle $\sigma$, is the the $\L$-module
tensor product $N\sstac \L B$ (where $B$ is a left $\L$-module via the source
map $s$), with associative and unital multiplication 
$$
(n\#b)(n'\#b')=n(b_{(1)}\cdot n')\sigma(b_{(2)},b'_{(1)})\#
b_{(3)}b'_{(2)},\qquad \textrm{for }n\#b, n'\#b'\in N\stac \L B.
$$
Equivalence classes of crossed products with a bialgebroid were classified in
\cite[Section 4]{BohmBrz:hgdcleft}.
\begin{theorem}
A right comodule algebra $M$ of a Hopf algebroid $(\hH_L,\hH_R,S)$ 
is a cleft extension of the $\hH_R$-coinvariant subalgebra ${\bf N}$ if and
only if 
$M$ is isomorphic, as a left ${\bf N}$-module and right 
$\hH$-comodule algebra, to a crossed product algebra $N\#_\sigma \hH_L$,
with respect to some invertible $N$-valued 2-cocycle $\sigma$ on $\hH_L$.
\end{theorem}

\subsubsection{The structure of Galois extensions by Hopf algebroids}
In the theory of Galois extensions by Hopf algebras (with a bijective
antipode) important tools are provided by theorems which state that in
appropriate situations {\em surjectivity} of the canonical map implies its
bijectivity, i.e. Galois property of an algebra extension ${\bf N}\subseteq
{\bf M}$. There are two big groups of such theorems. In the first group a 
Hopf algebra $H$ is assumed to be a flat module over its commutative base
ring $k$, and its regular comodule algebra 
is assumed to be a projective $H$-comodule. These properties hold in
particular if $H$ is a finitely generated and projective $k$-module, in which
case such a theorem was proven first by Kreimer and Takeuchi
\cite{KreTak:Gal}. In another group of such results, due to Schneider, $H$ is
assumed to be a projective $k$-module and its comodule algebra $M$ is assumed
to be a $k$-relative injective $H$-comodule \cite{Schn:Gal}. 

Analogous results for extensions by a Hopf algebroid $\hH=(\hH_L,\hH_R,S)$
were obtained in the papers \cite{Bohm:hgd.Gal} and \cite{BalSzl:fin.Gal}, and
in \cite{ArdBohmMen:Sch.type}, respectively. A common philosophy behind such
theorems originates from a work \cite{Scha:HbiGal} of Schauenburg (on the Hopf
algebra case). The key idea is to investigate a lifting \eqref{eq:T.lifted.can} 
of the canonical map \eqref{eq:two.canR}, introduced for an $\hH$-comodule
algebra $M$ below. By a general
result \cite[Theorem 2.1]{BrzTurWri} about Galois comodules, split
surjectivity of the lifted canonical map \eqref{eq:T.lifted.can}, as a
morphism of relative $(M,\hH_R)$-Hopf modules, implies $\hH_R$-Galois property,
i.e. bijectivity of \eqref{eq:two.canR}, whenever 
$(M\sstac {\bf T} M)^{co\hH_R}=M\sstac {\bf T} M^{co\hH_R}$, where
$(-)^{co\hH_R}$ denotes the $\hH_R$-coinvariants functor.

For a Hopf algebra $H$ over a commutative ring $k$, and a right $H$-comodule
algebra $M$ with coinvariant subalgebra ${\bf N}$, the canonical map $M\sstac
{\bf N} M \to M\sstac k H$ can be lifted to a map
\begin{equation}\label{eq:lifted.can}
\xymatrix{
M\stac k M \ar@{>>}[r]&M\stac {\bf N} M \ar[r]^{\mathrm{can}} &M\stac k H
}.
\end{equation}
More generally, consider a Hopf algebroid $\hH=(\hH_L,\hH_R,S)$ over base
$k$-algebras $\L$ and $\R$. For a right $\hH$-comodule algebra $M$, the
canonical map \eqref{eq:two.canR} can be lifted to 
\begin{equation}\label{eq:T.lifted.can}
M\stac {\bf T} M\to M\stac \R H,\qquad m\stac {\bf T} m'\mapsto m
m^{\prime[0]} \stac \R m^{\prime [1]},
\end{equation} 
for any $k$-algebra ${\bf T}$ such that the $\hH_R$-coinvariant subalgebra
${\bf N}$ of $M$ is a ${\bf T}$-ring.  
The map \eqref{eq:T.lifted.can} is a morphism of right-right
$(M,\hH)$-relative Hopf modules. 
\begin{theorem}\label{thm:KreTak}
Consider a Hopf algebroid $\hH$,
over base $k$-algebras $\L$ and $\R$, with a bijective antipode. Denote its
structure maps as in Definition \ref{def:hgd}. Assume that $H$ is a flat left
$\R$-module (via right multiplication by $t_R$) and a projective right
$\hH_R$-comodule (via $\Delta_R$). 
Let $M$ be a right $\hH$-comodule algebra with $\hH_R$-coinvariant subalgebra
${\bf N}$. Under these assumptions the following statements hold.
\begin{itemize}
\item[{(1)}]
If the $\hH_R$-coinvariants of the right $\hH$-comodule $M\sstac k M$
(with coactions given via the second factor) are precisely the elements of
$M\sstac k N$, 
then the canonical map \eqref{eq:two.canR} is bijective if and only if it is
surjective.
\item[{(2)}] If the canonical map \eqref{eq:two.canR} is bijective then
$M$ is a projective right ${\bf N}$-module. 
\end{itemize}
\end{theorem}
Since coinvariants are defined as a kernel, coinvariants of $M\sstac k M$ are
are precisely the elements of $M\sstac k N$ if e.g. $M$ is a flat $k$-module. 
In order to have an impression about the proof of part (1) of Theorem
\ref{thm:KreTak}, note that flatness of the left $\R$-module $H$ and
projectivity of the 
right regular $\hH_R$-comodule together imply that $M\sstac \R H$ is
projective as a right-right $(M,\hH_R)$-relative Hopf module. Hence if the the
canonical map \eqref{eq:two.canR} is surjective then the (surjective) 
lifted canonical map \eqref{eq:T.lifted.can} is a split epimorphism of
right-right $(M,\hH_R)$-relative Hopf modules, for any possible $k$-algebra
${\bf T}$. Thus bijectivity of the canonical map \eqref{eq:two.canR} follows
by \cite[Theorem 2.1]{BrzTurWri}. 
Part (2) of Theorem 
\ref{thm:KreTak} follows by exactness of the naturally equivalent 
functors $\mathrm{Hom}^{\hH_R}(H,-)\cong \mathrm{Hom}_{\bf
N}(M,(-)^{co\hH_R})$, which 
is a consequence of the projectivity of the right regular $\hH_R$-comodule.

If in a Hopf algebroid $\hH$ with a bijective antipode,
$H$ is a finitely generated and projective $\L$-, or $\R$-module in 
all of the four senses occurring in Proposition \ref{prop:fgp.hgd}, then it is 
obviously a flat left $\R$-module. Furthermore, in this case the right
regular $\hH_R$-comodule can be shown to be projective, by using Theorem
\ref{thm:fund.thm}. 
Following Weak Structure Theorem \ref{thm:hgd.mor} is thus  
based on  Remark \ref{rem:fgp.hgd}, part (2) of Theorem \ref{thm:KreTak} and
its application to the right comodule algebra $M^{op}$ of the opposite Hopf
algebroid $\hH^{op}$, and a theorem \cite[Theorem
  3.5]{CaeVerWang:Morita.cleft} about Galois corings.
\begin{theorem}\label{thm:hgd.mor}
Consider a Hopf algebroid $\hH=(\hH_L,\hH_R,S)$, over base algebras $\L$ and
$\R$, 
with a bijective antipode $S$. Assume that $H$ is a finitely generated and
projective $\L$-, or $\R$-module in all of the four senses occurring in
Proposition \ref{prop:fgp.hgd}. For a right $\hH$-comodule algebra $M$, with
$\hH_R$-coinvariant subalgebra ${\bf N}$, the following statements are
equivalent. 
\begin{itemize}
\item[{(i)}] The extension ${\bf N}\subseteq {\bf M}$ is $\hH_R$-Galois.
\item[{(ii)}] $M$ is a generator in the category $\cM^\hH_M\cong
  \cM^{\hH_R}_M\cong \cM^{\hH_L}_M$. 
\item[{(iii)}] The $\hH_R$-coinvariants functor $\cM^\hH_M\to \cM_{\bf N}$ is
  fully faithful. 
\item[{(iv)}] The extension ${\bf N}\subseteq {\bf M}$ is $\hH_L$-Galois.
\item[{(v)}] $M$ is a generator in the category ${}_M\cM^\hH\cong {}_M
  \cM^{\hH_R} \cong {}_M \cM^{\hH_L}$.  
\item[{(vi)}] The $\hH_L$-coinvariants functor ${}_M\cM^\hH\to {}_{\bf N}\cM$
  is fully faithful. 
\end{itemize}
Furthermore, if these equivalent conditions hold then $M$ is a projective left
and right ${\bf N}$-module.
\end{theorem}
It was a key observation by Doi that relative injectivity of a comodule
algebra $M$ of a Hopf algebra $H$ is equivalent to the existence of a so
called total integral -- meaning an $H$-comodule map $j:H\to M$,
such that $j(1_\bH)= 1_{\bf M}$ \cite{Doi:tot.int}. This fact extends to Hopf
algebroids as well.  
Recall (e.g. from \cite{GuRi:Kar}) that, for any functor ${\mathbb
  U}:{\mathcal A}\to {\mathcal B}$, 
between any categories ${\mathcal A}$ and ${\mathcal B}$, an object $A\in
{\mathcal A}$ is said to be {\em ${\mathbb U}$-injective}, if the map
$\Hom_{\mathcal A}(g,A): \Hom_{\mathcal A}(Y,A)\to \Hom_{\mathcal A}(X,A)$ is
surjective, for any objects $X,Y\in {\mathcal A}$, and all such morphisms
$g\in \Hom_{\mathcal A}(X,Y)$ for that ${\mathbb U}(g)$ is a split
monomorphism in ${\mathcal B}$. If ${\mathbb U}$ has a right adjoint, then
${\mathbb U}$-injectivity of an object $A$ is equivalent to the unit of the
adjunction, evaluated at $A$, being a split monomorphism in ${\mathcal A}$, see
\cite[Proposition 1]{GuRi:Kar}.
For example, for a Hopf algebra $H$ over a commutative ring $k$, injective
objects with respect to the forgetful functor $\cM^H \to \cM_k$ are precisely
relative injective $H$-comodules.  

A version of Theorem \ref{thm:doi.type} below was proven in \cite[Theorem
  4.1]{ArdBohmMen:Sch.type}, using the notion of {\em relative separability}
of a forgetful functor. 
Recall from Remark \ref{rem:antip.functor}
that the opposite of a right comodule algebra
$M$ of a Hopf algebroid $\hH$ has a canonical structure of a left
$\hH$-comodule algebra. 
\begin{theorem}\label{thm:doi.type}
For a right comodule algebra $M$ of a Hopf algebroid $\hH=(\hH_L,\hH_R,S)$,
the following statements are equivalent.
\begin{itemize}
\item[{(i)}] There exists a right $\hH$-comodule map (resp. right
  $\hH_L$-comodule map) $j:H \to M$, such that $j(1_\bH)=1_{\bf M}$.
\item[{(ii)}] $M$ is injective with respect to the forgetful functor
  $\cM^\hH\to \cM_\L$ (resp. with respect to the
  forgetful functor $\cM^{\hH_L}\to \cM_\L$, i.e. $M$ is a relative injective
  $\hH_L$-comodule). 
\item[{(iii)}] Any object in the category $\cM^\hH_M$ of right-right relative
  Hopf modules is injective with respect to the forgetful functor
  $\cM^\hH\to \cM_\L$ (resp. with respect to the
  forgetful functor $\cM^{\hH_L}\to \cM_\L$).
\end{itemize}
If in addition the antipode $S$ is bijective then assertions (i)-(iii) are
equivalent also to 
\begin{itemize}
\item[{(iv)}] There exists a left $\hH$-comodule map (resp. left
  $\hH_R$-comodule map) $j':H \to M$, such that $j'(1_\bH)=1_{\bf M}$.
\end{itemize}
Hence assertions (i)-(iv) are equivalent also to the 
symmetrical versions of (ii) and (iii).
\end{theorem}

The key idea behind Theorem \ref{thm:doi.type} is the observation that both
forgetful functors 
$\cM^\hH\to \cM_\L$ and $\cM^{\hH_L}\to \cM_\L$ possess left adjoints $-\ot_\L
H$ (cf. \eqref{eq:free.ff}). A correspondence can be established between 
comodule maps $j$ as in part (i), and natural retractions of the counit of the
adjunction, i.e. of the $\hH_L$-coaction. 

Based on \cite[Theorem 4.7]{Cae:Gal.coring&descent} and Theorem
\ref{thm:doi.type}, also the following Strong Structure Theorem holds.
\begin{theorem}
Consider a Hopf algebroid $(\hH_L,\hH_R,S)$, over base algebras $\L$ and $\R$, 
with a bijective antipode $S$. Assume that $H$ is a finitely generated and
projective $\L$-, or $\R$-module in any (hence all) of the senses listed in
Proposition \ref{prop:fgp.hgd}. For a right $\hH_R$-, (equivalently, right
$\hH_L$-) 
Galois extension ${\bf N}\subseteq {\bf M}$, the following statements are
equivalent. 
\begin{itemize}
\item[{(i)}] $M$ is a faithfully flat right ${\bf N}$-module.
\item[{(ii)}] The inclusion ${\bf N}\to {\bf M}$ splits in $\cM_{\bf N}$.
\item[{(ii)}] $M$ is a generator of right ${\bf N}$-modules.
\item[{(iii)}] $M$ is a faithfully flat left ${\bf N}$-module.
\item[{(iv)}] The inclusion ${\bf N}\to {\bf M}$ splits in ${}_{\bf N}\cM$.
\item[{(v)}] $M$ is a generator of left ${\bf N}$-modules.
\item[{(vi)}] The functor $-\sstac {\bf N} M:\cM_{\bf N}\to \cM^{\hH}_M$ is an
  equivalence. 
\item[{(vii)}] $M$ is a projective generator in $\cM^{\hH}_M$.
\item[{(viii)}] There exists a right $\hH$-comodule map $j:M\to H$, such that
  $j(1_\bH)=1_{\bf M}$. 
\end{itemize}
\end{theorem}
Note that if a Hopf algebra $H$, over a commutative ring $k$, is a projective
$k$-module, then $M\sstac k H$ is a projective left ${\bf M}$-module, for any
right $H$-comodule algebra $M$. Thus, denoting the subalgebra of coinvariants
in $M$ by ${\bf N}$, surjectivity of the canonical map $M\sstac {\bf N} M \to
M\sstac k H$ implies that its lifted version \eqref{eq:lifted.can}
is a split epimorphism of left ${\bf M}$-modules, so in particular of
$k$-modules. By Schneider's result \cite[Theorem I]{Schn:Gal}, if the antipode
of $H$ is bijective and 
$M$ is a $k$-relative injective right $H$-comodule, then bijectivity of the
canonical map follows from the $k$-module splitting of its lifted version 
\eqref{eq:lifted.can}. 
In order to formulate following generalisation Theorem \ref{thm:Sch.type} of
this result, note that the 
lifted canonical map \eqref{eq:T.lifted.can} is an $\L$-bimodule map, with
respect to  the $\L$-actions $l\cdot(m\sstac {\bf T} m')\cdot l':= 
m\cdot \epsilon_R(s_L(l))\sstac {\bf T} \epsilon_R(t_L(l'))\cdot m'$ 
(recall that the $\R$-, and ${\bf T}$-actions on $M$ commute by
virtue of \eqref{eq:n-r_comm})
and $l\cdot(m\sstac \R h)\cdot l':= m\sstac \R s_L(l) t_L(l')h$, on its domain 
and codomain, respectively, where notations in Definition \ref{def:hgd} are
used.  
\begin{theorem}\label{thm:Sch.type}
Consider a Hopf algebroid $\hH$ with a bijective antipode,
over base algebras $\L$ and $\R$. Denote its structure maps as in
Definition \ref{def:hgd}. Let $M$ be a right $\hH$-comodule algebra with
$\hH_R$-coinvariants ${\bf N}$. Let ${\bf T}$ be a $k$-algebra, such that
${\bf N}$ is a ${\bf T}$-ring. In this setting, if the
lifted canonical map \eqref{eq:T.lifted.can} is a split epimorphism of right
$\L$-modules then the following statements are equivalent.
\begin{itemize}
\item[{(i)}] ${\bf N}\subseteq {\bf M}$ is an $\hH_R$-Galois extension and 
the inclusion ${\bf N}\to {\bf M}$ splits in $\cM_{\bf N}$.
\item[{(ii)}] ${\bf N}\subseteq {\bf M}$ is an $\hH_R$-Galois extension and 
the inclusion ${\bf N}\to {\bf M}$ splits in ${}_{\bf N}\cM$.
\item[{(iii)}] There exists a right $\hH$-comodule map $j:H\to M$, such that
  $j(1_\bH)=1_{\bf M}$. 
\item[{(iv)}] $M\sstac {\bf N} -:{}_{\bf N}\cM\to {}^\hH_M\cM$ is an equivalence
  and the inclusion ${\bf N}\to {\bf M}$ splits in $\cM_{\bf N}$.
\end{itemize}
Furthermore, if the equivalent properties (i)-(iii) hold then $M$ is a ${\bf
T}$-relative projective right ${\bf N}$-module. 
\end{theorem}
Note that 
by Proposition \ref{prop:hgd.coinv} and Proposition
\ref{prop:bij.antip.Gal}, in parts (i) and (ii) $\hH_R$-Galois property can be
replaced equivalently by $\hH_L$-Galois property. Also, 
by Theorem \ref{thm:doi.type},
the existence of a unit preserving right $\hH$-comodule map in part (iii) can be
replaced equivalently by the existence of a unit preserving left
$\hH$-comodule map.  

The most interesting part of Theorem \ref{thm:Sch.type} is perhaps the claim
that if property (iii) holds then right $\L$-module splitting of the
lifted canonical map \eqref{eq:T.lifted.can} implies $\hH_R$-Galois property.
The proof of this fact is based on an observation, originated from
\cite{ArdBohmMen:Sch.type}, that assertion (iii) 
is equivalent to relative separability of the
forgetful functor $\cM^\hH \to \cM_\L$, with respect to the forgetful
functor $\cM^\hH_M \to \cM^\hH$. 
A relative separable functor reflects split epimorphisms in the 
sense that if $f$ is a morphism in $\cM^\hH_M$, which is a split epimorphism
of right $\L$-modules, then it is a split epimorphism of right
$\hH$-comodules. This proves that, under the assumptions made,
the lifted canonical map \eqref{eq:T.lifted.can} is a split epimorphism of
right $\hH$-comodules. Furthermore, the forgetful functor $\cM^\hH_M\to
\cM^\hH$ possesses a left adjoint $-\sstac \R M$. Hence the right-right
$(M,\hH)$-relative Hopf 
module $M\sstac \R H$, which is isomorphic to $H\sstac \R M$ by bijectivity of
the antipode, is relative projective in the sense that a split
epimorphism $g$ in $\cM^\hH$, of codomain $H\sstac \R M \cong M\sstac \R H$,
is a split epimorphism in $\cM^\hH_M$. This proves that in the situation
considered the lifted canonical
map \eqref{eq:T.lifted.can} is a split epimorphism of right-right
$(M,\hH)$-relative Hopf modules. 
Then it is a split epimorphism of right-right $(M,\hH_R)$-relative Hopf
modules. 
Moreover, in terms of a unit preserving right $\hH$-comodule map $H\to M$, a
left ${\bf N}$-module splitting of the equaliser of the $\hH_R$-coaction on
$M$, and the map $m\mapsto m\ot_\R 1_\bH$,
can be constructed. This implies that the equaliser is preserved by the
functor $M\ot_{\bf T} -$, i.e. $(M\ot_{\bf T} M)^{co\hH_R}= M \ot_{\bf T} N$.
Thus the canonical map \eqref{eq:two.canR} is bijective by
\cite[Theorem 2.1]{BrzTurWri}. 

Recall that a left module $V$ of a $k$-algebra ${\bf N}$ is $k$-relative
projective if and only if the left action $N\sstac k V \to V$ is a split
epimorphism of left ${\bf N}$-modules. If $V$ has an additional structure of a
right comodule for a $k$-coalgebra $C$, such that the ${\bf N}$-action is a
right $C$-comodule map, then it can be asked if the action $N\sstac k V \to V$
splits as a map of left ${\bf N}$-modules and right $C$-comodules too. In the
case when it does, $V$ is said to be a $C$-equivariantly projective left ${\bf
  N}$-module. For a Galois extension 
${\bf N}\subseteq {\bf M}$ by a $k$-Hopf algebra $H$, $H$-equivariant
projectivity of the left ${\bf N}$-module $M$ was shown by Hajac to be
equivalent to the existence of a strong connection \cite{Hajac}. Interpreting
a Hopf Galois 
extension as a non-commutative principal bundle, this means its local
triviality. In case of a Galois extension ${\bf N}\subseteq {\bf M}$ by a
$k$-Hopf algebra $H$, equivalent conditions (i)-(iii) in 
Theorem \ref{thm:Sch.type} are known to imply $H$-equivariant projectivity of
the left ${\bf N}$-module $M$. In order to obtain an analogous result for a
Galois extension by a Hopf algebroid, slightly stronger assumptions are needed,
see Theorem \ref{thm:eq.proj} below.
\begin{definition}
Consider a Hopf algebroid $\hH$ and consider a ${\bf T}$-ring $N$ for some
algebra ${\bf T}$. 
Let $V$ be a left ${\bf N}$-module and right $\hH$-comodule, such that the
left ${\bf N}$-action on $V$ is a right ${\hH}$-comodule map. $V$ is said to
be a {\em ${\bf T}$-relative $\hH$-equivariantly projective} left ${\bf
N}$-module provided that the left action $N\sstac 
{\bf T} V \to V$ is an epimorphism split by a left ${\bf N}$-module, right
$\hH$-comodule map. 
\end{definition}
\begin{theorem}\label{thm:eq.proj}
Consider a Hopf algebroid $\hH$ with a bijective antipode.
Let $M$ be a right $\hH$-comodule algebra with $\hH_R$-coinvariants ${\bf N}$.  
Let ${\bf T}$ be an algebra, such that ${\bf N}$ is a ${\bf T}$-ring. 
Assume that there exists a unit preserving right $\hH$-comodule map $H\to M$ 
and that the
lifted canonical map \eqref{eq:T.lifted.can} is a split epimorphism of
$\L$-bimodules. Then ${\bf N}\subseteq {\bf M}$ is a right $\hH_R$-, and
right $\hH_L$-Galois extension and $M$ is a ${\bf T}$-relative
$\hH$-equivariantly projective left ${\bf N}$-module.  
\end{theorem}
Under the premises of Theorem \ref{thm:eq.proj}, $\hH_R$-Galois property holds
by Theorem \ref{thm:Sch.type} and the $\hH_L$-Galois property follows by
Proposition \ref{prop:hgd.coinv} and Proposition \ref{prop:bij.antip.Gal}.
Proof of Theorem \ref{thm:eq.proj} is completed by constructing a required
left ${\bf N}$-module right $\hH$-comodule section of the action $N\sstac {\bf
T} M\to M$. The construction makes use of the relative Hopf module section of
the lifted canonical map \eqref{eq:T.lifted.can}, on the existence of which it
is concluded in the paragraph following Theorem \ref{thm:Sch.type}. Examples
of ${\bf L}$-relative $\hH$-equivariantly projective Galois extensions are
provided by cleft extensions by a Hopf algebroid $\hH$ 
with a bijective antipode, over base algebras $\L$ and $\R$. 

\subsection{Alternative notions}
In the literature there is an accord that the right generalisation of a
bialgebra to the case of a non-commutative base algebra is a bialgebroid. On
the contrary, there is some discussion about the structure to replace a Hopf
algebra. In current final section we revisit and compare the various
suggestions.  

\subsubsection{Lu's Hopf algebroid}
In Definition \ref{def:hgd} the antipode axioms are formulated for a compatible
pair of a left and a right bialgebroid. In following Definition
\ref{def:Lu.hgd}, quoted from \cite{Lu:hgd}, only a left bialgebroid is
used. While the first one of the antipode axioms in Definition \ref{def:hgd}
(iv) is easily formulated also in this case, in order to formulate 
the second one some additional assumption is needed.
\begin{definition}\label{def:Lu.hgd}
Consider a left bialgebroid $\cB$, over a $k$-algebra $\L$, with structure maps
denoted as in Definition \ref{def:left.bgd}. $\cB$ is a {\em Lu's Hopf
  algebroid} provided that there exists an anti-algebra map $S:\B\to \B$, and a
$k$-module section $\xi$ of the canonical epimorphism $B\sstac k B \to B\sstac
\L B$, such that the following axioms are satisfied.
\begin{itemize}
\item[{(i)}] $S\circ t =s$,
\item[{(ii)}] $\mu_B\circ(S\sstac \L B)\circ \Delta = t\circ \epsilon \circ
  S$,
\item[{(iii)}] $\mu_\B\circ (B\sstac k S)\circ \xi\circ \Delta = s\circ
  \epsilon$,
\end{itemize}
where $\mu_B$ denotes multiplication in the $\L$-ring $(B,s)$ and $\mu_\B$
denotes multiplication in the underlying $k$-algebra $\B$.
\end{definition}
None of the notions of a Hopf algebroid in Definition \ref{def:hgd}, or in
Definition \ref{def:Lu.hgd}, seems to be more general than the other
one. Indeed, a 
Hopf algebroid in the sense of Definition \ref{def:hgd}, which does not
satisfy the axioms in Definition \ref{def:Lu.hgd}, is constructed as
follows. Let $k$ be a commutative ring in which $2$ is invertible. For
the order 2 cyclic group $Z_2$, consider the group bialgebra $kZ_2$ as a left
bialgebroid over $k$. Equip it with the twisted (bijective) antipode $S$,
mapping the order 2 generator $t$ of $Z_2$ to $S(t):=-t$. Together with the
unique right bialgebroid, determined by the requirement that $S$ is a
bialgebroid anti-isomorphism in the sense of Proposition \ref{prop:antipode},
they constitute a Hopf algebroid as in Definition \ref{def:hgd}. However, for
this Hopf algebroid there exists no section $\xi$ as in Definition
\ref{def:Lu.hgd}. 

\subsubsection{A $\times_R$-Hopf algebra}\label{sec:x_R.hopf.alg}
The coinvariants of the (left or right) regular comodule of a bialgebra $H$,
over a commutative ring $k$, are precisely the multiples of the unit element
$1_\bH$. $H$ is known to be a Hopf algebra if and only if $\bH$ is an
$H$-Galois extension of $k$. Indeed, the hom-tensor relation
${}_H\mathrm{Hom}^H(H\sstac k H,H\sstac k H)\cong \mathrm{Hom}_k(H,H)$ relates
the inverse of the canonical map to the antipode. Motivated by this
characterisation of a Hopf algebra, in \cite{Scha:dual.double.bgd} Schauenburg
proposed the following definition.
\begin{definition} \label{def:x_L.ha}
Let $\cB$ be a left bialgebroid over an algebra $\L$, with structure maps 
denoted as in Definition \ref{def:left.bgd}. Consider the left regular
$\cB$-comodule, whose coinvariant subalgebra is $t(\L^{op})$. $\cB$ is said to
be a {\em $\times_L$-Hopf algebra} provided that the algebra extension
$t:\L^{op}\to \B$ is left $\cB$-Galois. 
\end{definition}
The notion of a $\times_L$-Hopf algebra in Definition \ref{def:x_L.ha} is more
general than that of a Hopf algebroid in Definition \ref{def:hgd}. Indeed,
consider a Hopf algebroid $(\hH_L,\hH_R,S)$, over the base algebras $\L$ and
$\R$, with structure maps denoted as in Definition \ref{def:hgd}. The
canonical map 
$H\sstac {\L^{op}} H \to H\sstac \L H$, $h\otimes h'\mapsto h_{(1)}\otimes
h_{(2)}h'$ is bijective, with inverse $h\otimes h'\mapsto h^{(1)}\otimes
S(h^{(2)})h'$. Hence $\hH_L$ is a $\times_L$-Hopf algebra. Although it is
believed that not every a $\times_L$-Hopf algebra is a constituent left
bialgebroid in a Hopf algebroid, we do not know about any counterexample.

Extending a result \cite{Sch:ff} by Schauenburg about Hopf algebras, the
following proposition was proven in \cite {Hobst:phd}.
\begin{proposition}\label{prop:Hobst}
Consider a left bialgebroid $\cB$ over a base algebra $\L$. If there is a 
left $\cB$-Galois extension ${\bf N}\subseteq {\bf M}$, such that $M$ is a
faithfully flat left $\L$-module, then $\cB$ is a $\times_L$-Hopf algebra.
\end{proposition}
Indeed, the canonical maps $\mathrm{can}:M\sstac {\bf N} M\to B\sstac \L M$
and $\vartheta:B\sstac {\L^{op}} B \to B \sstac \L B$ satisfy the pentagonal
identity 
$
(B\sstac \L \mathrm{can})\circ (\mathrm{can} \sstac {\bf N} M)=
(\vartheta\sstac \L M)\circ \mathrm{can}_{13}\circ (M\sstac {\bf N}
\mathrm{can}),
$
where the (well defined) map $\mathrm{can}_{13}: M \sstac {\bf N} (B\sstac \L
M)\to (B\sstac {\L^{op}} B) \sstac \L M$ is obtained by applying $\mathrm{can}$
in the first and third factors.

In Theorem \ref{thm:bgd.mod.cat} bialgebroids were characterised via strict
monoidality of a forgetful functor. A characterisation of a similar
flavour of $\times_L$-Hopf algebras 
was given in \cite[Theorem and Definition 3.5]{Scha:dual.double.bgd}. Recall
that a 
monoidal category $(\cM,\otimes, U)$ is said to be {\em right closed} if the
endofunctor $-\otimes X$ on $\cM$ possesses a right adjoint, denoted by
$\mathrm{hom}(X,-)$, for any object $X$ in $\cM$. The monoidal category of
bimodules of an algebra $\L$ is right closed with
$\mathrm{hom}(X,Y)=\mathrm{Hom}_\L (X,Y)$. It is slightly more involved to see
that so 
is the category of left $\B$-modules, for a left $\L$-bialgebroid $\cB$, with
$\mathrm{hom}(X,Y) ={}_\B\mathrm{Hom}(B\sstac \L X, Y)$ (where, for a left
$\B$-module $X$, $B\sstac \L X$
is a left $\B$-module via the diagonal action and a right $\B$-module via the
first factor). A strict monoidal functor
${\mathbb F}:\cM\to \cM'$ between right closed categories is called {\em
strong right closed} provided that a canonical morphism ${\mathbb
F}(\mathrm{hom}(X,Y))\to \mathrm{hom}({\mathbb F}(X),{\mathbb F}(Y))$ is an
isomorphism, for all objects $X,Y\in \cM$. 
\begin{theorem}\label{thm:str.closed}
A left bialgebroid $\cB$ over a base algebra $\L$ is a $\times_L$-Hopf algebra
if and only if the (strict monoidal) forgetful functor ${}_\B\cM\to {}_\L
\cM_\L$ is strong right closed. 
\end{theorem}
For a Hopf algebra $H$ over a commutative ring $k$, those (left or right)
$H$-modules, which are finitely generated and projective $k$-modules, possess
(right or left) duals in the monoidal category of (left or right)
$H$-modules. This property extends to $\times_L$-Hopf algebras (hence to Hopf
algebroids!) as follows.
\begin{proposition}\label{prop:dual.mod.x_Lha}
Let $\cB$ be a $\times_\L$-Hopf algebra over a base algebra $\L$. Denote its
structure maps as in Definition \ref{def:left.bgd}. For a left $\B$-module
$M$, the dual $M^*:=\mathrm{Hom}_\L(M,L)$ is a left $\B$-module, via the
action 
$$
(b\cdot \phi)(m):=\epsilon\big(b^{\langle 1\rangle} t(\phi(
b^{\langle 2\rangle}\cdot m))\big),\qquad \textrm{for }b\in B,\ \phi\in M^*,\
  m\in M,
$$
where for the inverse of the (right $\B$-linear) canonical map $B\sstac
{\L^{op}} B \to B\sstac \L B$ the index notation $b\otimes b'\mapsto
b^{\langle 1\rangle}\otimes b^{\langle 2\rangle} b'$ is used. Furthermore, if
$M$ is a finitely generated and projective right $\L$-module then $M^*$ is a
right dual of $M$ in the monoidal category ${}_\B\cM$.
\end{proposition}

\subsubsection{Hopf monad}\label{sec:Hopf.monad}
In Theorem \ref{thm:bimonad} bialgebroids were related to those bimonads on a
bimodule category which possess a right adjoint. In the paper
\cite{BrugVire:Hopf.monad} special bimonads, so called Hopf monads, on
autonomous categories were studied. Recall that a monoidal category is said to
be left (resp. right) autonomous provided that every object possesses a left
(resp. right) dual. In particular, a category of finitely generated and
projective bimodules is autonomous. Instead of a somewhat technical definition
in \cite[3.3]{BrugVire:Hopf.monad}, we adopt an equivalent description in 
\cite[Theorem 3.8]{BrugVire:Hopf.monad} as a definition.
\begin{definition}
A {\em left (resp. right) Hopf monad} is a bimonad $\mathbb{B}$ on a left
(resp. right) autonomous monoidal category $\cM$, such that the 
left (resp. right) autonomous structure of $\cM$ lifts to the category of
$\mathbb{B}$-algebras.
\end{definition}
The reader should be warned that, although the same term `Hopf monad' is used
in the papers \cite{BrugVire:Hopf.monad} and \cite{Moe:monad.tens.cat}, they
have different meanings (and a further totally different meaning of the same
term is used in \cite{MesWis:Hopf.mon}). Also, the notions of a comodule and a
corresponding (co)integral in  \cite{BrugVire:Hopf.monad} are different from
the notions used in these notes. 

\subsubsection{A $^*$-autonomous structure on a strong monoidal special
  opmorphism between pseudomonoids in a monoidal bicategory}
  \label{sec:Day.Street} 
In the paper \cite{DayStr:Qcat}, strong monoidal special opmorphisms $h$ in
monoidal bicategories, from a canonical pseudomonoid $R^{op}\otimes R$ to some
pseudomonoid $B$, were studied. The opmorphism $h$ was called {\em
Hopf} if in addition there is a $^*$-autonomous structure on 
$B$ and $h$ is strong $^*$-autonomous (where $R^{op}\otimes R$ is meant to
be $^*$-autonomous in a canonical way). In \cite[Section 3]{DayStr:Qcat} a
bialgebroid was described as a strong monoidal special opmorphism $h$ of
pseudomonoids in the monoidal bicategory of {\tt [Algebras; Bimodules;
    Bimodule maps]}. This opmorphism $h$ is strong
$^*$-autonomous if and only if the corresponding bialgebroid constitutes a
Hopf algebroid with a bijective antipode, see \cite[Section
  4.2]{BohmSzl:hgdax}. 

\section*{Acknowledgement}
The author's work is supported by the Bolyai J\'anos Fellowship and the
Hungarian Scientific Research Fund OTKA F67910. She is
grateful for their helpful comments to D\'enes Bajnok, Imre B\'alint, Tomasz
Brzezi\'nski and Korn\'el Szlach\'anyi.

\end{document}